\documentclass[10pt]{amsart}

\usepackage{xypic,amscd,amssymb,latexsym}
\usepackage{amsmath}	
\usepackage{graphicx}
\usepackage{subcaption}

\usepackage{pdfsync}

\usepackage{hyperref}

\usepackage{cancel}

\usepackage{mathrsfs}

\usepackage{soul}

\usepackage{url}
\usepackage{pifont}

\DeclareMathAlphabet{\mathpzc}{OT1}{pzc}{m}{it}

\usepackage{setspace}

\usepackage{changepage}

\usepackage{xcolor}

\renewcommand{\thesubfigure}{\arabic{subfigure}}

\begin{document}
\pagenumbering{arabic}
\newtheorem{theorem}{Theorem}[section]
\newtheorem{proposition}[theorem]{Proposition}
\newtheorem{lemma}[theorem]{Lemma}
\newtheorem{corollary}[theorem]{Corollary}
\newtheorem{remark}[theorem]{Remark}
\newtheorem{definition}[theorem]{Definition}
\newtheorem{question}[theorem]{Question}
\newtheorem{claim}[theorem]{Claim}
\newtheorem{conjecture}[theorem]{Conjecture}
\newtheorem{defprop}[theorem]{Definition and Proposition}
\newtheorem{example}[theorem]{Example}
\newtheorem{deflem}[theorem]{Definition and Lemma}

\newcommand{\II}[0]{{\rm \it II}}

\definecolor{delete}{rgb}{0.9, 0.6, 0.9}

\def\qed{{\quad \vrule height 8pt width 8pt depth 0pt}}

\newcommand{\cplx}[0]{\mathbb{C}}

\newcommand{\fr}[1]{\mathfrak{#1}}

\newcommand{\vs}[0]{\vspace{2mm}}

\newcommand{\til}[1]{\widetilde{#1}}

\newcommand{\mcal}[1]{\mathcal{#1}}

\renewcommand{\ul}[1]{\underline{#1}}

\newcommand{\ol}[1]{\overline{#1}}

\newcommand{\wh}[1]{\widehat{#1}}

\renewcommand\thesubfigure{\roman{subfigure}}


\author[H. Kim]{Hyun Kyu Kim}
\email{hkim@kias.re.kr}

\author[T. T. Q. L\^{e}]{Thang T. Q. L\^{e}}
\email{letu@math.gatech.edu}

\author[M. Son]{Miri Son}
\email{ms235@rice.edu}

\address{School of Mathematics, Korea Institute for Advanced Study (KIAS), 85 Hoegiro Dongdaemun-gu, Seoul 02455, Republic of Korea}

\address{School of Mathematics, 686 Cherry Street, Georgia Tech, Atlanta, GA 30332, USA}

\address{Department of Mathematics, Rice University, 6100 Main Street, Houston, TX 77005, USA}

\thanks{This paper appears in: Algebr. Geom. Topol. {\bf 23}(1) (2023), 339-418. \quad DOI: 10.2140/agt.2023.23.339}

\numberwithin{equation}{section}

\title[${\rm SL}_2$ quantum trace in quantum Teichm\"uller theory via writhe]{${\rm SL}_2$ quantum trace in quantum Teichm\"uller theory via writhe}

\begin{abstract}
Quantization of the Teichm\"uller space of a punctured Riemann surface $S$ is an approach to $3$-dimensional quantum gravity, and is a prototypical example of quantization of cluster varieties. Any simple loop $\gamma$ in $S$ gives rise to a natural trace-of-monodromy function $\mathbb{I}(\gamma)$ on the Teichm\"uller space. For any ideal triangulation $\Delta$ of $S$, this function $\mathbb{I}(\gamma)$ is a Laurent polynomial in the square-roots of the exponentiated shear coordinates for the arcs of $\Delta$. An important problem was to construct a quantization of this function $\mathbb{I}(\gamma)$, namely to replace it by a noncommutative Laurent polynomial in the quantum variables. This problem, which is closely related to the framed protected spin characters in physics, has been solved by Allegretti and Kim using Bonahon and Wong's ${\rm SL}_2$ quantum trace for skein algebras, and by Gabella using Gaiotto, Moore and Neitzke's Seiberg-Witten curves, spectral networks, and writhe of links. We show that these two solutions to the quantization problem coincide. We enhance Gabella's solution and show that it is a twist of the Bonahon-Wong quantum trace. 
\end{abstract}

\maketitle

\vspace{-9mm}

\tableofcontents

\vspace{-9mm}

\section{Introduction}

\subsection{Quantization of Teichm\"uller space}

Quantization of the Teichm\"uller space of a Riemann surface appeared in the late 1990s as an approach to $(2+1)$-quantum gravity in mathematical physics. Certain versions of quantum Teichm\"uller spaces were first established by Kashaev \cite{Kash} and independently by Chekhov and Fock \cite{CF} \cite{F}. The Chekhov-Fock formulation was later generalized by Fock and Goncharov to the theory of quantum cluster varieties \cite{FG09}. Here we review the Chekhov-Fock(-Goncharov) quantum Teichm\"uller theory.

\vs

\def\fS{{\mathfrak S}}
\newcommand{\red}[1]{{\color{red}#1}}
\newcommand{\blue}[1]{[[{\small \color{blue}#1}]]}

Let $\fS$ be a {\em punctured surface}, which is the result of removing a finite number of points, called {\em punctures}, from a closed oriented surface of finite genus. Among several versions of Teichm\"uller spaces of $\fS$, we will consider the {\em enhanced Teichm\"uller space} $\mathcal{X}^+(\frak{S})$, also known as the {\em holed Teichm\"uller space}; see Fock and Goncharov \cite{F, FG07} as well as Bonahon and Liu \cite{BL, Liu}. By definition, $\mathcal{X}^+(\frak{S})$ is the set of all complete hyperbolic metrics $h$ on $\fS$ up to isotopy, enhanced with a choice of an orientation of a small loop surrounding each puncture whose monodromy with respect to $h$ is hyperbolic. The space $\mathcal{X}^+(\frak{S})$ is a contractible smooth manifold equipped with the {\em Weil-Petersson} Poisson structure \cite{F, FG07}. A nice coordinate system requires the choice of an {\em ideal triangulation} $\Delta$ of $\frak{S}$, which is a collection of isotopy classes of mutually nonintersecting simple paths running between punctures, called {\em (ideal) arcs}, dividing $\frak{S}$ into {\em (ideal) triangles}. Per each arc $e$ of $\Delta$, Thurston's {\em shear coordinate} function $x_e = x_{e;\Delta}$ is associated \cite{BL, F, FG07, Liu, P, Thurston}, and they together comprise a global coordinate system of $\mathcal{X}^+(\frak{S})$. Their exponentials $X_e := \exp(x_e)$ are particularly convenient. The Poisson bracket is given by $\{X_e,X_f\} = \varepsilon_{ef} X_e X_f$, for all $e,f\in \Delta$, where $\varepsilon_{ef} = a_{ef} - a_{fe} \in \mathbb{Z}$, with $a_{ef}$ being the number of corners of ideal triangles of $\Delta$ formed by $e$ and $f$ appearing clockwise in this order. Per change of ideal triangulations, the coordinate change formulas for the exponentiated shear coordinates are rational, providing examples of so-called cluster $\mathcal{X}$-mutations; see Fock and Goncharov \cite{FG06, FG07}. 

\vs

For a chosen triangulation $\Delta$, a quantum deformed version of the classical Poisson algebra generated by the exponentiated shear coordinates $X_e$ and their inverses had been constructed as the $\mathbb{C}$-algebra generated by the symbols $\wh{X}_e$ and their inverses, with $e$ running through all arcs of $\Delta$, modulo the relations $\wh{X}_e \wh{X}_f = q^{2\varepsilon_{ef}} \wh{X}_f \wh{X}_e$ for all $e,f\in \Delta$, where $q \in \mathbb{C}^*$ is the quantum parameter. This algebra is called the {\em Chekhov-Fock algebra} $\mathcal{T}^q_\Delta$; see Liu \cite{Liu} and Bonahon and Liu \cite{BL}.  A key aspect of quantum Teichm\"uller theory is consistency under change of triangulations. It has been shown that, for each pair of ideal triangulations $\Delta$ and $\Delta'$, there exists a quantum coordinate change isomorphism $\Phi^q_{\Delta,\Delta'} : {\rm Frac}(\mathcal{T}^q_{\Delta'}) \to {\rm Frac}(\mathcal{T}^q_\Delta)$ between the skew-fields of fractions of the Chekhov-Fock algebras that recovers the classical formula as $q=1$ and satisfies the consistency equations $\Phi^q_{\Delta,\Delta''} = \Phi^q_{\Delta,\Delta'} \circ \Phi^q_{\Delta',\Delta''}$ \cite{CF, Liu}. The algebras ${\rm Frac}(\mathcal{T}^q_\Delta)$ for different $\Delta$ can then be identified by these isomorphisms in a consistent manner, forming a {\em quantum Teichm\"uller space} \cite{BL, Liu}.

\vs

The algebraic aspects of \cite{BL, CF, F, Kash, FG09, Liu} were mostly on establishing this consistent quantum Teichm\"uller space as just described, but not so much on finding a {\em deformation quantization map}, a map sending a classical function to a quantum counterpart. One must first determine which functions to quantize. Fock and Goncharov \cite{FG06} suggested to consider the functions that can be written as Laurent polynomials in the exponentiated shear coordinates $(X_e)_{e\in \Delta}$ for {\em every} ideal triangulation $\Delta$, i.e. the regular functions on $\mathcal{X}^+(\frak{S})$ viewed as ($\mathbb{R}_{>0}$-points of) a Fock-Goncharov cluster $\mathcal{X}$-variety. This important class of functions, called the {\em universally Laurent functions}, forms the classical Poisson algebra to be quantized, and is denoted by ${\bf L}(\mathcal{X}^+(\frak{S}))$. By a deformation quantization we mean a map
\begin{equation}
\mathbb{I}^q_\Delta : {\bf L}(\mathcal{X}^+(\frak{S})) \to \mathcal{T}^q_\Delta \label{eq.qu}
\end{equation}
satisfying favorable properties, e.g. when $q=1$ this should be the tautological classical map, and this map (for general $q$) must be consistent under the quantum coordinate change $\mathbb{I}^q_{\Delta'} = \mathbb{I}^q_\Delta \circ \Phi^q_{\Delta,\Delta'}$. Existence of such a map was conjectured in \cite{FG06} and a solution was first given by Allegretti and Kim \cite{AK}. Later, another solution was given by Gabella \cite{G}. Our goal is to compare these two solutions.

\subsection{Allegretti-Kim quantization using Bonahon-Wong quantum trace} A basic example of a universally Laurent function on the enhanced Teichm\"uller space $\mathcal{X}^+(\frak{S})$ is the trace-of-monodromy function along a loop. Let $\gamma$ be any non-null-homotopic simple loop in $\frak{S}$. Note that a point of the Teichm\"uller space yields a monodromy representation $\rho : \pi_1(\frak{S}) \to {\rm PSL}(2,\mathbb{R})$, defined up to conjugation in ${\rm PSL}(2,\mathbb{R})$. The function $\mathbb{I}(\gamma) := |{\rm trace}(\rho(\gamma))|$ is almost universally Laurent on $\mathcal{X}^+(\frak{S})$; it is universally Laurent not in the exponentiated shear coordinates $X_e$ in general, but only in their square-roots $Z_e := \sqrt{X_e}$. It is known (see Fock and Goncharov \cite{FG06}, Shen \cite{Shen} and Gross, Hacking and Keel \cite{GHK15}) that any element of ${\bf L}(\mathcal{X}^+(\frak{S}))$ is generated by these basic functions, together with relatively easy extra functions associated to punctures. So the problem of quantization almost boils down to quantizing these trace-of-monodromy functions $\mathbb{I}(\gamma)$. 

\def\hZ{\widehat Z}
\def\cT{\mathcal T}
\def\BC{{\mathbb C}}
\def\SS{\mathcal{S}^A (\frak{S})  }
\def\TT{\cT^\omega_\Delta}
\def\be{ \begin{equation}   }
\def\ee{ \end{equation}   }
\def\Tr{\mathrm{Tr}}
\def\ZZ{\mathcal{Z}^\omega_\Delta}

\vs

The appearance of the square-roots of the exponentiated shear coordinates in $\mathbb{I}(\gamma)$ forces us to replace the target space of the map in \eqref{eq.qu} by the bigger space $\cT^\omega_\Delta$, which is the $\BC$-algebra generated by the symbols $\hZ_e$ and their inverses, modulo the relations $\hZ_e \hZ_f = \omega^{2 \varepsilon_{ef}} \hZ_f \hZ_e$. Here $q= \omega^4$, and $\mathcal{T}^q_\Delta$ embeds into $\mathcal{T}^\omega_\Delta$ by $\wh{X}_e \mapsto \wh{Z}_e^2$, for all $e\in T$.

\vs

Bonahon-Wong {\em quantum trace} \cite{BW} (see also \cite{CL,Le, Le18}) is an algebra homomorphism 
\be  \Tr_\Delta^\omega: \SS \to \TT, \qquad \text{with} \ A= \omega^{-2}.
\label{eq.qt}
\ee
Here $\SS$ is the {\em skein algebra} (see Przytycki \cite{Przy} and Turaev \cite{Turaev91}) of $\fS$ at quantum parameter $A \in \mathbb{C}^*$, which is the $\BC$-module generated by isotopy classes of framed unoriented links in the thickened surface $\fS \times (-1,1)$ subject to certain local relations. The product of two framed links is given by stacking one above the other. The skein algebra $\SS$ is known to be a quantization of the ${\rm SL}_2$-character variety of $\fS$ in the direction of the Weil-Petersson-Goldman Poisson bracket; see Przytycki and Sikora \cite{PS}. Actually the image of the quantum trace map \eqref{eq.qt} is in a smaller subalgebra $\ZZ$ of $\TT$, and we can extend the quantum coordinate change isomorphism of the Chekhov-Fock algebras to the skew-field of fractions of $\ZZ$; see Hiatt \cite{Hiatt} and Bonahon and Wong \cite{BW}. The quantum trace map is compatible with this extended quantum coordinate change. We will recall the definitions of 
 the skein algebra, the Bonahon-Wong quantum trace map, as well as all these facts in Section \ref{sec:BW}.

\vs

Any simple loop $\gamma$ in $\fS$ can be considered as a framed link in the thickened surface $\fS \times (-1,1)$, using the standard blackboard framing. 
In Allegretti and Kim \cite{AK}, the function $\mathbb{I}(\gamma)$ is quantized by setting 
$$\mathbb{I}_\Delta^q(\gamma)= \Tr^\omega_\Delta(\gamma).$$
We note that a full description of $\mathbb{I}^q_\Delta$ involves still more non-trivial ideas; see \cite{AK}. For example, the above recipe $\mathbb{I}_\Delta^q(\gamma)= \Tr^\omega_\Delta(\gamma)$ does not immediately fit into the original setting as in eq.\eqref{eq.qu}, but only into a certain square-root version. We also note that the Bonahon-Wong quantum trace enjoys some more properties in accordance with the conjecture of Fock-Goncharov \cite{FG06} for $\mathbb{I}^q_\Delta$; see \cite{AK}, and also Cho, Kim, Kim and Oh \cite{CKKO} for an important positivity property.

\def\GH{{\rm TrHol}^\omega_\Delta(\gamma)}
\subsection{Another quantization using Gabella quantum holonomy} 
Gabella \cite{G} constructed, for every oriented simple loop $\gamma$ in $\fS$, a {\em quantum holonomy} 
$\GH\in \TT$. 
Gabella's construction works for general Lie groups ${\rm SL}_N$, where the algebra $\TT$ is replaced by the Fock-Goncharov quantum moduli space. The construction is
based on works of Gaiotto, Moore and Neitzke \cite{G09, GMN11, GMN13, GMN14} and Galakhov, Longhi and Moore \cite{GLM}, which give a correspondence between the ${\rm GL}_N$ holonomies for the surface $\frak{S}$ and the abelian ${\rm GL}_1$ holonomies for certain $N$-fold branched covers of $\frak{S}$.   
Building on these works, and combining with the idea of Bonahon and Wong about going to three dimensions, Gabella constructed the quantum holonomy $\GH$. In this paper we  
 focus only on the case $N=2$. 

\vs

\def\cV{\mathcal V}
\def\TS{\widetilde {\fS}_\Delta}
\def\TrHol{\mathrm{TrHol}}
\def\al{\alpha}
\def\gm{\gamma}
\def\tgm{\til{\gamma}}
\def\wri{\mathrm{wr}}
Let us briefly and informally discuss the main ingredients in Gabella's construction. Details will be given in Section \ref{sec:Gabella}. Let $\pi: \TS \to \fS$ be a branched double covering of $\fS$ branched over a set $\cV\subset \fS$ consisting of exactly one point in the interior of each triangle of the triangulation $\Delta$. Let $\fS'=\fS\setminus \cV$. Let $\gm:[0,1] \to \fS$ be a (not necessarily closed) smooth oriented curve in general position. When the image of $\gamma$ lies in $\fS'$, there are exactly two lifts of $\gamma$ in $\TS$. If $\gamma$ is considered up to homotopy in $\frak{S}$ fixing its endpoints, then it has many more mutually nonhomotopic lifts in $\TS$. In Gaiotto, Moore and Neitzke's construction, based on the triangulation one chooses a certain finite number of lifts, called {\em admissible lifts}, of the homotopy class of $\gm$ in $\fS$. Then, the ${\rm GL}_2$ parallel transport for $\gamma$ is expressed as sum over all admissible lifts $\tgm$ of a formal variable $X_{\tgm}$. Galakhov, Longhi and Moore proposed a quantum parallel transport$$ \sum_{\text{admissible lifts}\,\, \tgm} q^{-\wri(\tgm) } \wh{X}_{\tgm},$$
where $\wri(\tgm)\in \mathbb{Z}$ is the {\em writhe} of the lift $\tgm$, i.e. the number of self-intersections of $\tgm$ counted with signs, and $\wh{X}_{\tgm}$ is a certain formal variable associated to $\tgm$. This construction works well for open curves $\gamma$. When $\gamma$ is a closed oriented curve in $\frak{S}$, one can represent it by a map $\gamma:[0,1] \to \frak{S}$ by choosing a basepoint, then apply the above construction, where the variables $\wh{X}_{\tgm}$ become monomials in $\mathcal{T}^\omega_\Delta$; in general, the result depends on the choice of basepoint. Gabella suggested a modification in order to remove this dependency, and hence to construct a quantum holonomy of a simple closed oriented curve $\gamma$ in $\frak{S}$. Assume $\gamma$ is in general position with respect to the arcs of $\Delta$. Choose an arc $e$ of $\Delta$ and cut $\gamma$ along $e$ to get a collection of open curves $\gamma_i$. Apply the Galakhov-Longhi-Moore construction to each of the $\gamma_i$. The intersection of $\gamma$ and a small neighborhood of $e$  is then interpreted as a certain type of tangle in the thickening of this neighborhood, which can be viewed as a {\em biangle}, and a modification of the operator invariant \`{a} la Reshetikhin-Turaev \cite{RT, Turaev16} of this tangle is used to ``glue'' the quantum holonomy of the arcs $\gamma_i$ to give the Gabella quantum holonomy $\TrHol_\Delta^\omega(\gamma)\in \TT$ \cite{G}. Gabella showed that $\TrHol_\Delta^\omega(\gamma)$ is an isotopy invariant of simple closed oriented curves $\gamma$ in $\fS$.  It should be noted that while the Bonahon-Wong quantum trace map is defined for unoriented framed links in the thickening of $\fS$, the Gabella quantum holonomy is defined for oriented simple loops in $\fS$. The definition of the Gabella quantum holonomy, where lifting of curves and writhe numbers are used, is very different from that of the Bonahon-Wong quantum trace.
 
 \vs
 
 A special case of the main theorem is:
 \begin{theorem}[Part of Theorem \ref{thm:main}] Suppose $\fS$ is a punctured surface with an ideal triangulation $\Delta$. Let $\gamma$ be a non-null-homotopic oriented simple loop in $\fS$. 
 The Bonahon-Wong quantum trace and the Gabella quantum holonomy for $\gamma$ coincide with each other:
$$ \Tr^\omega_\Delta(\gamma) = \GH.$$
Here in the left hand side we consider $\gamma$ as an unoriented link with the blackboard framing.
 \end{theorem}
 In particular, the deformation quantization map $\mathbb{I}^q_\Delta$ of Allegretti and Kim using the Bonahon-Wong quantum trace coincides with Gabella's solution (after a certain modification to be explained) for the trace-of-monodromy function $\mathbb{I}(\gamma)\in {\bf L}(\mathcal{X}^+(\frak{S}))$ of a non-null-homotopic simple loop $\gamma$ in $\frak{S}$.
 
 \def\cP{\mathcal P}
 \def\cS{\mathcal S}
\subsection{More general result, and proof} Actually we will prove a much more general result. First we show that one can extend the Gabella quantum holonomy not only to all framed oriented links in the thickening of $\fS$, but also to a larger class of what we call {\em stated VH-tangles} in surfaces with boundary, and for VH-tangles we show that the Bonahon-Wong quantum trace and the Gabella quantum holonomy are not exactly equal, but related by a certain relation. The more general result for VH-tangles is actually needed for the proof, as we decompose the surface into triangles and biangles, and establish the relation in each biangle and each triangle. Then we show that we can ``glue'' the relations obtained in the biangles and triangles to get the desired global relation.

\vs

First we allow the surface $\fS$ to have boundary. So $\fS=\Sigma\setminus \cP$, where $\Sigma$ is a compact oriented surface with (possibly empty) boundary and $\cP$ is a finite set that intersects every connected component of the boundary of $\Sigma$. By a {\em tangle} $K$ in the thickened surface $\frak{S}\times(-1,1)$ we mean an unoriented compact properly embedded framed $1$-dimensional submanifold of $\frak{S} \times (-1,1)$. A {\em V-tangle} is a tangle $K$ such that for each boundary arc $b$ of $\frak{S}$, the endpoints of $K$ in $b\times (-1,1)$ have distinct vertical coordinates, i.e. distinct elevations in $(-1,1)$. A {\em VH-tangle} is a V-tangle $K$ such that for each boundary arc $b$, the endpoints of $K$ in $b\times (-1,1)$ have distinct horizontal coordinates, i.e. project down to distinct points in $b$.  By {\em V-isotopy} and {\em VH-isotopy} we mean the isotopies within the respective classes. A {\em state} of a V-tangle or VH-tangle is a map from its boundary points to the set $\{+, -\}$. Bonahon and Wong \cite{BW} extend the definition of the skein algebra to the stated skein algebra $\cS^A_s(\fS)$ based on stated V-tangles, and show that the quantum trace map $\Tr^\omega_\Delta: \cS^A_s(\fS) \to \TT$ can be defined  as an algebra homomorphism.

\vs

In Gabella's definition the writhe plays an important role, but the V-isotopy can change the writhe. This is the reason why  we have to consider the finer class of VH-tangles, as VH-isotopy does not change the writhe. For each VH-isotopy class of a stated  VH-tangle $\al$ in $\frak{S}\times(-1,1)$ we will define the Gabella quantum holonomy $\TrHol_\Delta^\omega(\al)\in \TT$, extending Gabella's definition for simple closed curves in $\frak{S}$.

\vs

If we forget the H-structure in a stated VH-tangle $\al$, we get a stated V-tangle. There is  a simple quantity $\partial \mathcal{C}_\frak{S}(\al)\in \mathbb Z$ (Def.\ref{def:signed_order_correction_amount}) which records some information from this forgetting process.  The main result of the present paper is the following.

\begin{theorem}[main theorem]
For a stated VH-tangle $\al$ in the thickened surface $\frak{S}\times(-1,1)$,
\begin{align}
\nonumber
{\rm Tr}^\omega_\Delta(\al) = \omega^{2{\rm wr}(\alpha)} \, \omega^{\partial \mathcal{C}_{\frak{S}}(\al)} \, {\rm TrHol}^\omega_\Delta(\al),
\end{align}
where ${\rm wr}(\al)$ is the usual writhe of the tangle $\al$ in $\frak{S}\times(-1,1)$.
\end{theorem}
We briefly explain the idea of the proof. 
For a biangle, we start from the well-known Reshetikhin-Turaev operator invariant \cite{RT, Turaev16}, and show that a modification by a certain boundary twist and writhe still yields an operator invariant. We then show that the matrix elements of these operator invariants are precisely the Bonahon-Wong biangle quantum trace and Gabella biangle quantum holonomy. For a triangle, we first write the Bonahon-Wong's triangle factor as the so-called Weyl-ordered monomial times some power of $\omega$ which we call a {\em deviation}. Lemma \ref{lem:equality_of_the_triangle_factors}, our main technical lemma, is the formula expressing this deviation in terms of a certain combination of writhe in $\frak{S}$, writhe of (Gabella) lifts in $\til{\frak{S}}_\Delta$, and a sign correction at the sides of the triangle. In the entire surface, we show that the boundary twists of biangles and the sign corrections of triangles cancel each other at inner arcs of $\Delta$, whereas the two kinds of writhes add up throughout the surface. 

\vs

While the twist factor $\omega^{2{\rm wr}(\alpha)} \, \omega^{\partial \mathcal{C}_{\frak{S}}(\al)}$ is indeed what makes possible the forgetting process from VH-tangles to V-tangles, it also records the orientation information of the tangle $\alpha$. In particular, even if $\alpha$ is a closed tangle, so that it has no V-, H-, or VH-structure, nor a state, the Bonahon-Wong quantum trace and the Gabella quantum holonomy are not exactly the same, but are related by ${\rm Tr}^\omega_\Delta(\alpha) = \omega^{2{\rm wr}(\alpha)} {\rm TrHol}^\omega_\Delta(\alpha)$.

\vs

Finally, we note that the present version of the paper substantially improves the previous one \cite{KSv1}, with a major new input by L\^e who joined the authorship. The construction and properties of the biangle quantum holonomy were stated only as a conjecture in \cite{KSv1}, and now we resolve that conjecture by using the Reshetikhin-Turaev invariant. Thus we provide a mathematically rigorous treatment of Gabella's ${\rm SL}_2$ quantum holonomy, whose properties and even the well-definedness had not been completely proved. We reformulate the construction of Gabella \cite{G} so that it is more easily understandable to the mathematics community, and provide suitable necessary corrections. The new results of the present paper enable us to extend the construction of Gabella to a larger class of tangles (namely all stated VH-tangles), as well as to show that the result is essentially equal to Bonahon-Wong's quantum trace \cite{BW}. This equality was proved in \cite{KSv1} only for lifts of simple loops, by using an argument that is not extendable to general tangles. We prove the equality for general tangles by a more direct and conceptual proof.

\vs

\noindent{\bf Acknowledgments.} L\^e is partially supported by the National Science Foundation under grant \ul{DMS} 1811114. This research was supported by Basic Science Research Program through the National Research Foundation of Korea(NRF) funded by the Ministry of Education(grant number 2017R1D1A1B03030230). This work was supported by the National Research Foundation of Korea(NRF) grant
funded by the Korea government(MSIT) (No. 2020R1C1C1A01011151). Kim has been supported by a KIAS Individual Grant (MG047203) at Korea Institute for Advanced Study. Kim thanks Dylan Allegretti, Carlos Scarinci, Linhui Shen, Maxime Gabella and Francis Bonahon for help, discussion, questions, and comments. The authors thank the anonymous referees and the editorial staff for helpful suggestions.

\section{Quantum Teichm\"uller theory}

We establish basic definitions, and briefly introduce quantum Teichm\"uller theory.

\subsection{Ideal triangulation of a noncompact surface}

In this subsection we choose to mostly follow the terminology conventions of \cite{Le}, with some modification.
\begin{definition}[\cite{Le}]
\label{def:surface}
$\bullet$ A \ul{\em generalized marked surface} is a pair $(\Sigma,\mathcal{P})$, where $\Sigma$ is a compact oriented surface with boundary $\partial \Sigma$ and $\mathcal{P}$ is a finite subset of $\Sigma$, such that each component of $\partial \Sigma$ contains at least one point of $\mathcal{P}$.

\vs

$\bullet$ The surface orientation on $\Sigma$ determines the notion of clockwise/counterclockwise rotation on (the tangent space of each point of) the surface, and the convention of the \ul{\em boundary-orientation} on the boundary $\partial \Sigma$ coming from the surface orientation is chosen to be compatible with the {\em clockwise} rotation on (points of) the surface near its boundary.

\vs

$\bullet$ Elements of $\mathcal{P}$ are called \ul{\em marked points}. Elements of $\mathcal{P}$ in $\Sigma \setminus \partial \Sigma$ are called \ul{\em interior marked points}, or \ul{\em punctures}.

\vs

$\bullet$ If $\partial \Sigma = {\O}$, then $(\Sigma,\mathcal{P})$ is called a \ul{\em punctured surface}.

\vs

$\bullet$ By the \ul{\em boundary} of the generalized marked surface $(\Sigma,\mathcal{P})$ we mean $(\partial \Sigma) \setminus \mathcal{P}$, and write it as $\partial \Sigma \setminus \mathcal{P}$.

\vs

$\bullet$ A \ul{\em diffeomorphism} between generalized marked surfaces $(\Sigma,\mathcal{P})$ and $(\Sigma',\mathcal{P}')$ is an orientation-preserving diffeomorphism $\Sigma \to \Sigma'$ that restricts to a bijection $\mathcal{P} \to \mathcal{P}'$. We say $(\Sigma,\mathcal{P})$ and $(\Sigma',\mathcal{P}')$ are \ul{\em diffeomorphic} if there exists a diffeomorphism between them.
\end{definition}
As mentioned in \cite{Le}, a generalized marked surface defined this way corresponds to a ``punctured surface with boundary" of \cite{BW}.

\begin{definition}
By ``with boundary" we always mean ``with (possibly empty) boundary". 
\end{definition}
In particular, in Def.\ref{def:surface}, we allow $\partial \Sigma$ to be empty. Notice that, in Def.\ref{def:surface}, when $\partial \Sigma$ is empty the boundary $\partial \Sigma \setminus \mathcal{P}$ of $(\Sigma,\mathcal{P})$ is also empty.

\begin{definition}[\cite{Le}]
\label{def:P-link}
Let $(\Sigma,\mathcal{P})$ be a generalized marked surface.

$\bullet$ A \ul{\em $\mathcal{P}$-link} is an immersion $\alpha : C \to \Sigma$, where $C$ is a compact unoriented $1$-manifold with boundary, such that
\begin{enumerate}
\item[\rm (1)] restriction of $\alpha$ onto the interior of $C$ is an embedding into $\Sigma\setminus \mathcal{P}$, and

\item[\rm (2)] $\alpha$ maps the boundary of $C$ into $\mathcal{P}$.
\end{enumerate}
The \ul{\em interior} of this $\mathcal{P}$-link is the image under $\alpha$ of the interior of $C$, and is denoted by $\mathring{\alpha}$.

\vs

$\bullet$ When $C$ is $[0,1]$, we call $\alpha$ a \ul{\em $\mathcal{P}$-arc}. When $C$ is $S^1$, we call $\alpha$ a \ul{\em $\mathcal{P}$-knot}.

\vs

$\bullet$ When $C$ has no boundary, we call $\alpha$ a \ul{\em closed} $\mathcal{P}$-link.

\vs

$\bullet$ Two $\mathcal{P}$-links are \ul{\em $\mathcal{P}$-isotopic} if they are isotopic in the class of $\mathcal{P}$-links.

\vs

$\bullet$ A $\mathcal{P}$-arc is called a \ul{\em boundary arc} if it is contained in $\partial \Sigma$.

\vs

$\bullet$ A $\mathcal{P}$-link is often identified with its image in $\Sigma$. 
\end{definition}
In the literature, the surface is commonly considered to be $\Sigma\setminus \mathcal{P}$, in which case a $\mathcal{P}$-arc corresponds to a so-called {\em ideal arc}, a $\mathcal{P}$-knot to a simple loop in $\Sigma\setminus \mathcal{P}$, and a closed $\mathcal{P}$-link to a simple closed curve in $\Sigma\setminus\mathcal{P}$. The last line of Def.\ref{def:P-link} indicates that a $\mathcal{P}$-link is an unoriented object. One easy observation is that the boundary $\partial \Sigma\setminus \mathcal{P}$ is the disjoint union of interiors of boundary arcs, which are the connected components of $\partial \Sigma\setminus \mathcal{P}$. The interior of a boundary arc may also be called a boundary arc, by abuse of notation. We now go on to triangulations. 

\begin{definition}[\cite{Le, FST}]
\label{def:triangulability}
A generalized marked surface $(\Sigma,\mathcal{P})$ is \ul{\em triangulable} if and only if $\mathcal{P}$ is nonempty and $(\Sigma,\mathcal{P})$ is not diffeomorphic to one of the following:

-- a sphere (with no boundary) with one or two marked points,

-- a monogon with no interior marked point, (genus zero, with one boundary component, with one marked point on the boundary, and no interior marked point), or

-- a bigon with no interior marked point (genus zero, with one boundary component, with two marked points on the boundary, and no interior marked point), also called a \ul{\em biangle}.
\end{definition}

We do not always assume that $(\Sigma,\mathcal{P})$ is triangulable, although we do when we talk about triangulations. Later, it is crucial that we consider the third case of the above definition.

\begin{deflem}[\cite{Le, FST}]
\label{def:P-triangulation}
Let $(\Sigma,\mathcal{P})$ be a triangulable generalized marked surface.

\vs

$\bullet$ For integer $n\ge 3$, a \ul{\em $\mathcal{P}$-$n$-gon} is a smooth map $\beta : \sigma \to \Sigma$ from a regular $n$-gon $\sigma$ in the standard plane $\mathbb{R}^2$ to $\Sigma$, such that
\begin{enumerate} 
\item[\rm (1)] the restriction of $\beta$ onto the interior $\mathring{\sigma}$ is a diffeomorphism onto its image, and

\item[\rm (2)] the restriction of $\beta$ onto each of the $n$ sides of $\sigma$ is a $\mathcal{P}$-arc, called a \ul{\em side} (or \ul{\em edge}) of $\beta$.
\end{enumerate}
In particular, a $\mathcal{P}$-$3$-gon is called a \ul{\em $\mathcal{P}$-triangle}.

\vs

$\bullet$ If two sides of a $\mathcal{P}$-$n$-gon coincide as a $\mathcal{P}$-arc, such a side, as well as that $\mathcal{P}$-$n$-gon, is said to be \ul{\em self-folded}.

\vs

$\bullet$ A $\mathcal{P}$-$n$-gon is \ul{\em oriented} if the relevant map $\beta:\sigma\to \Sigma$ is orientation preserving, where $\sigma$ is oriented according to the standard orientation on $\mathbb{R}^2$.

\vs

$\bullet$ A \ul{\em $\mathcal{P}$-triangulation} (or, a \ul{\em triangulation}, when $\mathcal{P}$ is clear) of $\Sigma$ (or, an \ul{\em ideal triangulation} or just a triangulation of $(\Sigma,\mathcal{P})$) is a collection $\Delta$ of $\mathcal{P}$-arcs such that
\begin{enumerate}
\item[\rm (1)] no $\mathcal{P}$-arc in $\Delta$ bounds a disk whose interior is in $\Sigma\setminus\mathcal{P}$,

\item[\rm (2)] no two $\mathcal{P}$-arcs in $\Delta$ intersect in $\Sigma\setminus \mathcal{P}$ and no two are $\mathcal{P}$-isotopic, and

\item[\rm (3)] $\Delta$ is a maximal collection of $\mathcal{P}$-arcs with the above properties (1) and (2).
\end{enumerate}
For a triangulation $\Delta$, one can replace $\mathcal{P}$-arcs in $\Delta$ by $\mathcal{P}$-arcs in their respective $\mathcal{P}$-isotopy classes so that every arc in $\Delta$ isotopic to a boundary arc is a boundary arc. We always assume that $\Delta$ satisfies this condition.

\vs

$\bullet$ An element of $\Delta$ is called a \ul{\em (constituent) edge} of the triangulation $\Delta$. An edge of $\Delta$ is called a \ul{\em boundary edge} of $\Delta$ if it is a boundary arc, and an \ul{\em internal edge} of $\Delta$ otherwise. Let $\mathring{\Delta}$ be the set of all internal edges of $\Delta$.

\vs

$\bullet$ For a $\mathcal{P}$-triangulation $\Delta$, the closure of each connected component of $\Sigma\setminus (\bigcup_{e\in \Delta}e)$ can be naturally given a structure of an oriented $\mathcal{P}$-triangle. By a \ul{\em (constituent) ($\mathcal{P}$-)triangle of $\Delta$} we mean one of these triangles coming from $\Delta$. Denote by $\mathcal{F}(\Delta)$ the set of all triangles of $\Delta$.

\vs

$\bullet$ We say $(\Sigma,\mathcal{P})$ is \ul{\em triangulated} if it is equipped with a $\mathcal{P}$-triangulation.
\end{deflem}
The notions of $\mathcal{P}$-triangulation and $\mathcal{P}$-triangle correspond to those of {\em ideal triangulation} and {\em ideal triangle} in the literature. Implicitly or explicitly, we will identify two triangulations if one can be obtained by simultaneous isotopy of (edges of) the other; when extra care is needed, we shall make it clear.  For convenience of the present paper, we define the notion of the corner of a triangle as follows.
\begin{deflem}
\label{def:corners}
Let $t$ be a constituent triangle of a triangulation $\Delta$ of a generalized marked surface $(\Sigma,\mathcal{P})$. Let $\beta : \sigma \to \Sigma$ be a map giving $t$ an oriented $\mathcal{P}$-triangle structure. Let $e$, $f$ and $g$ be the sides of $\sigma$, appearing clockwise in this order. Denote by same labels $e$, $f$ and $g$ the corresponding sides of $t$. The three pairs of sides $(e,f)$, $(f,g)$, and $(g,e)$ are called \ul{\em corners} of $t$. The three corners of $t$ are well-defined and mutually distinct.
\end{deflem}
Visually, one can consider a corner $(e,f)$ of $t$ as being a small part of $t \setminus (e\cup f\cup g)$ close to the vertex of $t$ shared by $e$ and $f$, lying in between $e$ and $f$.  There may arise some confusion when dealing with a self-folded triangle. In the above definition, we give three different labels $e$, $f$ and $g$ for the three sides; however sometimes we identify a side with its image, in which case we will need only two distinct labels for the sides of a self-folded triangle, say $e$, $e$ and $f$. Notice that, still in such a case, the three corners $(e,e)$, $(e,f)$ and $(f,e)$ of this triangle are unambiguously defined. However, a pair of edges of $\Delta$ may not be sufficient to determine a corner, because sometimes such a pair may represent corners of two different triangles. 

\vs

We find it convenient to define the following notion here:
\begin{definition}
\label{def:minimal_position}
Let $\Delta$ be a triangulation of a generalized marked surface $(\Sigma,\mathcal{P})$, and let $\gamma$ be a closed $\mathcal{P}$-link in $\Sigma$. Denote by $\#(\gamma\cap \Delta)$ the total number of intersections of $\gamma$ and the edges of $\Delta$; it can be infinite. We say that $\gamma$ is in a \ul{\em minimal position} with respect to $\Delta$ if $\#(\gamma \cap \Delta)$ is minimal among the numbers $\#(\gamma\hspace{0,2mm}' \cap \Delta)$, where $\gamma\hspace{0,2mm}'$ runs through all $\mathcal{P}$-links that are $\mathcal{P}$-isotopic to $\gamma$.
\end{definition}

It is easy to see that if $\gamma$ is in a minimal position with respect to $\Delta$, then $\#(\gamma\cap \Delta)$ is finite. It is well-known that the intersection numbers $\#(\gamma \cap e)$ for the edges $e \in \Delta$ completely determine the $\mathcal{P}$-link $\gamma$ in a minimal position, up to $\mathcal{P}$-isotopy.
\begin{lemma}
\label{lem:a_e}
Let $\Delta$ be a triangulation of a generalized marked surface $(\Sigma,\mathcal{P})$, and let $\gamma$ be a closed $\mathcal{P}$-link in $\Sigma$. Let $\gamma'$ be a closed $\mathcal{P}$-link that is $\mathcal{P}$-isotopic to $\gamma$ and is in a minimal position with respect to $\Delta$. For each $e\in \Delta$, denote by $a_e(\gamma) \in \mathbb{Z}_{\ge 0}$ the number of intersections of $\gamma'$ and $e$, called the \ul{\em intersection number} of $\gamma$ and $e$.

\vs

The intersections numbers completely determine a closed $\mathcal{P}$-link $\gamma$ up to $\mathcal{P}$-isotopy.
\end{lemma}
For a proof, we refer to \cite{FG06}; the intersection numbers are two times the Fock-Goncharov tropical $\mathcal{A}$-coordinate of $\gamma$, and it is straightforward to reconstruct $\gamma$ out of these numbers. 

\vs

As is widely used in the literature, counting of corners effectively encodes the combinatorics of $\Delta$.
\begin{definition}
\label{def:exchange_matrix}
Let $\Delta$ be a triangulation of a generalized marked surface $(\Sigma,\mathcal{P})$. Let $e$ and $f$ be constituent edges of $\Delta$, and let $t$ be a constituent triangle of $\Delta$. Define
\begin{align*}
c_{ef}(t) & = c^{\Delta}_{ef} (t) := \left\{
\begin{array}{rl}
1 & \mbox{if $(e,f)$ is a corner of $t$}, \\
0 & \mbox{otherwise},
\end{array}
\right. \\
\varepsilon_{ef}(t) & = \varepsilon^\Delta_{ef}(t) := c_{ef}(t) - c_{fe}(t), \\
\varepsilon_{ef} & = \varepsilon^\Delta_{ef} := \sum_{t\in \mathcal{F}(\Delta)} \varepsilon_{ef}(t).
\end{align*}
The matrix $\varepsilon = \varepsilon^\Delta = (\varepsilon_{ef})_{e,f\in \Delta}$ is called the \ul{\em exchange matrix} of $\Delta$.
\end{definition}
The matrix $\varepsilon$ is sometimes called the signed adjacency matrix (see \cite{FST}), or a face matrix (see \cite{Le}). Notice that the matrix $B = (b_{ij})$ of \cite{FST} coincides with the above $\varepsilon =  (\varepsilon_{ij})$ when $\Delta$ has no self-folded triangles, and otherwise it differs from $\varepsilon$ in general.

\subsection{Teichm\"uller space and shear coordinates}

The basic geometric objects of study are certain versions of the so-called Teichm\"uller space of a surface $(\Sigma,\mathcal{P})$. In general, the Teichm\"uller space of a surface refers to the space of all complete hyperbolic metrics on the surface, considered up to pullback by diffeomorphisms isotopic to the identity. In the case of a generalized marked surface, some care is needed for the behavior near the boundary and marked points. The most general case appears in Dylan Allegretti's thesis \cite{A}, which we follow. Like in \cite{A}, we formulate using the monodromy representations instead of directly dealing with hyperbolic metrics.
\begin{definition}[enhanced Teichm\"uller space of a punctured surface \cite{BL, F, FG06, FG07, Liu}]
\label{def:enhanced_Teichmuller_space_punctured}
Let $(\Sigma,\mathcal{P})$ be a triangulable generalized marked surface, and suppose that $(\Sigma,\mathcal{P})$ is a punctured surface, i.e. $\partial \Sigma = {\O}$.

\vs

$\bullet$ The \ul{\em Teichm\"uller space} $\mathcal{T}(\Sigma,\mathcal{P})$ is the set
\begin{align}
\label{eq:Teichmuller_space}
{\rm Hom}^{\rm df,\Sigma}(\pi_1(\Sigma\setminus\mathcal{P}), {\rm PSL}(2,\mathbb{R}))/{\rm PSL}(2,\mathbb{R}).
\end{align}
More precisely, $\mathcal{T}(\Sigma,\mathcal{P})$ consists of all faithful group homomorphisms from $\pi_1(\Sigma\setminus\mathcal{P})$ to ${\rm PSL}(2,\mathbb{R})$ that have discrete image such that the quotient $\mathbb{H}/\rho(\pi_1(\Sigma\setminus\mathcal{P}))$ of the upper half-plane is homeomorphic to $\Sigma\setminus\mathcal{P}$, where an element of $\mathcal{T}(\Sigma,\mathcal{P})$ is defined up to conjugation by an element of ${\rm PSL}(2,\mathbb{R})$. 

\vs 

$\bullet$ Given a point of $\mathcal{T}(\Sigma,\mathcal{P})$, i.e. an equivalence class of a group homomorphism $\rho : \pi_1(\Sigma\setminus\mathcal{P})\to {\rm PSL}(2,\mathbb{R})$, an element $p$ of $\mathcal{P}$ is called a \ul{\em cusp} (with respect to $\rho$) if the image under $\rho$ of a small loop surrounding $p$ is a parabolic element of ${\rm PSL}(2,\mathbb{R})$, and called a \ul{\em hole} (with respect to $\rho$) if this image is a hyperbolic element of ${\rm PSL}(2,\mathbb{R})$.

\vs

$\bullet$ The \ul{\em enhanced Teichm\"uller space} $\mathcal{X}^+(\Sigma,\mathcal{P})$ is the set of all pairs $(\rho,O)$, where $\rho$ is a point of $\mathcal{T}(\Sigma,\mathcal{P})$, and $O$ is the choice of an orientation for (a loop surrounding) each hole with respect to $\rho$; using the surface orientation of $\Sigma$, the data of $O$ can be represented as a sign for each hole (clockwise or counterclockwise orientation of a loop surrounding the hole).
\end{definition}

\begin{definition}[enhanced Teichm\"uller space for a triangulable generalized marked surface; \cite{A}]
\label{def:enhanced_Teichmuller_space}
Let $(\Sigma,\mathcal{P})$ be a triangulable generalized marked surface with nonempty boundary $\partial \Sigma \neq {\O}$. 

\vs

$\bullet$ Let $(\Sigma^{\rm op}, \mathcal{P}^{\rm op})$ be the same generalized marked surface as $(\Sigma,\mathcal{P})$, except equipped with the opposite orientation. Choose a parametrization for each boundary component of $\Sigma$, which induces a parametrization for the corresponding boundary component of $\Sigma^{\rm op}$. Glue $\Sigma$ and $\Sigma^{\rm op}$ along these parametrized boundary components to construct a smooth oriented surface $\Sigma_{\rm D}$ without boundary; denote by $\mathcal{P}_{\rm D}$ the marked points on $\Sigma_{\rm D}$ coming from $\mathcal{P}$ and $\mathcal{P}^{\rm op}$. The resulting generalized marked surface $(\Sigma_{\rm D}, \mathcal{P}_{\rm D}) =: (\Sigma,\mathcal{P})_{\rm D}$, which is a punctured surface, is called the \ul{\em doubled surface} for $(\Sigma,\mathcal{P})$. Let $\iota : \Sigma_{\rm D} \to \Sigma_{\rm D}$ be the natural involutive diffeomorphism exchanging $\Sigma$ and $\Sigma^{\rm op}$, induced by the identification $\Sigma \leftrightarrow \Sigma^{\rm op}$.

\vs

$\bullet$ The \ul{\em (generalized) enhanced Teichm\"uller space} $\mathcal{X}^+(\Sigma,\mathcal{P})$ is defined as the $\iota$-invariant subspace of the enhanced Teichm\"uller space $\mathcal{X}^+(\Sigma_{\rm D}, \mathcal{P}_{\rm D})$ of the doubled surface $(\Sigma_{\rm D}, \mathcal{P}_{\rm D})$, where $\mathcal{X}^+(\Sigma_{\rm D}, \mathcal{P}_{\rm D})$ is as in Def.\ref{def:enhanced_Teichmuller_space_punctured}, and $\iota$ acts on $\mathcal{X}^+(\Sigma_{\rm D}, \mathcal{P}_{\rm D})$ naturally on the hyperbolic metrics and by reversing the (signs of the) orientations of the holes.
\end{definition}
As a consequence of the above definition, for a point of $\mathcal{X}^+(\Sigma_{\rm D},\mathcal{P}_{\rm D})$ yielding a point of $\mathcal{X}^+(\Sigma,\mathcal{P})$, i.e. for an $\iota$-invariant point of $\mathcal{X}^+(\Sigma_{\rm D},\mathcal{P}_{\rm D})$, each marked point of $\mathcal{P}_{\rm D}$ corresponding to a marked point of $\mathcal{P}$ lying in the boundary of $\Sigma$ is a cusp\footnote{We thank the referee and Dylan Allegretti for pointing this out.}.

\vs

A starting point of many problems related to Teichm\"uller spaces is a construction of suitable coordinate systems for $\mathcal{X}^+(\Sigma,\mathcal{P})$. Among various kinds, what is relevant to the current situation is Thurston's shear coordinate function \cite{BL, F, FG07, Liu, P, Thurston}, which makes use of the choice of an ideal triangulation.
\begin{proposition}[Thurston-Fock theorem \cite{A, F, FG07, P, Thurston}]
Let $\Delta$ be a triangulation of a triangulable generalized marked surface $(\Sigma,\mathcal{P})$. There exists a global coordinate system for the enhanced Teichm\"uller space $\mathcal{X}^+(\Sigma,\mathcal{P})$, whose coordinate functions are enumerated by the internal edges $e$ of the triangulation $\Delta$. These coordinate functions $X_e = X_e^\Delta : \mathcal{X}^+(\Sigma,\mathcal{P}) \to \mathbb{R}_{>0}$ associated to $e\in \mathring{\Delta}$, each of which is called the \ul{\em exponentiated shear coordinate} for the internal edge $e \in \Delta$, provide a bijection
$$
\mathcal{X}^+(\Sigma,\mathcal{P}) \to (\mathbb{R}_{>0})^{\mathring{\Delta}} ~ : ~ (\rho,O) \mapsto X_e(\rho,O).
$$
\end{proposition}
Roughly speaking, in terms of hyperbolic metrics, these coordinates are the ``shearing" amounts that measure how the hyperbolic ideal triangles are glued along their edges to one another to form the (hyperbolic) surface $\Sigma\setminus\mathcal{P}$. A key point of studying these coordinates is how they transform under change of ideal triangulations.

\begin{proposition}[\cite{F, FG07, Liu, P}]
Let $\Delta$ and $\Delta'$ be ideal triangulations of a triangulable generalized marked surface $(\Sigma,\mathcal{P})$. Then the exponentiated shear coordinate functions $X_{e'}'$ for the internal edges $e'$ of $\Delta'$ can be expressed as rational functions in terms of those $X_e$ for the internal edges $e$ of $\Delta$.
\end{proposition}

\subsection{Quantum Teichm\"uller space}
\label{subsec:quantum_Teichmuller_space}

We now review quantum Teichm\"uller space, which was first constructed in the 1990s by Kashaev \cite{Kash} and independently by Chekhov and Fock \cite{F,CF}. For our purposes we follow the latter works, and their modern developments as appearing in \cite{BL, BW, FG09, Hiatt, Liu}.

\vs

For each triangulation $\Delta$, we considered the exponentiated shear coordinate functions $X_e$ associated to internal edges $e \in \mathring{\Delta}$. They provide identification of the (generalized) enhanced Teichm\"uller space $\mathcal{X}^+(\Sigma,\mathcal{P})$ with $(\mathbb{R}_{>0})^{\mathring{\Delta}}$, making $\mathcal{X}^+(\Sigma,\mathcal{P})$ a smooth manifold. This manifold $\mathcal{X}^+(\Sigma,\mathcal{P})$ is equipped with a natural Poisson structure, named the Weil-Petersson Poisson structure, whose Poisson brackets among the coordinate functions are given by
$$
\{X_e, X_f\} = \varepsilon_{ef} X_e X_f.
$$
Thus one can verify that $\{X_e^{\pm 1} : e\in \mathring{\Delta}\}$ generates a Poisson subalgebra of the ring $C^\infty(\mathcal{X}^+(\Sigma,\mathcal{P}))$ of smooth functions on $\mathcal{X}^+(\Sigma,\mathcal{P})$. The following noncommutative algebras serve as quantum deformed versions of this Poisson subalgebra:
\begin{definition}[\cite{BL,Liu}]
\label{def:Chekhov-Fock_algebra}
Let $\Delta$ be a triangulation of a generalized marked surface $(\Sigma,\mathcal{P})$, and $\mathring{\Delta}$ be the set of all internal edges of $\Delta$. Let $\hbar \in \mathbb{R}$ be a quantum parameter, and let $q := \exp(\pi {\rm i} \hbar) \in \mathbb{C}^*$. 

\vs

The \ul{\em internal Chekhov-Fock algebra} $\mathring{\mathcal{T}}^q_\Delta$ is an algebra over $\mathbb{C}$ defined by generators and relations:
\begin{align*}
\mbox{generating set} ~ & : ~ \{\, \wh{X}_e, \, \wh{X}_e^{-1} ~| ~ e\in \mathring{\Delta} \,\}, \\
\mbox{relations} ~ & : ~ \wh{X}_e \wh{X}_f = q^{2\varepsilon_{ef}} \wh{X}_f \wh{X}_e, \quad \forall e,f\in \mathring{\Delta}, \qquad
\wh{X}_e \wh{X}_e^{-1} = \wh{X}_e^{-1} \wh{X}_e =1, \quad \forall e\in \mathring{\Delta}.
\end{align*}

\vs

The \ul{\em Chekhov-Fock algebra} $\mathcal{T}^q_\Delta$ is an algebra over $\mathbb{C}$ defined by generators and relations as follows:
\begin{align*}
\mbox{generating set} ~ & : ~ \{\, \wh{X}_e, \, \wh{X}_e^{-1} ~| ~ e\in \Delta \,\}, \\
\mbox{relations} ~ & : ~ \wh{X}_e \wh{X}_f = q^{2\varepsilon_{ef}} \wh{X}_f \wh{X}_e, \quad \forall e,f\in \Delta, \qquad
\wh{X}_e \wh{X}_e^{-1} = \wh{X}_e^{-1} \wh{X}_e =1, \quad \forall e\in \Delta.
\end{align*}
\end{definition}
Notice that $\mathring{\mathcal{T}}^q_\Delta$ and $\mathcal{T}^q_\Delta$ coincide when $(\Sigma,\mathcal{P})$ is a punctured surface, which is often the case in the literature. As the generators of the classical algebra are enumerated by internal edges, of course $\mathring{\mathcal{T}}^q_\Delta$ is more natural to consider than $\mathcal{T}^q_\Delta$, but we shall also need $\mathcal{T}^q_\Delta$ later. As noted in \cite{BL,BW, Hiatt,Liu}, these algebras $\mathring{\mathcal{T}}^q_\Delta$ and $\mathcal{T}^q_\Delta$ satisfy the so-called Ore Condition of ring theory \cite{BG, Cohn}, and therefore their skew-fields (i.e. division algebras) of fractions ${\rm Frac}(\mathring{\mathcal{T}}^q_\Delta)$ and ${\rm Frac}(\mathcal{T}^q_\Delta)$ can be considered. As written in \cite{BL, Liu}, elements of ${\rm Frac}(\mathcal{T}^q_\Delta)$ are formal fraction expressions $PQ^{-1}$, with $P,Q\in \mathcal{T}^q_\Delta$ and $Q\neq 0$, where two such expressions $P_1 Q_1^{-1}$ and $P_2 Q_2^{-1}$ represent the same element of ${\rm Frac}(\mathcal{T}^q_\Delta)$ if there exist nonzero $S_1,S_2 \in \mathcal{T}^q_\Delta$ such that $P_1 S_1 = P_2S_2$ and $Q_1 S_1 =Q_2S_2$. A product of such expressions can also be written as one fraction $PQ^{-1}$ by algebraic manipulation. Similar holds for ${\rm Frac}(\mathring{\mathcal{T}}^q_\Delta)$.

\begin{proposition}[\cite{CF,Kash, Liu}]
\label{prop:skew-field_isomorphisms}
There exists a skew-field isomorphism
$$
\Phi^q_{\Delta,\Delta'} ~ : ~ {\rm Frac}(\mathcal{T}^q_{\Delta'}) \to {\rm Frac}(\mathcal{T}^q_\Delta)
$$
associated to each pair $\Delta,\Delta'$ of triangulations of a generalized marked surface $(\Sigma,\mathcal{P})$, satisfying:

(1) When $q=1$ it recovers the coordinate change formulas for exponentiated shear coordinates for internal edges,

(2) The consistency relation
$$
\Phi^q_{\Delta,\Delta''} = \Phi^q_{\Delta,\Delta'} \circ \Phi^q_{\Delta',\Delta''}
$$
holds for each triple $\Delta,\Delta',\Delta''$ of triangulations.
\end{proposition}
In fact, the above proposition was first established for the internal Chekhov-Fock algebras only, i.e. for the maps ${\rm Frac}(\mathring{\mathcal{T}}^q_{\Delta'}) \to {\rm Frac}(\mathring{\mathcal{T}}^q_\Delta)$. We note that these original results naturally extend to the Chekhov-Fock algebras as written above.

\vs

Some authors \cite{BL, BW, Liu} use the term \ul{\em quantum Teichm\"uller space} as referring to the quotient of the disjoint union of all ${\rm Frac}(\mathcal{T}^q_\Delta)$ by the equivalence relation given by the identifications $\Phi^q_{\Delta,\Delta'}$. We might also do so from time to time, but we do not make serious use of this term.
 
\subsection{Deformation quantization of the Teichm\"uller space and the quantum ordering problem}
\label{subsec:deformation_quantization}

The isomorphisms in Prop.\ref{prop:skew-field_isomorphisms} let us identify ${\rm Frac}(\mathcal{T}^q_\Delta)$ for different $\Delta$ in a consistent way. However, it is ${\rm Frac}(\mathcal{T}^q_\Delta)$ that is being identified with others, instead of the original Chekhov-Fock algebra $\mathcal{T}^q_\Delta$ which we began with. For certain reasons elements of ${\rm Frac}(\mathcal{T}^q_\Delta)$ belonging to $\mathcal{T}^q_\Delta$ are considered more important. First, in order to be physically relevant, one must realize elements of the quantum algebra as operators on a Hilbert space. It is relatively easy to deal with representation of the algebra $\mathcal{T}^q_\Delta$ on the Hilbert space using functional analysis, but hard to do for its skew-field of fractions. The second reason comes from the viewpoint of considering the enhanced Teichm\"uller space as a cluster $\mathcal{X}$-variety. Each triangulation yields a chart, and regular functions on this chart are defined as Laurent polynomials in the exponentiated shear coordinates; we note that Laurent polynomials play a crucial role in the theory of cluster algebras and cluster varieties. Quantum regular functions for the quantized chart are noncommutative Laurent polynomials in the quantum counterpart of the exponentiated shear coordinates.

\vs

Hence, for both reasons coming from representation theory and cluster variety theory, the nice functions to deal with, for both the classical and the quantum cases, are functions that are regular for every chart, i.e. that are Laurent polynomials in every triangulation. 
\begin{definition}
\label{def:bf_L_q}
Define
$$
{\bf L}^q_\Delta := \bigcap_{\Delta'} \Phi^q_{\Delta,\Delta'}(\mathcal{T}^q_{\Delta'}) \quad \subset \quad \mathcal{T}^q_\Delta \quad \subset \quad {\rm Frac}(\mathcal{T}^q_\Delta).
$$
where the intersection runs through all ideal triangulations $\Delta'$ of the relevant generalized marked surface $(\Sigma,\mathcal{P})$, including $\Delta$. Elements of ${\rm Frac}(\mathcal{T}^q_\Delta)$ that belong to ${\bf L}^q_\Delta$ are said to be \ul{\em universally Laurent}.
\end{definition}

\begin{lemma}
$\Phi^q_{\Delta,\Delta'}$ induces an isomorphism from ${\bf L}^q_{\Delta'}$ to ${\bf L}^q_\Delta$. \qed
\end{lemma}
This means that the algebras ${\bf L}^q_\Delta$ for different $\Delta$ can be consistently identified, so they can be collectively denoted by ${\bf L}^q = {\bf L}^q_{(\Sigma,\mathcal{P})}$. That is, in the style used in \cite{BL, BW, Liu}, we can view ${\bf L}^q = \bigsqcup_{\Delta} {\bf L}^q_\Delta/\hspace{-1mm}\sim$, where $\bigsqcup$ is the set disjoint union and $\sim$ is the equivalence relation coming from $\Phi^q_{\Delta,\Delta'}$.

\vs

The discussion so far is only on the construction of a consistent quantum system related to the classical system, which is the enhanced Teichm\"uller space with the Weil-Petersson Poisson structure. That is, we now have an algebra of quantum observables, namely ${\bf L}^q$, independent of the choice of an ideal triangulation $\Delta$. Our next step is to establish a quantization, which is a map from the algebra of classical observables to the algebra of quantum observables. Namely, given a classical observable, i.e. a smooth function on the enhanced Teichm\"uller space $\mathcal{X}^+(\Sigma,\mathcal{P})$, what quantum observable do we assign to it? A place to begin with is each of the exponentiated shear coordinate functions $X_e$ for a chosen ideal triangulation $\Delta$. A natural candidate for a quantization map is to send each $X_e$ to $\wh{X}_e \in \mathcal{T}^q_\Delta$. Then, what about other functions? At the moment, let us only focus on the functions on $\mathcal{X}^+(\Sigma,\mathcal{P})$ that can be written as Laurent polynomials in $X_e$ for a given $\Delta$, such as $X_e X_f + X_g^{-2}$. To each such Laurent polynomial in $X_e$, we'd like to assign a noncommutative Laurent polynomial in $\wh{X}_e$'s that recovers the classical one when we put $q=1$ and replace each $\wh{X}_e$ by $X_e$. For a fixed $\Delta$, building such an assignment so that it satisfies the axiom of a ``deformation quantization map" is not so hard. The difficult part is to make sure that such a deformation quantization map does not depend on the choice of $\Delta$, in a sense. 

\vs

The first step is to restrict our attention to (classical) functions on $\mathcal{X}^+(\Sigma,\mathcal{P})$ that can be written as Laurent polynomials in the exponentiated shear coordinates for {\em every} ideal triangulation; finding all of them is already a highly nontrivial task, and is accomplished in \cite{FG06}. Then, for each $\Delta$, we must devise a way to assign to each such universally Laurent function a quantum Laurent polynomial, and prove that the resulting quantum Laurent polynomials for different $\Delta$ are related to each other by the quantum mutation maps $\Phi^q_{\Delta,\Delta'}$. 

\vs

To give an intuition, let us consider some simple toy model. To a function $X_e X_f$, what should we assign? Options are $\wh{X}_e \wh{X}_f$, $\wh{X}_f \wh{X}_e$, $q^r \wh{X}_e \wh{X}_f$ for some $r\in \mathbb{Z}$, or maybe it could be even more complicated. Which one is the best choice? For a more general Laurent polynomial, we should choose how to quantize each monomial term. Finding a good choice of quantum Laurent polynomial so that it satisfies certain favorable properties is sometimes referred to as the \ul{\em quantum ordering problem}, as if it is the problem of choosing the order of a product of noncommutative quantum functions.

\vs

For a monomial, there is a well-known standard answer, namely the Weyl-ordered product; we formulate this in a general setting:
\begin{definition}[Weyl-ordered product]
\label{def:Weyl-ordered_product}
Let ${\bf X}_1,\ldots,{\bf X}_n$ be elements living in an algebra, and satisfying
$$
{\bf X}_i {\bf X}_j = A^{2m_{ij}} {\bf X}_j {\bf X}_i, \quad \forall i,j \in \{1,\ldots,n\},
$$
for some invertible scalar $A$ and integers $m_{ij} \in \mathbb{Z}$. In such a situation, we say that these $n$ elements \ul{\em $A$-commute} with one another. Define the \ul{\em Weyl-ordered product} of these $A$-commuting elements as
$$
[{\bf X}_1 {\bf X}_2 \cdots {\bf X}_n]_{\rm Weyl} := A^{-\sum_{1\le i<j\le n} m_{ij} } {\bf X}_1 {\bf X}_2 \cdots {\bf X}_n.
$$
\end{definition}
It is a straightforward exercise to show that the Weyl-ordered product is invariant under permutation, namely $[{\bf X}_1  \cdots {\bf X}_n]_{\rm Weyl}$ equals $[{\bf X}_{\sigma(1)}  \cdots {\bf X}_{\sigma(n)}]_{\rm Weyl}$  for each permutation $\sigma$ of $\{1,\ldots,n\}$; it makes the Weyl-ordered product a standard answer to the quantum ordering problem of a monomial. However, for our case, for a classical function given by a Laurent polynomial, the quantum Laurent polynomial obtained by replacing each constituent monomial $X_e X_f \cdots$ by the Weyl-ordered product $[\wh{X}_e \wh{X}_f \cdots]_{\rm Weyl}$ turns out {\em not} to satisfy the desired property, namely the compatibility under the quantum mutations $\Phi^q_{\Delta,\Delta'}$. So the quantum ordering problem for universally Laurent functions cannot be solved just by term-by-term Weyl-ordered products. In the present paper, we will review two solutions to this problem, one by Allegretti and Kim \cite{AK} and the other by Gabella \cite{G}, and finally show that these two answers are the same.

\section{Bonahon-Wong quantum trace}
\label{sec:BW}

\subsection{Tangles and skein algebra}
\label{label:tangles_and_skein_algebra}

The known answers to the quantum ordering problem mentioned in the last section are heavily based on the work of Bonahon and Wong \cite{BW}. Their construction requires us to consider {\em tangles}, which are $1$-dimensional submanifolds in a 3-dimensional manifold subject to certain conditions and with extra data. Depending on the authors and situations, these objects can be defined in different ways, and here we recollect some versions relevant to our paper.

\begin{definition}
$\bullet$ For an oriented surface $\frak{S}$ with boundary, the \ul{\em thickening} of $\frak{S}$ is the $3$-manifold 
$$
\frak{S} \times (-1,1),
$$
whose boundary is $\partial (\frak{S}\times(-1,1)) = (\partial \frak{S}) \times (-1,1)$. We also say that $\frak{S} \times (-1,1)$ is a \ul{\em thickened surface}.

\vs

$\bullet$ The \ul{\em elevation} of a point $(p,t) \in \frak{S} \times (-1,1)$ is the $(-1,1)$-coordinate, namely $t$.

\vs

$\bullet$ If $A$ is a subset of $\frak{S}$, we say that a point $x$ of $\frak{S}\times(-1,1)$ \ul{\em lies over $A$} if $x\in A\times(-1,1)$.
\end{definition}
We will mostly apply this definition to
$$
\frak{S} = \Sigma\setminus\mathcal{P}
$$
for a generalized marked surface $(\Sigma,\mathcal{P})$. Note that connected components $b$ of $\partial \frak{S}$ are interiors of boundary arcs of $(\Sigma,\mathcal{P})$, hence in particular are diffeomorphic to an open interval in $\mathbb{R}$; as mentioned already, we refer to these $b$ as boundary arcs of $(\Sigma,\mathcal{P})$ or of $\frak{S}$. Note that the connected components of $\partial \frak{S} \times (-1,1)$, i.e. the boundary components of $\frak{S}\times(-1,1)$, are the thickenings $b\times(-1,1)$ of boundary arcs $b$ of $\frak{S}$. The {\em projection} to $\frak{S}$ means the usual projection map $\frak{S} \times (-1,1) \to \frak{S}$, or the image of some subset under this map.
\begin{definition}[various tangles]
Let $(\Sigma,\mathcal{P})$ be a generalized marked surface, not necessarily triangulable. Let $\frak{S} = \Sigma\setminus\mathcal{P}$.

\vs

$\bullet$ A \ul{\em tangle} in a thickened surface $\frak{S} \times (-1,1)$ is a $1$-dimensional compact manifold $K$ with boundary properly embedded into $\frak{S} \times (-1,1)$ such that
\begin{enumerate}
\item[\rm (T1)] $K$ is equipped with the choice of a \ul{\em framing} on it, i.e. a continuous choice of a vector in $T_x (\frak{S}\times(-1,1)) \setminus T_x K$ for each point $x$ of $K$, and

\item[\rm (T2)] the framing at each \ul{\em endpoint} of $K$ (a point of $\partial K \subset \partial \frak{S}\times(-1,1)$) is \ul{\em upward vertical}, i.e. is parallel to the $(-1,1)$ factor and points toward $1$ of $(-1,1)$.
\end{enumerate}

$\bullet$ For a tangle $K$ in $\frak{S} \times (-1,1)$, consider the conditions:
\begin{enumerate}
\item[\rm (T3)] for each boundary arc $b$ of $\frak{S}$, the elevations of the endpoints of $K$ lying over $b$ are mutually distinct;

\item[\rm (T4)] for each boundary arc $b$ of $\frak{S}$, no two endpoints of $K$ lying over $b$ project to the same point in $b$.
\end{enumerate}
The tangle $K$ is called a \ul{\em V-tangle} if it satisfies (T3), an \ul{\em H-tangle} if it satisfies (T4), and a \ul{\em VH-tangle} if it satisfies both (T3) and (T4).

\vs

$\bullet$ Define each of \ul{\em V-isotopy}, \ul{\em H-isotopy}, and \ul{\em VH-isotopy} as an isotopy within the class of V-tangles, H-tangles, and VH-tangles, respectively.

\vs

$\bullet$ A V-tangle, H-tangle, or VH-tangle is \ul{\em closed} if the underlying tangle is closed.

\vs

$\bullet$ An oriented version of a tangle, V-tangle, H-tangle and  VH-tangle, as well as the respective isotopies, can be defined.
\end{definition}
In particular, by a tangle we always mean a framed tangle. In the three versions of tangles, the letter V indicates that a V-isotopy class remembers the vertical ordering of endpoints lying over each boundary arc of $\frak{S}$, while H indicates that an H-isotopy class remembers their horizontal ordering. A tangle $K$ in $\frak{S}\times(-1,1)$ can be represented by an embedding map $C \to \frak{S}\times(-1,1)$ of a compact $1$-dimensional manifold $C$ with boundary into $\frak{S}\times(-1,1)$, together with the data of framing. 
The tangles are often studied through their projections in the surface $\frak{S}$, called {\em tangle diagrams}.
\begin{definition}
\label{def:tangle_diagram}
$\bullet$ A \ul{\em tangle diagram} in an oriented surface $\frak{S}$ with boundary is a $1$-dimensional compact manifold $D$ with boundary properly immersed into $\frak{S}$ so that the interior of $D$ lies in the interior of $\frak{S}$, the boundary $\partial D$ of $D$ lies in the boundary $\partial \frak{S}$, and the only possible singularities are transverse double self-intersections lying in the interior $\frak{S} \setminus\partial \frak{S}$ of $\frak{S}$, together with the data:
\begin{enumerate}
\item[\rm (TD1)] for each self-intersection of $D$, called a \ul{\em crossing}, the over- or under-passing information, i.e. an ordering of the two small pieces (``strands") of $D$ forming the crossing (one is under and the other is over).
\end{enumerate}
In a picture, for each crossing $x$, the underpassing part of $D$ is drawn as broken lines, i.e. a small neighborhood of $x$ in this part is deleted. The elements of $\partial D$ are called the \ul{\em endpoints} of $D$.

\vs

An \ul{\em isotopy} of tangle diagrams is an isotopy within the class of tangle diagrams, preserving the over- and under-passing information of the crossings.

\vs

$\bullet$ A \ul{\em boundary-ordering} on a tangle diagram $D$ is the following data:
\begin{enumerate}
\item[\rm (TD2)] for each boundary arc $b$ of $\frak{S}$, the choice of a (total) ordering on the set $\partial_b D := \partial D \cap b$, called a \ul{\em vertical ordering} on $\partial_b D$; if $x\in \partial_b D$ is higher than $y\in \partial_b D$ in this ordering, we write $x\succ y$.
\end{enumerate}
A tangle diagram $D$ equipped with a boundary-ordering is called a \ul{\em boundary-ordered tangle diagram}. When the boundary-ordering is clear from the context, it is denoted just by $D$.

\vs

An \ul{\em isotopy} of boundary-ordered tangle diagrams is an isotopy within the class of boundary-ordered tangle diagrams, preserving the over- and under-passing information of the crossings and the boundary-ordering.

\vs

$\bullet$ An \ul{\em oriented tangle diagram}, a \ul{\em boundary-ordered oriented tangle diagram} and their respective isotopies are defined analogously.
\end{definition}
For a tangle $K$ in $\frak{S} \times (-1,1)$, we consider the following conditions in order to represent it via its projection in $\frak{S}$:
\begin{enumerate}
\item[\rm (P1)] {\em the framing at {\em every} point of $K$ is upward vertical, and}

\item[\rm (P2)] {\em the projection $D$ of $K$ via the projection $\frak{S} \times (-1,1) \to \frak{S}$ is a tangle diagram, where the over- and under-passing information of crossings are induced by the elevations of the points of $K$}.
\end{enumerate}
It is known that each of a V-tangle, H-tangle, or VH-tangle can be isotoped within the respective class to another one satisfying (P1) and (P2). For each case, the over- and under-passing information (TD1) for each crossing is naturally induced from $K$. For a V-tangle or a VH-tangle, the vertical ordering (TD2) on $\partial_b D$ for each boundary arc $b$ of $\frak{S}$ is induced by the ordering on $\partial K \cap (b\times (-1,1))$ coming from elevations; hence the name {\em vertical} ordering. We find it convenient to define also the notion of a {\em horizontal} ordering on $\partial_b D$:

\begin{definition}
\label{def:horizontal_ordering}
Let $b$ be a boundary arc of a surface $\frak{S}$ (with boundary).

$\bullet$ The choice of an orientation on $b$ induces an ordering on the points of $b$ called the \ul{\em horizontal ordering} on $b$ as follows: a point $x\in b$ is (horizontally) higher than $y\in b$ if the direction from $y$ to $x$ matches the chosen orientation on $b$.

\vs

$\bullet$ Let $D$ be a tangle diagram in $\frak{S}$. The \ul{\em horizontal ordering} on the set $\partial_b D$ with respect to an orientation on $b$ is the ordering induced by the horizontal ordering on $b$.
\end{definition}
If $\frak{S}$ is an oriented surface and no particular orientation on $b$ is chosen, we use the boundary-orientation on $b$ (Def.\ref{def:surface}) coming from the surface orientation on $\frak{S}$.

\vs

Note that, once (P1) and (P2) are satisfied, a V-tangle or a VH-tangle in $\frak{S}\times(-1,1)$ yields a boundary-ordered tangle diagram in $\frak{S}$ via projection, while an H-tangle yields a tangle diagram. We say that such a diagram is a diagram \ul{\em of} the respective tangle. Conversely, from an isotopy class of one of these three kinds of tangle diagrams one can recover a unique isotopy class of a respective kind of tangle in $\frak{S}\times (-1,1)$. However, nonisotopic tangle diagrams may yield isotopic tangles, if they differ by certain ``moves".

\begin{figure}[htbp!]
\hspace{-5mm}
\begin{subfigure}[b]{0.22\textwidth}
\hspace{5mm} 
\begingroup%
  \makeatletter%
  \providecommand\color[2][]{%
    \errmessage{(Inkscape) Color is used for the text in Inkscape, but the package 'color.sty' is not loaded}%
    \renewcommand\color[2][]{}%
  }%
  \providecommand\transparent[1]{%
    \errmessage{(Inkscape) Transparency is used (non-zero) for the text in Inkscape, but the package 'transparent.sty' is not loaded}%
    \renewcommand\transparent[1]{}%
  }%
  \providecommand\rotatebox[2]{#2}%
  \newcommand*\fsize{\dimexpr\f@size pt\relax}%
  \newcommand*\lineheight[1]{\fontsize{\fsize}{#1\fsize}\selectfont}%
  \ifx\svgwidth\undefined%
    \setlength{\unitlength}{99.21259843bp}%
    \ifx\svgscale\undefined%
      \relax%
    \else%
      \setlength{\unitlength}{\unitlength * \real{\svgscale}}%
    \fi%
  \else%
    \setlength{\unitlength}{\svgwidth}%
  \fi%
  \global\let\svgwidth\undefined%
  \global\let\svgscale\undefined%
  \makeatother%
  \begin{picture}(1,0.71428571)%
    \lineheight{1}%
    \setlength\tabcolsep{0pt}%
    \put(0,0){\includegraphics[width=\unitlength,page=1]{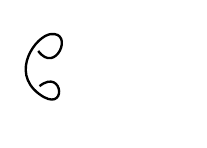}}%
    \put(3.85714286,-0.46524179){\color[rgb]{0,0,0}\makebox(0,0)[lt]{\begin{minipage}{1.54285714\unitlength}\raggedright \end{minipage}}}%
    \put(3.18571429,-0.35095607){\color[rgb]{0,0,0}\makebox(0,0)[lt]{\begin{minipage}{1.08571429\unitlength}\raggedright \end{minipage}}}%
    \put(0.45970355,0.34517114){\color[rgb]{0,0,0}\makebox(0,0)[lt]{\lineheight{1.25}\smash{\begin{tabular}[t]{l}$\leftrightarrow$\end{tabular}}}}%
    \put(0,0){\includegraphics[width=\unitlength,page=2]{moves_I_new.pdf}}%
  \end{picture}%
\endgroup%

\caption{framed Reidemeister\\ \hspace*{3mm} \, move I}
\label{subfig:move_I}
\end{subfigure}
\hfill
\hspace{-3mm}
\begin{subfigure}[b]{0.23\textwidth}
\hspace{4mm} 
\begingroup%
  \makeatletter%
  \providecommand\color[2][]{%
    \errmessage{(Inkscape) Color is used for the text in Inkscape, but the package 'color.sty' is not loaded}%
    \renewcommand\color[2][]{}%
  }%
  \providecommand\transparent[1]{%
    \errmessage{(Inkscape) Transparency is used (non-zero) for the text in Inkscape, but the package 'transparent.sty' is not loaded}%
    \renewcommand\transparent[1]{}%
  }%
  \providecommand\rotatebox[2]{#2}%
  \newcommand*\fsize{\dimexpr\f@size pt\relax}%
  \newcommand*\lineheight[1]{\fontsize{\fsize}{#1\fsize}\selectfont}%
  \ifx\svgwidth\undefined%
    \setlength{\unitlength}{99.21259843bp}%
    \ifx\svgscale\undefined%
      \relax%
    \else%
      \setlength{\unitlength}{\unitlength * \real{\svgscale}}%
    \fi%
  \else%
    \setlength{\unitlength}{\svgwidth}%
  \fi%
  \global\let\svgwidth\undefined%
  \global\let\svgscale\undefined%
  \makeatother%
  \begin{picture}(1,0.71428571)%
    \lineheight{1}%
    \setlength\tabcolsep{0pt}%
    \put(0,0){\includegraphics[width=\unitlength,page=1]{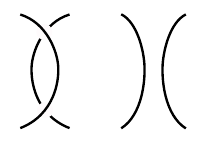}}%
    \put(3.85714286,-0.46524179){\color[rgb]{0,0,0}\makebox(0,0)[lt]{\begin{minipage}{1.54285714\unitlength}\raggedright \end{minipage}}}%
    \put(3.18571429,-0.35095607){\color[rgb]{0,0,0}\makebox(0,0)[lt]{\begin{minipage}{1.08571429\unitlength}\raggedright \end{minipage}}}%
    \put(0.38410841,0.34517114){\color[rgb]{0,0,0}\makebox(0,0)[lt]{\lineheight{1.25}\smash{\begin{tabular}[t]{l}$\leftrightarrow$\end{tabular}}}}%
  \end{picture}%
\endgroup%

\caption{framed Reidemeister\\ \hspace*{4,5mm} move II}
\label{subfig:move_II}
\end{subfigure}
\hfill
\begin{subfigure}[b]{0.23\textwidth}
\hspace{4mm} 
\begingroup%
  \makeatletter%
  \providecommand\color[2][]{%
    \errmessage{(Inkscape) Color is used for the text in Inkscape, but the package 'color.sty' is not loaded}%
    \renewcommand\color[2][]{}%
  }%
  \providecommand\transparent[1]{%
    \errmessage{(Inkscape) Transparency is used (non-zero) for the text in Inkscape, but the package 'transparent.sty' is not loaded}%
    \renewcommand\transparent[1]{}%
  }%
  \providecommand\rotatebox[2]{#2}%
  \newcommand*\fsize{\dimexpr\f@size pt\relax}%
  \newcommand*\lineheight[1]{\fontsize{\fsize}{#1\fsize}\selectfont}%
  \ifx\svgwidth\undefined%
    \setlength{\unitlength}{107.71653543bp}%
    \ifx\svgscale\undefined%
      \relax%
    \else%
      \setlength{\unitlength}{\unitlength * \real{\svgscale}}%
    \fi%
  \else%
    \setlength{\unitlength}{\svgwidth}%
  \fi%
  \global\let\svgwidth\undefined%
  \global\let\svgscale\undefined%
  \makeatother%
  \begin{picture}(1,0.65789474)%
    \lineheight{1}%
    \setlength\tabcolsep{0pt}%
    \put(0,0){\includegraphics[width=\unitlength,page=1]{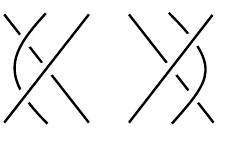}}%
    \put(0.39662231,0.32902216){\color[rgb]{0,0,0}\makebox(0,0)[lt]{\lineheight{1.25}\smash{\begin{tabular}[t]{l}$\leftrightarrow$\end{tabular}}}}%
  \end{picture}%
\endgroup%

\caption{framed Reidemeister\\ \hspace*{4,5mm} move III}
\label{subfig:move_III}
\end{subfigure}
\hfill
\begin{subfigure}[b]{0.23\textwidth}
\hspace{4mm} 
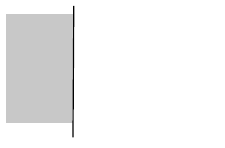
\caption{boundary exchange\\ \hspace*{5mm} move}
\label{subfig:boundary_exchange_move}
\end{subfigure}      
\hspace{-2mm}
\vspace{-2mm}
\caption{Elementary moves for tangles}
\vspace{-4mm}
\label{fig:elementary_moves_for_tangles}
\end{figure}

\begin{enumerate}
\item[\rm (M1)] {\em framed Reidemeister moves I, II and III, as in Fig.\ref{fig:elementary_moves_for_tangles}(I), (II) and (III); in each case, the left tangle diagram and the right tangle diagram differ only in a small disc where they differ as in the picture;}

\item[\rm (M2)] {\em boundary-exchange move as in Fig.\ref{fig:elementary_moves_for_tangles}(IV), where the shaded region indicates the inside of the surface, the vertical thin line part of a boundary arc $b$, the left tangle diagram $D_1$ and the right tangle diagram $D_2$ differ only in a small half-disc near the boundary as in the picture, and the vertical orderings on $\partial_b D_1$ and $\partial_b D_2$ match each other with respect to the natural bijection $\partial_b D_1 \leftrightarrow \partial_b D_2$ sending $x$ to $x$ and $y$ to $y$.}
\end{enumerate}

\begin{proposition}[correspondence between tangle diagrams and tangles {\cite{O}}] Let $\frak{S}$ be an oriented surface (with boundary).
\label{prop:tangle_diagrams_and_tangles}

\begin{enumerate}
\item[\rm (1)] Tangle diagrams in $\frak{S}$ yield isotopic tangles in $\frak{S}\times(-1,1)$ if and only if they are related by a sequence of isotopies of tangle diagrams and moves from (M1) and (M2).

\item[\rm (2)] Tangle diagrams in $\frak{S}$ yield isotopic H-tangles in $\frak{S}\times(-1,1)$ if and only if they are related by a sequence of isotopies of tangle diagrams and moves from (M1).

\item[\rm (3)] Boundary-ordered tangle diagrams in $\frak{S}$ yield isotopic V-tangles in $\frak{S}\times(-1,1)$ if and only if they are related by a sequence of isotopies of tangle diagrams and moves from (M1) and (M2).

\item[\rm (4)] Boundary-ordered tangle diagrams in $\frak{S}$ yield isotopic VH-tangles in $\frak{S}\times(-1,1)$ if and only if they are related by a sequence of isotopies of tangle diagrams and moves from (M1).
\end{enumerate}
\end{proposition}
The statements for oriented versions hold analogously.

\vs

We finally define the stated skein algebra through the following definitions.

\begin{definition}
\label{def:stated_tangles}
Let $(\Sigma,\mathcal{P})$ be a (not necessarily triangulable) generalized marked surface. Let $\frak{S} = \Sigma\setminus\mathcal{P}$.
   
\vs

$\bullet$ A \ul{\em state} on a tangle $K$ in $\frak{S} \times (-1,1)$ is an assignment of a sign to each endpoint of $K$, i.e. a map $s : \partial K \to \{+,-\}$. A \ul{\em stated tangle} is a pair $(K,s)$ of a tangle $K$ and a state $s$ on $K$. The stated versions of various kinds of tangles, such as \ul{\em stated V-tangles}, are defined analogously.

\vs

$\bullet$ A \ul{\em state} on a tangle diagram $D$ is a map $s:\partial D \to \{+,-\}$. A \ul{\em stated tangle diagram} is a pair $(D,s)$ of a tangle diagram $D$ and a state $s$ on $D$. The stated versions of various kinds of tangle diagrams, such as \ul{\em stated boundary-ordered tangle diagrams}, are defined analogously.

\vs

$\bullet$ If $\partial K={\O}$ or $\partial D = {\O}$, the only possible state is denoted by ${\O}$, and the stated tangle $(K,{\O})$ and the stated tangle diagram $(D,{\O})$ are said to be \ul{\em closed}.
\end{definition}

\begin{figure}
\centering 
\scalebox{0.7}{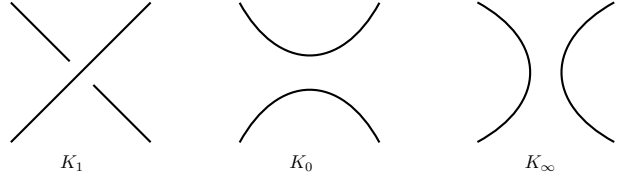}
\vspace{-3mm}
\caption{Kauffman triple of tangles}
\vspace{-5mm}
\label{fig:Kauffman_triple}
\end{figure}

\begin{definition}[the skein algebra; {\cite{BW, Przy, Turaev91}}]
\label{def:skein_algebra}

Let $(\Sigma,\mathcal{P})$ and $\frak{S}$ be as in Def.\ref{def:stated_tangles}. Let $A\in \mathbb{C}^*$.

\vs

The \ul{\em (Kauffman bracket) stated skein algebra} $\mathcal{S}_{\rm s}^A(\Sigma,\mathcal{P})$ is the associative $\mathbb{C}$-algebra defined as follows:

\vs

$\bullet$ The underlying $\mathbb{C}$-vector space is freely spanned by the isotopy classes of all possible stated V-tangles in $\frak{S}\times (-1,1)$, modulo the following relations, where an element of $\mathcal{S}^A_{\rm s}(\Sigma,\mathcal{P})$ represented by the stated tangle $(K,s)$ is denoted by $[K,s]$, and is called a \ul{\em stated skein}:
\begin{enumerate}
\item[\rm (SA1)] (Kauffman bracket skein relation) For each triple of stated V-tangles $(K_1,s_1)$, $(K_0,s_0)$ and $(K_\infty,s_\infty)$ that differ only over a small open disc in the interior of $\frak{S}$ as in Fig.\ref{fig:Kauffman_triple}, one has
$$
[K_1,s_1] = A^{-1}[K_0,s_0] + A [K_\infty,s_\infty]
$$

\item[\rm (SA2)] (the trivial loop relation) For each pair of stated V-tangles $(K,s)$ and $(K',s')$ that differ only over a small disc $U$ such that the intersection of the tangle diagram of $K$ with $U$ is a contractible loop while the intersection of the tangle diagram of $K'$ with $U$ is empty, one has:
$$
[K,s] = -(A^2 + A^{-2}) [K',s'].
$$
\end{enumerate}

\vs

$\bullet$ Let $[K_1,s_1]$ and $[K_2,s_2]$ be stated skeins. Let $[K_1,s_1] = [K_1',s_1']$ and $[K_2,s_2]=[K_2',s_2']$ be such that $K_1' \subset \frak{S} \times (-1,0)$ and $K_2' \subset \frak{S}\times(0,1)$, obtained, for example, by vertically rescaling and translating $K_1$ and $K_2$. Then the product $[K_1,s_1] [K_1,s_2]$ of the two stated skeins in $\mathcal{S}^A_{\rm s}(\Sigma,\mathcal{P})$ is defined as the stated skein $[K,s]$, where the stated V-tangle $(K,s)$ is the union of $(K_1',s_1')$ and $(K_2',s_2')$.

\vs

$\bullet$ A stated skein $[K,{\O}]$ represented by a closed stated V-tangle $(K,{\O})$ is said to be \ul{\em closed}.
\end{definition}
By forgetting the states, one can also define the (Kauffman bracket) skein algebra $\mathcal{S}^A(\Sigma,\mathcal{P})$ spanned by \ul{\em skeins} $[K]$, which goes back to \cite{Turaev91}; the above is one of the several slightly different definitions of what can be called a skein algebra, suited to our purpose; see also \cite{CL,Le, Le18}. It is easy to see that the subspace spanned by all closed skeins is a subalgebra, and often this subalgebra is the only focus of attention (see \cite{Le}); however, we shall see that we need the full stated skein algebra.

\subsection{Square-root quantum Teichm\"uller space}

The skein algebra $\mathcal{S}^A(\Sigma,\mathcal{P})$ is in general noncommutative, and is commutative when $A = \pm 1$. It has been observed that the family of algebras $S^A(\Sigma,\mathcal{P})$, with fixed $(\Sigma,\mathcal{P})$ and $A$ being a varying parameter, yields a version of quantization of the so-called ${\rm SL}_2(\mathbb{C})$-character variety of the surface $\Sigma\setminus \mathcal{P}$ with respect to the Weil-Petersson-Goldman Poisson structure; see \cite{BW} for references. On the other hand, this character variety is closely related to the (enhanced) Teichm\"uller space equipped with the Weil-Petersson Poisson structure, whose corresponding quantum version is established by Chekhov and Fock \cite{F, CF}  and by Kashaev \cite{Kash}, as explained in \S\ref{subsec:quantum_Teichmuller_space}. Thus it is a natural expectation that the skein algebras $\mathcal{S}^A(\Sigma,\mathcal{P})$ should be related to the quantum Teichm\"uller space, and a precise map relating them is what Bonahon and Wong constructed in \cite{BW}. However, the image of the Bonahon-Wong map lies not in the quantum Teichm\"uller space as presented in \S\ref{subsec:quantum_Teichmuller_space}, but in its square-root version, which we now recall.

\vs

Following \cite{BW} and \cite{Hiatt}, the Chekhov-Fock algebra $\mathcal{T}^q_\Delta$ shall be redefined in a slightly different way than in Def.\ref{def:Chekhov-Fock_algebra}. Another quantum parameter $\omega$ is introduced, and two algebras $\mathcal{T}^\omega_\Delta$ and $\mathcal{T}^q_\Delta$ will be constructed, such that $\mathcal{T}^q_\Delta$ is embedded in $\mathcal{T}^\omega_\Delta$.
\begin{definition}[\cite{BW, Hiatt}]
\label{def:square-root_Chekhov-Fock_algebra}
Let $\Delta$ be a triangulation of a generalized marked surface $(\Sigma,\mathcal{P})$.

$\bullet$ Let $t$ be an oriented $\mathcal{P}$-triangle of $\Delta$, and let $e_{t,1}$, $e_{t,2}$ and $e_{t,3}$ be the three sides of $t$ appearing clockwise in this order. In this definition each side is viewed as a map from a closed interval into the triangle, not just as its image, hence these three sides are distinct, whether or not $t$ is self-folded. The \ul{\em triangle algebra} $\mathcal{T}^\omega_t$ associated to $t$ with a nonzero complex parameter $\omega$ is defined as the algebra over $\mathbb{C}$ generated by three elements $\wh{Z}_{t,1}$, $\wh{Z}_{t,2}$ and $\wh{Z}_{t,3}$, and their inverses $\wh{Z}_{t,1}^{-1}$, $\wh{Z}_{t,2}^{-1}$ and $\wh{Z}_{t,3}^{-1}$, with relations
$$
\wh{Z}_{t,1} \wh{Z}_{t,2}  = \omega^2 \wh{Z}_{t,2} \wh{Z}_{t,1}, \qquad
\wh{Z}_{t,2} \wh{Z}_{t,3}  = \omega^2 \wh{Z}_{t,3} \wh{Z}_{t,2}, \qquad
\wh{Z}_{t,3} \wh{Z}_{t,1}  = \omega^2 \wh{Z}_{t,1} \wh{Z}_{t,3},
$$
together with the trivial relations $\wh{Z}_{t,i} \wh{Z}_{t,i}^{-1} = \wh{Z}_{t,i}^{-1} \wh{Z}_{t,i}=1$ for $i=1,2,3$. 

\vs

$\bullet$ The generators $\wh{Z}_{t,1}$, $\wh{Z}_{t,2}$ and $\wh{Z}_{t,3}$ of $\mathcal{T}^\omega_t$ are thought of as being associated to the three sides $e_{t,1}$, $e_{t,2}$ and $e_{t,3}$. In particular, when these sides are named $e$, $f$ and $g$, then $\wh{Z}_{t,e}$, $\wh{Z}_{t,f}$ and $\wh{Z}_{t,g}$ denote $\wh{Z}_{t,1}$, $\wh{Z}_{t,2}$ and $\wh{Z}_{t,3}$ respectively.

\vs

$\bullet$ Write all triangles of $\Delta$ as $t_1,t_2,\ldots,t_m$; so, $m = |\mathcal{F}(\Delta)|$. Consider the tensor product algebra
$$
\textstyle \underset{j=1}{\overset{m}{\bigotimes}} \mathcal{T}^\omega_{t_j} = \mathcal{T}^\omega_{t_1} \otimes \cdots \otimes \mathcal{T}^\omega_{t_m}.
$$
When referring to an element $Z_1 \otimes Z_2 \otimes \cdots \otimes Z_m$ of $\bigotimes_{j=1}^m \mathcal{T}^\omega_{t_j}$, if some factor $Z_j$ equals to $1\in \mathcal{T}^\omega_{t_j}$ then we may omit that factor, provided that one can still see clearly which factor of the element lives in which factor of the tensor product algebra.  The element $1\otimes 1\otimes \cdots \otimes 1$ is denoted by $1$.

\vs

$\bullet$ To each edge $e$ of $\Delta$,  we associate an element $\wh{Z}_e$ of $\bigotimes_{j=1}^m \mathcal{T}^\omega_{t_j}$ defined as:
\begin{enumerate}
\item[\rm (1)] If $e$ is an internal edge that is not a self-folded edge, then $e$ is a side of two distinct triangles $t_j$ and $t_k$; in this case, let $\wh{Z}_e := \wh{Z}_{t_j,e} \otimes \wh{Z}_{t_k,e}$.

\item[\rm (2)] If $e$ is a self-folded edge of a triangle $t_j$, let $e_1$, $e_2$ and $f$ be the three sides of $t_j$ appearing clockwise in this order with images of $e_1$ and $e_2$ coinciding hence forming the self-folded edge; in this case, let $\wh{Z}_e := \omega^{-1} \wh{Z}_{t_j,e_1} \wh{Z}_{t_j,e_2} = \omega \wh{Z}_{t_j,e_2} \wh{Z}_{t_j,e_1}$.

\item[\rm (3)] If $e$ is a boundary edge, and belongs to a triangle $t_j$, then let $\wh{Z}_e := \wh{Z}_{t_j,e}$.
\end{enumerate}
We also define an element $\wh{X}_e$ of $\bigotimes_{j=1}^m \mathcal{T}^\omega_{t_j}$ as
$$
\wh{X}_e := \wh{Z}_e^2.
$$

$\bullet$ Let
$$
q = \omega^4.
$$
Inside $\bigotimes_{j=1}^m \mathcal{T}^\omega_{t_j}$, define 
\begin{align}
\label{eq:BW_CF_algebras}
\left. \begin{array}{l}
\mathcal{T}^\omega_\Delta := \mbox{the subalgebra of $\bigotimes_{j=1}^m \mathcal{T}^\omega_{t_j}$ generated by $\{ \wh{Z}_e, \wh{Z}_e^{-1} \, |\,  e\in \Delta\}$}, \\
\mathcal{T}^q_\Delta := \mbox{the subalgebra of $\bigotimes_{j=1}^m \mathcal{T}^\omega_{t_j}$ generated by $\{ \wh{X}_e, \wh{X}_e^{-1} \, |\,  e\in \Delta\}$,}
\end{array}
\right.
\end{align}
by a slight abuse of notation. Both are called the \ul{\em Chekhov-Fock algebras} associated to $\Delta$. 

\end{definition}
It is easy to observe that the generators of the Chekhov-Fock algebras satisfy the relations
$$
\wh{Z}_e \wh{Z}_f = \omega^{2\varepsilon_{ef}} \wh{Z}_f \wh{Z}_e, \qquad \forall e,f\in {\Delta},
$$
and hence
$$
\wh{X}_e \wh{X}_f = q^{2\varepsilon_{ef}} \wh{X}_f \wh{X}_e, \qquad \forall e,f\in {\Delta}.
$$
Keep in mind that we have injective algebra homomorphism
$$
\mathcal{T}^q_\Delta \hookrightarrow \mathcal{T}^\omega_\Delta \quad : \quad \wh{X}_e \mapsto \wh{Z}_e^2, \qquad \forall e \in {\Delta}.
$$
From now on, we will only use Def.\ref{def:square-root_Chekhov-Fock_algebra}, instead of Def.\ref{def:Chekhov-Fock_algebra}.

\vs

In order to obtain quantum coordinate change maps for the square-root generators as rational formulas, we need to consider:

\begin{definition}[\cite{BW, Hiatt}]
\label{def:CF_balanced_square-root_algebra}
$\bullet$ Denote the generators of the Chekhov-Fock algebra $\mathcal{T}^\omega_\Delta$ associated to the edges of $\Delta$ labeled by $e_1,e_2,\ldots,e_n$ by $\wh{Z}_{e_1},\wh{Z}_{e_2},\ldots,\wh{Z}_{e_n}$; in particular, $n = |{\Delta}|$. A \ul{\em monomial} in $\mathcal{T}^\omega_\Delta$ is an element $\omega^N \wh{Z}_{e_1}^{k_1} \wh{Z}_{e_2}^{k_2} \cdots \wh{Z}_{e_n}^{k_n}$ for some $N,k_1,k_2,\ldots,k_n \in \mathbb{Z}$. A monomial $\omega^N \wh{Z}_{e_1}^{k_1} \wh{Z}_{e_2}^{k_2} \cdots \wh{Z}_{e_n}^{k_n}$ is said to be \ul{\em balanced} if, for every triangle $t_j$ of $\Delta$, the exponents $k_i$ of the generators $\wh{Z}_{e_i}$ associated to the three sides of $t_j$ add up to an even number, where the exponent $k_i$ for a self-folded edge $e_i$ is counted with multiplicity two in this sum.

\vs

$\bullet$ The \ul{\em Chekhov-Fock (balanced) square-root algebra} $\mathcal{Z}^\omega_\Delta$ is defined as the subspace of the Chekhov-Fock algebra $\mathcal{T}^\omega_\Delta$ spanned by all balanced monomials of $\mathcal{T}^\omega_\Delta$.

$\bullet$ The \ul{\em balanced square-root fraction algebra}, denoted by $\wh{\rm Frac}(\mathcal{Z}^\omega_\Delta)$, is defined as the subalgebra of ${\rm Frac}(\mathcal{T}^\omega_\Delta)$ consisting of elements that can be written as $P Q^{-1}$, with $P \in \mathcal{Z}^\omega_\Delta$ and $Q\in \mathcal{T}^q_\Delta \subset \mathcal{T}^\omega_\Delta$.
\end{definition}

As mentioned in \cite{BW}, one can check that ${\rm Frac}(\mathcal{T}^q_\Delta)$ is naturally contained in $\wh{\rm Frac}(\mathcal{Z}^\omega_\Delta)$. Hiatt \cite{Hiatt} constructed quantum coordinate change maps for the (balanced) square-root algebras extending $\Phi^q_{\Delta,\Delta'}$ for the usual Chekhov-Fock algebras, in the following sense:

\begin{proposition}[Hiatt \cite{Hiatt}, see also \cite{BW}]
\label{prop:square-root_skew-field_isomorphisms}
Let $(\Sigma,\mathcal{P})$ be a triangulable generalized marked surface. Let $q$ and $\omega$ be nonzero complex numbers such that $q = \omega^4$. There exists an algebra isomorphism
$$
\Theta^\omega_{\Delta,\Delta'} ~ : ~ \wh{\rm Frac}(\mathcal{Z}^\omega_{\Delta'}) \to \wh{\rm Frac}(\mathcal{Z}^\omega_{\Delta})
$$
associated to each pair $\Delta,\Delta'$ of triangulations of $(\Sigma,\mathcal{P})$, satisfying:
\begin{enumerate}
\item[\rm (1)] the restriction of $\Theta^\omega_{\Delta,\Delta'}$ to ${\rm Frac}(\mathcal{T}^q_\Delta)$ coincides with the map $\Phi^q_{\Delta,\Delta'}$ from Prop.\ref{prop:skew-field_isomorphisms}, and

\item[\rm (2)] the consistency relation
$$
\Theta^\omega_{\Delta,\Delta''} = \Theta^\omega_{\Delta,\Delta'} \circ \Theta^\omega_{\Delta',\Delta''}
$$
holds for each triple $\Delta,\Delta',\Delta''$ of triangulations.
\end{enumerate}

\end{proposition}

\begin{remark}
The algebra $\mathcal{Z}^\omega_\Delta$ is a quantum torus \cite[Cor.13]{BW2}, and hence admits a skew-field of fractions ${\rm Frac}(\mathcal{Z}^\omega_\Delta)$. It can be shown that ${\rm Frac}(\mathcal{Z}^\omega_\Delta)$ coincides with $\wh{\rm Frac}(\mathcal{Z}^\omega_\Delta)$, simplifying the situation. Here is a sketch of proof. Viewing $\wh{\rm Frac}(\mathcal{Z}^\omega_\Delta)$ as a module over ${\rm Frac}(\mathcal{T}^q_\Delta)$, it is a free module, as any module over a skew-field is free, and one can easily see that this module is of finite rank $d$. For any nonzero element $x$ of $\wh{\rm Frac}(\mathcal{Z}^\omega_\Delta)$, note that $1,x,x^2,\ldots,x^d$ are linearly dependent over ${\rm Frac}(\mathcal{T}^q_\Delta)$, hence $\sum_{i=1}^d c_i x^i=0$ for some $c_0,\ldots,c_d \in {\rm Frac}(\mathcal{T}^q_\Delta)$ that are not all zero. Since $x$ is not a zero divisor, $c_0$ must be nonzero, and hence can be assumed to be $1$. It follows that $1 = -(c_1+c_2x + \cdots +c_dx^{d-1})x$, showing that $x$ is invertible. So $\wh{\rm Frac}(\mathcal{Z}^\omega_\Delta)$ is a skew-field, and hence coincides with the skew-field of fractions of $\mathcal{Z}^\omega_\Delta$.
\end{remark}
As pointed out in \cite{BW}, this map $\Theta^\omega_{\Delta,\Delta'}$ is much better understood in terms of operators on Hilbert spaces. This operator-theoretic viewpoint will be made clear and explicit in an upcoming work \cite{KS} (see \cite{Son} for an attempt).

\vs

Similarly as in the case of ${\rm Frac}(\mathcal{T}^q_\Delta)$, by the \ul{\em square-root quantum Teichm\"uller space} we mean the quotient of the disjoint union of all $\wh{\rm Frac}(\mathcal{Z}^\omega_\Delta)$ by the equivalence relation given by the identifications $\Theta^\omega_{\Delta,\Delta'}$. Also, as an analog of Def.\ref{def:bf_L_q}, the ring of \ul{\em square-root quantum regular functions} can be defined as
\begin{align}
\label{eq:bf_L_omega}
{\bf L}^\omega_\Delta := \bigcap_{\Delta'} \Theta^\omega_{\Delta,\Delta'}(\mathcal{Z}^\omega_{\Delta'}) \quad \subset \quad \mathcal{Z}^\omega_\Delta \quad \subset \quad \wh{{\rm Frac}}(\mathcal{Z}^\omega_\Delta) \quad \subset \quad {\rm Frac}(\mathcal{T}^\omega_\Delta),
\end{align}
by a slight abuse of notation with the previously defined symbol ${\bf L}^q_\Delta$. As $\Theta^\omega_{\Delta,\Delta'}$ induces a natural isomorphism from ${\bf L}^\omega_{\Delta'}$ to ${\bf L}^\omega_\Delta$, we may denote ${\bf L}^\omega_\Delta$ for all $\Delta$ collectively by ${\bf L}^\omega = {\bf L}^\omega_{(\Sigma,\mathcal{P})}$, or understand this situation as ${\bf L}^\omega = \bigsqcup_\Delta {\bf L}^\omega_\Delta / \hspace{-1mm}\sim$, where $\sim$ is the equivalence relation coming from the maps $\Theta^\omega_{\Delta,\Delta'}$.

\subsection{Quantum trace map}
\label{subsec:quantum_trace_map}

If we restrict our attention only to triangulable generalized marked surfaces $(\Sigma,\mathcal{P})$ with empty boundary, i.e. punctured surfaces, we are ready to state the result of Bonahon and Wong, namely a map from the skein algebra of this surface to the square-root quantum Teichm\"uller space. In the meantime, one of the major defining properties of this map is the cutting/gluing property, which is a certain compatibility that holds when cutting the surface (together with a skein) along a $\mathcal{P}$-arc, i.e. an edge of some triangulation. In order to fully reflect this property, it is not just a luxury but rather a must, that we state Bonahon and Wong's result in complete generality for any triangulable generalized marked surface, instead of only for punctured surfaces without boundary.

\vs

Bonahon and Wong \cite{BW} described the process of gluing surfaces and skeins along two boundary arcs of a (not necessarily connected) surface. To conveniently attain uniqueness of such a process, we instead formulate everything in terms of cutting.

\begin{deflem}[cutting construction]
\label{def:cutting}
Let $(\Sigma,\mathcal{P})$ be a generalized marked surface, not necessarily connected nor triangulable. Let $b$ be a $\mathcal{P}$-arc in $\Sigma$, and assume that the interior $\mathring{b}$ of $b$ lies in the interior $\Sigma\setminus\partial\Sigma$ of the surface. 

\vs

$\bullet$ Let $(\Sigma',\mathcal{P}')$ be the unique (up to diffeomorphism) generalized marked surface obtained from $(\Sigma,\mathcal{P})$ by cutting along $b$. In particular, there is a natural (gluing) map $\Sigma' \to \Sigma$, whose restriction to $\mathcal{P}'$ yields a correspondence $\mathcal{P}' \to \mathcal{P}$. The preimage of $b$ under $\Sigma'\to \Sigma$ is the union of two boundary arcs of $(\Sigma',\mathcal{P}')$, denoted by $b_1$ and $b_2$. The gluing map $\Sigma' \to \Sigma$ restricts to a diffeomorphism $\Sigma' \setminus (b_1\cup b_2) \to \Sigma\setminus b$.

\vs

$\bullet$ Suppose that $K$ is a tangle in $(\Sigma\setminus\mathcal{P})\times(-1,1)$ satisfying:
\begin{enumerate}
\item[\rm (C1)] $K$ is transverse to $b\times (-1,1)$, and

\item[\rm (C2)] the framing of $K$ at each point in $K\cap(b\times (-1,1))$ is upward vertical.
\end{enumerate}
The above process of cutting along $b$ uniquely yields a tangle $K'$ in $(\Sigma'\setminus\mathcal{P}')\times(-1,1)$, equipped with a map $K' \to K$. The number of preimages of $x\in K$ under $K' \to K$ is two if $x\in b$, and is one otherwise.

\vs

$\bullet$ Suppose further that $(\Sigma,\mathcal{P})$ is triangulable, and that $\Delta$ is a triangulation of $(\Sigma,\mathcal{P})$ such that $b$ is an internal edge of $\Delta$, i.e. $b\in \mathring{\Delta}$. The process of cutting along $b$ uniquely yields a triangulation $\Delta'$ of $(\Sigma',\mathcal{P}')$.

\vs

$\bullet$ In each of the above three situations, we say that the new object (for $(\Sigma',\mathcal{P}')$) is obtained from the original by \ul{\em cutting along $b$}.

\end{deflem}

\begin{lemma}[\cite{BW}]
\label{lem:injection_of_sq_algebras_for_gluing}
Suppose the situation of Def.\ref{def:cutting}, and $\omega \in \mathbb{C}^*$. The triangles of $\Delta$ are naturally in one-to-one correspondence with those of $\Delta'$, and the sides of each triangle of $\Delta$ are naturally in one-to-one correspondence with those of the corresponding triangle of $\Delta'$. Hence we have a natural isomorphism $\bigotimes_{t \in \mathcal{F}(\Delta)} \mathcal{T}^\omega_{t} \to \bigotimes_{t' \in \mathcal{F}(\Delta')} \mathcal{T}^\omega_{t'}$ of algebras, and this restricts to an injective algebra homomorphism
\begin{align}
\label{eq:injection_of_sq_algebras_for_gluing}
\mathcal{Z}^\omega_{\Delta} \to \mathcal{Z}^\omega_{\Delta'}
\end{align}
between the Chekhov-Fock (balanced) square-root algebras.
\end{lemma}

\begin{definition}[\cite{BW}]
\label{def:compatible_stated_skeins}
Let $(\Sigma,\mathcal{P})$ be a generalized marked surface, not necessarily triangulable. Let $A \in \mathbb{C}^*$, and let $\mathcal{S}^A_{\rm s}(\Sigma,\mathcal{P})$ be the stated skein algebra defined in Def.\ref{def:skein_algebra}.

\vs

$\bullet$ Suppose that $(\Sigma',\mathcal{P}')$ is obtained from $(\Sigma,\mathcal{P})$ by cutting along a $\mathcal{P}$-arc $b$ of $\Sigma$, and that a tangle $K'$ in $(\Sigma'\setminus\mathcal{P}')\times(-1,1)$ is obtained from a tangle $K$ in $(\Sigma\setminus\mathcal{P})\times(-1,1)$ through this cutting process, as described in Def.\ref{def:cutting}. For each $x' \in \partial K'$, denote by $x\in K$ the image of $x'$ under the map $K' \to K$ in Def.\ref{def:cutting}. We say that the states $s' : \partial K' \to \{+,-\}$ and $s:\partial K \to\{+,-\}$ for these tangles are \ul{\em compatible} if the following hold: 
\begin{enumerate}
\item[\rm (CS1)] If $x \in \partial K$, then $s'(x') = s(x)$;

\item[\rm (CS2)] If $x \notin \partial K$, i.e. $x\in b$, so that the preimage of $x$ under $K'\to K$ is $\{x',x''\}$, then $s'(x') = s'(x'')$.
\end{enumerate}

\end{definition}

We can finally state the result of Bonahon and Wong in full generality.

\begin{figure}
\centering 
\begin{subfigure}[b]{0.4\textwidth}
\hspace{13mm} \scalebox{0.8}{
\begingroup%
  \makeatletter%
  \providecommand\color[2][]{%
    \errmessage{(Inkscape) Color is used for the text in Inkscape, but the package 'color.sty' is not loaded}%
    \renewcommand\color[2][]{}%
  }%
  \providecommand\transparent[1]{%
    \errmessage{(Inkscape) Transparency is used (non-zero) for the text in Inkscape, but the package 'transparent.sty' is not loaded}%
    \renewcommand\transparent[1]{}%
  }%
  \providecommand\rotatebox[2]{#2}%
  \newcommand*\fsize{\dimexpr\f@size pt\relax}%
  \newcommand*\lineheight[1]{\fontsize{\fsize}{#1\fsize}\selectfont}%
  \ifx\svgwidth\undefined%
    \setlength{\unitlength}{113.38582677bp}%
    \ifx\svgscale\undefined%
      \relax%
    \else%
      \setlength{\unitlength}{\unitlength * \real{\svgscale}}%
    \fi%
  \else%
    \setlength{\unitlength}{\svgwidth}%
  \fi%
  \global\let\svgwidth\undefined%
  \global\let\svgscale\undefined%
  \makeatother%
  \begin{picture}(1,1)%
    \lineheight{1}%
    \setlength\tabcolsep{0pt}%
    \put(0,0){\includegraphics[width=\unitlength,page=1]{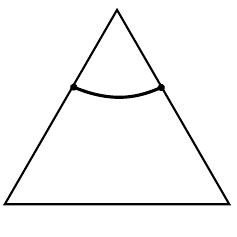}}%
    \put(0.21079075,0.63504734){\color[rgb]{0,0,0}\makebox(0,0)[lt]{\lineheight{1.25}\smash{\begin{tabular}[t]{l}$x_1$\end{tabular}}}}%
    \put(0.71184502,0.63008637){\color[rgb]{0,0,0}\makebox(0,0)[lt]{\lineheight{1.25}\smash{\begin{tabular}[t]{l}$x_2$\end{tabular}}}}%
    \put(0,0){\includegraphics[width=\unitlength,page=2]{two_a_new2.pdf}}%
  \end{picture}%
\endgroup%
}
\vspace{-3mm}
\caption*{(a)}
\end{subfigure}
~
\begin{subfigure}[b]{0.4\textwidth}
\hspace{15mm}  \scalebox{0.8}{
\begingroup%
  \makeatletter%
  \providecommand\color[2][]{%
    \errmessage{(Inkscape) Color is used for the text in Inkscape, but the package 'color.sty' is not loaded}%
    \renewcommand\color[2][]{}%
  }%
  \providecommand\transparent[1]{%
    \errmessage{(Inkscape) Transparency is used (non-zero) for the text in Inkscape, but the package 'transparent.sty' is not loaded}%
    \renewcommand\transparent[1]{}%
  }%
  \providecommand\rotatebox[2]{#2}%
  \newcommand*\fsize{\dimexpr\f@size pt\relax}%
  \newcommand*\lineheight[1]{\fontsize{\fsize}{#1\fsize}\selectfont}%
  \ifx\svgwidth\undefined%
    \setlength{\unitlength}{113.38582677bp}%
    \ifx\svgscale\undefined%
      \relax%
    \else%
      \setlength{\unitlength}{\unitlength * \real{\svgscale}}%
    \fi%
  \else%
    \setlength{\unitlength}{\svgwidth}%
  \fi%
  \global\let\svgwidth\undefined%
  \global\let\svgscale\undefined%
  \makeatother%
  \begin{picture}(1,1)%
    \lineheight{1}%
    \setlength\tabcolsep{0pt}%
    \put(0,0){\includegraphics[width=\unitlength,page=1]{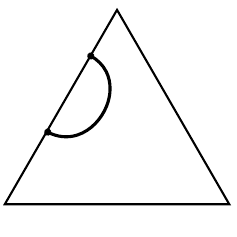}}%
    \put(0.2818975,0.78222184){\color[rgb]{0,0,0}\makebox(0,0)[lt]{\lineheight{1.25}\smash{\begin{tabular}[t]{l}$x_1$\end{tabular}}}}%
    \put(0.08842052,0.45479994){\color[rgb]{0,0,0}\makebox(0,0)[lt]{\lineheight{1.25}\smash{\begin{tabular}[t]{l}$x_2$\end{tabular}}}}%
    \put(0.2264047,0.57954812){\color[rgb]{0,0,0}\rotatebox{57.591302}{\makebox(0,0)[lt]{\lineheight{1.25}\smash{\begin{tabular}[t]{l}$\prec$\end{tabular}}}}}%
    \put(0,0){\includegraphics[width=\unitlength,page=2]{two_b_new2.pdf}}%
  \end{picture}%
\endgroup%
}
\vspace{-3mm}
\caption*{(b)}
\end{subfigure}
\vspace{-2mm}
\caption{Boundary-ordered tangle diagrams for elementary cases of a skein in a triangle (arrows on the sides indicate the boundary-orientation (Def.\ref{def:surface}))}
\vspace{-3mm}
\label{fig:elementary_skeins_in_the_triangle}
\end{figure}

\begin{proposition}[quantum trace map; the main theorem of \cite{BW}, see also \cite{CL, Le, Le18}]
\label{prop:BW_full}
Let $A$ and $\omega$ be nonzero complex numbers such that
$$
A = \omega^{-2}.
$$
Then there is a unique family of algebra homomorphisms
\begin{align}
\label{eq:BW_map}
{\rm Tr}^\omega_\Delta = {\rm Tr}^\omega_{(\Sigma,\mathcal{P});\Delta} ~ : ~ \mathcal{S}^A_{\rm s}(\Sigma,\mathcal{P}) \to \mathcal{Z}^\omega_\Delta
\end{align}
from the stated skein algebra to the Chekhov-Fock (balanced) square-root algebra, defined for each triangulable generalized marked surface $(\Sigma,\mathcal{P})$ and each triangulation $\Delta$ of $(\Sigma,\mathcal{P})$, satisfying:

\begin{enumerate}
\item[\rm (1)] (Cutting Property\footnote{Bonahon and Wong \cite{BW} call this `State Sum Property'}) Suppose that $(\Sigma',\mathcal{P}')$, $\Delta'$ and $K'$ are related to $(\Sigma,\mathcal{P})$, $\Delta$ and $K$ as in Def.\ref{def:cutting}, i.e. the former are obtained by cutting the latter along an internal edge of $\Delta$. Assume that $K'$ is a V-tangle. Then for each state $s$ for the tangle $K$, we have
\begin{align}
\nonumber
{\rm Tr}^\omega_{(\Sigma,\mathcal{P});\Delta}([K,s]) = \sum_{\mbox{\tiny compatible $s'$}} {\rm Tr}^\omega_{(\Sigma',\mathcal{P}');\Delta'}([K',s']),
\end{align}
where the sum is over all states $s'$ for the tangle $K'$ for the cut surface $(\Sigma',\mathcal{P}')$ that are compatible with $s$ in the sense of Def.\ref{def:compatible_stated_skeins}. The left-hand-side, which is a priori an element of $\mathcal{Z}^\omega_{\Delta}$, is viewed as an element of $\mathcal{Z}^\omega_{\Delta'}$ via the embedding map in eq.\eqref{eq:injection_of_sq_algebras_for_gluing}.

\vs

\item[\rm (2)] (Elementary Cases) Suppose that $(\Sigma,\mathcal{P})$ is a non-self-folded triangle, i.e. $\Sigma$ is diffeomorphic to a closed disc, and $\mathcal{P}$ consists of three marked points on the boundary; let $\Delta$ denote its unique triangulation. If $[K,s] \in \mathcal{S}^A_{\rm s}(\Sigma,\mathcal{P})$ is a stated skein where the boundary-ordered tangle diagram of $K$ is one of the two cases (a) and (b) in Fig.\ref{fig:elementary_skeins_in_the_triangle}, we have:

\vs

\begin{enumerate}
\item[\rm (a)] Suppose that $\wh{Z}_1$ and $\wh{Z}_2$ are the generators of the Chekhov-Fock algebra $\mathcal{T}^\omega_\Delta$ of eq.\eqref{eq:BW_CF_algebras} associated to the edges of $\Delta$ containing the endpoints $x_1$ and $x_2$ in Fig.\ref{fig:elementary_skeins_in_the_triangle}(a), respectively. In particular, $\wh{Z}_1 \wh{Z}_2 = \omega^2 \wh{Z}_2 \wh{Z}_1$. Then,
\begin{align}
\label{eq:BW_triangle_case_a}
{\rm Tr}^\omega_{(\Sigma,\mathcal{P});\Delta} ([K,s])  = \left\{
\begin{array}{ll}
0 & \mbox{if $s(x_1) = -$ and $s(x_2)=+$}, \\
\left[\wh{Z}_1^{s(x_1)} \wh{Z}_2^{s(x_2)}\right]_{\rm Weyl} & \mbox{otherwise}, 
\end{array}
\right.
\end{align}
where a sign $\varepsilon = \pm$ in the exponent means $\pm1$ respectively.

\vs

\item[\rm (b)] For $K$ as in Fig.\ref{fig:elementary_skeins_in_the_triangle}(b), one has
\begin{align*}
{\rm Tr}^\omega_{(\Sigma,\mathcal{P});\Delta} ([K,s]) = \left\{
\begin{array}{ll}
0 & \mbox{if $s(x_1)=s(x_2)$}, \\
-\omega^{-5} & \mbox{if $s(x_1) = +$ and $s(x_2) = -$}, \\
\omega^{-1} & \mbox{if $s(x_1) = -$ and $s(x_2) = +$}.
\end{array}
\right.
\end{align*}

\end{enumerate}

\end{enumerate}
\end{proposition}

For a triangulated generalized marked surface, cutting along all the internal edges of the triangulation yields the disjoint union of non-self-folded triangles. Also, using the skein relations repeatedly, any stated skein $[K,s] \in \mathcal{S}^A_{\rm s}(\Sigma,\mathcal{P})$ can be written as a linear combination of stated skeins whose tangle diagrams have no crossing at all. Along this line, one can show that the properties of ${\rm Tr}^\omega_{(\Sigma,\mathcal{P});\Delta}$ in Prop.\ref{prop:BW_full} completely determine the values of the maps ${\rm Tr}^\omega_{(\Sigma,\mathcal{P});\Delta}$, if the existence of these maps is assumed.

\vs

A crucial property of the Bonahon-Wong quantum trace is its compatibility with the quantum coordinate-change maps.

\begin{proposition}[\cite{BW}]
\label{prop:BW_under_mutations}
Let $\Delta,\Delta'$ be triangulations of a triangulable generalized marked surface $(\Sigma,\mathcal{P})$. Then we have
$$
{\rm Tr}^\omega_\Delta = \Theta^\omega_{\Delta,\Delta'} \circ {\rm Tr}^\omega_{\Delta'},
$$
where $\Theta^\omega_{\Delta,\Delta'} : \wh{\rm Frac}(\mathcal{Z}^\omega_{\Delta'}) \to \wh{\rm Frac}(\mathcal{Z}^\omega_\Delta)$ is the (square-root) quantum coordinate-change isomorphism in Prop.\ref{prop:square-root_skew-field_isomorphisms}. More precisely, for each stated skein $[K,s] \in \mathcal{S}^A_{\rm s}(\Sigma,\mathcal{P})$, with $A = \omega^{-2}$, the balanced Laurent polynomial ${\rm Tr}^\omega_{\Delta'}([K,s]) \in \mathcal{Z}^\omega_{\Delta'}$ is sent by $\Theta^\omega_{\Delta,\Delta'}$ to the balanced Laurent polynomial ${\rm Tr}^\omega_\Delta([K,s]) \in \mathcal{Z}^\omega_\Delta$.
\end{proposition}
In particular, the images of the Bonahon-Wong quantum trace map are balanced Laurent polynomials in the square-root generators for {\em any} triangulation $\Delta$, and hence belong to ${\bf L}^\omega_\Delta$, defined in eq.\eqref{eq:bf_L_omega}.

\vs

We now briefly introduce Allegretti and Kim's idea \cite{AK} of using the Bonahon-Wong quantum trace to obtain a solution to the quantum ordering problem discussed in \S\ref{subsec:deformation_quantization}. Consider a $\mathcal{P}$-knot $\gamma$ in a triangulable generalized marked surface $(\Sigma,\mathcal{P})$,  that is, $\gamma$ is a simple loop in the surface $\frak{S}=\Sigma\setminus\mathcal{P}$. Suppose that $\gamma$ is not a contractible loop. Then $\gamma$ represents a nontrivial element $[\gamma]$ of $\pi_1(\frak{S})$, say, if we choose a basepoint and an orientation. Given any point $\rho$ of the Teichm\"uller space $\mathcal{T}(\Sigma,\mathcal{P})$ defined by eq.\eqref{eq:Teichmuller_space}, one considers the monodromy $\rho([\gamma])$, which is an element of ${\rm PSL}(2,\mathbb{R})$ defined up to conjugation in ${\rm PSL}(2,\mathbb{R})$; then $| {\rm trace}(\rho([\gamma])) |$, the absolute value of the trace of this monodromy, is a well-defined real number. This provides a smooth function $\mathbb{I}(\gamma)$ on $\mathcal{X}^+(\Sigma,\mathcal{P})$, whose value at each point $(\rho,O)$ is defined to be $| {\rm trace}(\rho([\gamma])) |$. It turns out that, for {\em any} chosen ideal triangulation $\Delta$ of $(\Sigma,\mathcal{P})$, the function $\mathbb{I}(\gamma)$ is a (positive-)integer-coefficient Laurent polynomial in the square-roots $X_e^{1/2}$ of the exponentiated shear coordinate functions of edges $e$ of $\Delta$. Based on these functions, which naturally arise geometrically, Fock and Goncharov \cite{FG06} proposed a basis of the ring of all universally Laurent functions, i.e. functions that are Laurent polynomials in the exponentiated shear coordinate functions (resp. in their square-roots) for every ideal triangulation, where this basis is enumerated by ``even integral laminations" (resp. ``integral laminations") on the surface $\frak{S}$, which are multicurves with integer weights subject to certain conditions. For a proof that this basis indeed spans all universally Laurent functions, see \cite{FG06,GHK15,Shen}. One can view a simple loop $\gamma$ as a special example of an integral lamination. In order to construct a quantum version $\wh{\mathbb{I}}^\omega(\gamma)$ that deforms the classical function $\mathbb{I}(\gamma)$, we first lift the curve $\gamma$ living in the surface $\frak{S}$ to a tangle in the 3-dimensional manifold $\frak{S}\times(-1,1)$, the thickening of $\frak{S}$.
\begin{definition}
\label{def:K_gamma}
Let $(\Sigma,\mathcal{P})$ be a generalized marked surface. Let $\gamma$ be a simple closed curve in $\frak{S}=\Sigma\setminus\mathcal{P}$, i.e. a closed $\mathcal{P}$-link in $\Sigma$.

\vs

Denote by $K^c_\gamma$ the tangle in $\frak{S} \times(-1,1)$ obtained as the lift of $\gamma$ at a constant elevation $c\in (-1,1)$ with upward vertical framing everywhere. We call $K^c_\gamma$ a \ul{\em constant-elevation lift} of $\gamma$. We may denote $K^c_\gamma$ by $K_\gamma$ without specifying $c$.
\end{definition}
Note that $K_\gamma$ is well-defined up to isotopy, for a given $\gamma$. An easy observation:
\begin{lemma}
A closed tangle $K$ in $\frak{S} \times (-1,1)$ is isotopic to a tangle whose tangle diagram in $\Sigma$ has no crossing at all if and only if $K$ is isotopic to a constant-elevation lift $K_\gamma$ of a simple closed curve $\gamma$ in $\Sigma\setminus\mathcal{P}$. $\qed$
\end{lemma}

Note that $K_\gamma$ is a closed tangle, hence a closed V-tangle. Choose any ideal triangulation $\Delta$ of $(\Sigma,\mathcal{P})$, and apply the Bonahon-Wong quantum trace map in eq.\eqref{eq:BW_map} to the stated skein $[K_\gamma,{\O}] \in \mathcal{S}^A_{\rm s}(\Sigma,\mathcal{P})$; the resulting element of $\mathcal{Z}^\omega_\Delta$ is the Allegretti-Kim quantum element associated to $\mathbb{I}(\gamma)$:
\begin{align}
\label{eq:AK_solution}
\wh{\mathbb{I}}^\omega_\Delta(\gamma) := {\rm Tr}^\omega_{(\Sigma,\mathcal{P});\Delta}([K_\gamma,{\O}]) \quad \in \quad {\bf L}^\omega_\Delta \quad \subset\quad \mathcal{Z}^\omega_\Delta.
\end{align}
Treatment for more general integral laminations needs several more crucial ideas. See \cite{AK} for this, and also for various favorable properties enjoyed by these quantum elements; see also \cite{CKKO} for an important positivity property. To present a couple of examples of its nice properties, we have the quantum mutation compatibility $\wh{\mathbb{I}}^\omega_\Delta(\gamma)  = \Theta^\omega_{\Delta,\Delta'} (\wh{\mathbb{I}}^\omega_{\Delta'}(\gamma))$, and it is relatively easy to see that $\wh{\mathbb{I}}^\omega_\Delta(\gamma)$ indeed recovers the classical function $\mathbb{I}(\gamma)$ when $\omega=1$; see \cite{BW}. In particular, the assignment $\mathbb{I}(\gamma) \mapsto \wh{\mathbb{I}}^\omega(\gamma) = \wh{\mathbb{I}}^\omega_\Delta(\gamma) = {\rm Tr}^\omega_{(\Sigma,\mathcal{P});\Delta}([K_\gamma,{\O}])$ provides a partial answer to the quantum ordering problem mentioned in \S\ref{subsec:deformation_quantization}; see \cite{AK}  for a full answer, which requires algebraic manipulations including Chebyshev polynomials, as well as certain control of parity of powers of monomials.

\subsection{Biangle quantum trace}
\label{subsec:biangles}

Although Propositions \ref{prop:BW_full} and \ref{prop:BW_under_mutations} completely describe and determine the Bonahon-Wong quantum trace map, they are not very convenient when it comes to actual computation of the values. For any given stated skein $[K,s] \in \mathcal{S}^A_{\rm s}(\Sigma,\mathcal{P})$, there is a more direct algorithm that enables us to compute the quantum trace ${\rm Tr}^\omega_\Delta([K,s])$ called the ``state-sum formula", which we shall recall in \S\ref{subsec:state-sum_model}. As a preliminary step for that formula, we first recall the quantum trace for biangles in the present subsection.

\vs

Recall that a \ul{\em biangle} $B$ is a generalized marked surface $(\Sigma,\mathcal{P})$, with $\Sigma$ diffeomorphic to a closed disc, which in particular has one boundary component, and where $\mathcal{P}$ consists of two marked points on the boundary. As noted in Def.\ref{def:triangulability}, it is not triangulable, hence there is no Bonahon-Wong quantum trace map that Prop.\ref{prop:BW_full} associates to it; in particular, the quantum Teichm\"uller space is not defined. However, its stated skein algebra makes sense, because Definition \ref{def:skein_algebra} applies. Bonahon and Wong \cite{BW} defined and studied a quantum trace map for biangles, separately from Prop.\ref{prop:BW_full}.

\vs

Recall that in our notations, we have $\partial B = \partial \Sigma \setminus \mathcal{P}$, and $\partial B$ is a disjoint union of two boundary arcs of $B$. Let's say that a $\mathcal{P}$-arc in $B$ is an \ul{\em internal arc} if its interior is contained in the interior of $\Sigma$. Recall that each generalized marked surface $(\Sigma,\mathcal{P})$ is equipped with an orientation on $\Sigma$. Hence, by a biangle we automatically mean an oriented biangle. Notice that any two biangles are diffeomorphic.

\vs

\begin{figure}[htbp!]
\hfill
\begin{subfigure}[b]{0.13\textwidth}
\hspace*{-4mm} 
\begingroup%
  \makeatletter%
  \providecommand\color[2][]{%
    \errmessage{(Inkscape) Color is used for the text in Inkscape, but the package 'color.sty' is not loaded}%
    \renewcommand\color[2][]{}%
  }%
  \providecommand\transparent[1]{%
    \errmessage{(Inkscape) Transparency is used (non-zero) for the text in Inkscape, but the package 'transparent.sty' is not loaded}%
    \renewcommand\transparent[1]{}%
  }%
  \providecommand\rotatebox[2]{#2}%
  \newcommand*\fsize{\dimexpr\f@size pt\relax}%
  \newcommand*\lineheight[1]{\fontsize{\fsize}{#1\fsize}\selectfont}%
  \ifx\svgwidth\undefined%
    \setlength{\unitlength}{70.86614173bp}%
    \ifx\svgscale\undefined%
      \relax%
    \else%
      \setlength{\unitlength}{\unitlength * \real{\svgscale}}%
    \fi%
  \else%
    \setlength{\unitlength}{\svgwidth}%
  \fi%
  \global\let\svgwidth\undefined%
  \global\let\svgscale\undefined%
  \makeatother%
  \begin{picture}(1,1.04)%
    \lineheight{1}%
    \setlength\tabcolsep{0pt}%
    \put(0,0){\includegraphics[width=\unitlength,page=1]{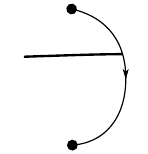}}%
    \put(0.02921547,0.65210342){\color[rgb]{0,0,0}\makebox(0,0)[lt]{\lineheight{1.25}\smash{\begin{tabular}[t]{l}$y$\end{tabular}}}}%
    \put(0,0){\includegraphics[width=\unitlength,page=2]{tangles_unoriented_identity.pdf}}%
    \put(0.87236481,0.66107025){\color[rgb]{0,0,0}\makebox(0,0)[lt]{\lineheight{1.25}\smash{\begin{tabular}[t]{l}$x$\end{tabular}}}}%
  \end{picture}%
\endgroup%

\caption{
\\\hspace*{1mm}}
\end{subfigure}
\hfill
\hspace{2mm}
\begin{subfigure}[b]{0.13\textwidth}
\hspace*{-3mm} 
\begingroup%
  \makeatletter%
  \providecommand\color[2][]{%
    \errmessage{(Inkscape) Color is used for the text in Inkscape, but the package 'color.sty' is not loaded}%
    \renewcommand\color[2][]{}%
  }%
  \providecommand\transparent[1]{%
    \errmessage{(Inkscape) Transparency is used (non-zero) for the text in Inkscape, but the package 'transparent.sty' is not loaded}%
    \renewcommand\transparent[1]{}%
  }%
  \providecommand\rotatebox[2]{#2}%
  \newcommand*\fsize{\dimexpr\f@size pt\relax}%
  \newcommand*\lineheight[1]{\fontsize{\fsize}{#1\fsize}\selectfont}%
  \ifx\svgwidth\undefined%
    \setlength{\unitlength}{70.86614173bp}%
    \ifx\svgscale\undefined%
      \relax%
    \else%
      \setlength{\unitlength}{\unitlength * \real{\svgscale}}%
    \fi%
  \else%
    \setlength{\unitlength}{\svgwidth}%
  \fi%
  \global\let\svgwidth\undefined%
  \global\let\svgscale\undefined%
  \makeatother%
  \begin{picture}(1,1.04)%
    \lineheight{1}%
    \setlength\tabcolsep{0pt}%
    \put(0,0){\includegraphics[width=\unitlength,page=1]{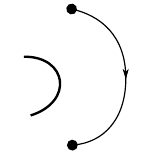}}%
    \put(0.02921547,0.65210342){\color[rgb]{0,0,0}\makebox(0,0)[lt]{\lineheight{1.25}\smash{\begin{tabular}[t]{l}$y_2$\end{tabular}}}}%
    \put(0,0){\includegraphics[width=\unitlength,page=2]{tangles_unoriented_cup.pdf}}%
    \put(0.07276626,0.22511337){\color[rgb]{0,0,0}\makebox(0,0)[lt]{\lineheight{1.25}\smash{\begin{tabular}[t]{l}$y_1$\end{tabular}}}}%
    \put(0.11562454,0.40571621){\color[rgb]{0,0,0}\rotatebox{90.38895}{\makebox(0,0)[lt]{\lineheight{1.25}\smash{\begin{tabular}[t]{l}$\prec$\end{tabular}}}}}%
  \end{picture}%
\endgroup%

\caption{
\\\hspace*{1mm}}
\end{subfigure}
\hfill
\begin{subfigure}[b]{0.15\textwidth}
\hspace*{-2mm} 
\begingroup%
  \makeatletter%
  \providecommand\color[2][]{%
    \errmessage{(Inkscape) Color is used for the text in Inkscape, but the package 'color.sty' is not loaded}%
    \renewcommand\color[2][]{}%
  }%
  \providecommand\transparent[1]{%
    \errmessage{(Inkscape) Transparency is used (non-zero) for the text in Inkscape, but the package 'transparent.sty' is not loaded}%
    \renewcommand\transparent[1]{}%
  }%
  \providecommand\rotatebox[2]{#2}%
  \newcommand*\fsize{\dimexpr\f@size pt\relax}%
  \newcommand*\lineheight[1]{\fontsize{\fsize}{#1\fsize}\selectfont}%
  \ifx\svgwidth\undefined%
    \setlength{\unitlength}{70.86614173bp}%
    \ifx\svgscale\undefined%
      \relax%
    \else%
      \setlength{\unitlength}{\unitlength * \real{\svgscale}}%
    \fi%
  \else%
    \setlength{\unitlength}{\svgwidth}%
  \fi%
  \global\let\svgwidth\undefined%
  \global\let\svgscale\undefined%
  \makeatother%
  \begin{picture}(1,1.04)%
    \lineheight{1}%
    \setlength\tabcolsep{0pt}%
    \put(0,0){\includegraphics[width=\unitlength,page=1]{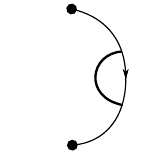}}%
    \put(0.84120651,0.68948082){\color[rgb]{0,0,0}\makebox(0,0)[lt]{\lineheight{1.25}\smash{\begin{tabular}[t]{l}$x_2$\end{tabular}}}}%
    \put(0,0){\includegraphics[width=\unitlength,page=2]{tangles_unoriented_cap.pdf}}%
    \put(0.84512581,0.28861343){\color[rgb]{0,0,0}\makebox(0,0)[lt]{\lineheight{1.25}\smash{\begin{tabular}[t]{l}$x_1$\end{tabular}}}}%
    \put(0.97133508,0.4321475){\color[rgb]{0,0,0}\rotatebox{89.006007}{\makebox(0,0)[lt]{\lineheight{1.25}\smash{\begin{tabular}[t]{l}$\prec$\end{tabular}}}}}%
  \end{picture}%
\endgroup%

\caption{
\\\hspace*{1mm}}
\end{subfigure}
\hfill
\begin{subfigure}[b]{0.20\textwidth}
\hspace*{2mm} 
\begingroup%
  \makeatletter%
  \providecommand\color[2][]{%
    \errmessage{(Inkscape) Color is used for the text in Inkscape, but the package 'color.sty' is not loaded}%
    \renewcommand\color[2][]{}%
  }%
  \providecommand\transparent[1]{%
    \errmessage{(Inkscape) Transparency is used (non-zero) for the text in Inkscape, but the package 'transparent.sty' is not loaded}%
    \renewcommand\transparent[1]{}%
  }%
  \providecommand\rotatebox[2]{#2}%
  \newcommand*\fsize{\dimexpr\f@size pt\relax}%
  \newcommand*\lineheight[1]{\fontsize{\fsize}{#1\fsize}\selectfont}%
  \ifx\svgwidth\undefined%
    \setlength{\unitlength}{70.86614173bp}%
    \ifx\svgscale\undefined%
      \relax%
    \else%
      \setlength{\unitlength}{\unitlength * \real{\svgscale}}%
    \fi%
  \else%
    \setlength{\unitlength}{\svgwidth}%
  \fi%
  \global\let\svgwidth\undefined%
  \global\let\svgscale\undefined%
  \makeatother%
  \begin{picture}(1,1.04)%
    \lineheight{1}%
    \setlength\tabcolsep{0pt}%
    \put(0,0){\includegraphics[width=\unitlength,page=1]{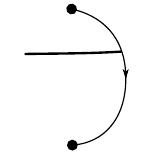}}%
    \put(0.83467506,0.67968369){\color[rgb]{0,0,0}\makebox(0,0)[lt]{\lineheight{1.25}\smash{\begin{tabular}[t]{l}$x_2$\end{tabular}}}}%
    \put(0,0){\includegraphics[width=\unitlength,page=2]{tangles_unoriented_height_exchange.pdf}}%
    \put(0.84512581,0.28861343){\color[rgb]{0,0,0}\makebox(0,0)[lt]{\lineheight{1.25}\smash{\begin{tabular}[t]{l}$x_1$\end{tabular}}}}%
    \put(0,0){\includegraphics[width=\unitlength,page=3]{tangles_unoriented_height_exchange.pdf}}%
    \put(0.04558784,0.29500092){\color[rgb]{0,0,0}\makebox(0,0)[lt]{\lineheight{1.25}\smash{\begin{tabular}[t]{l}$y_1$\end{tabular}}}}%
    \put(0.02647766,0.67545753){\color[rgb]{0,0,0}\makebox(0,0)[lt]{\lineheight{1.25}\smash{\begin{tabular}[t]{l}$y_2$\end{tabular}}}}%
  \end{picture}%
\endgroup%

\caption{$x_1\succ x_2$, $y_2\succ y_1$,\\
\hspace*{2mm} or $x_2\succ x_1$, $y_1\succ y_2$}
\end{subfigure}
\hfill
\begin{subfigure}[b]{0.17\textwidth}
\hspace{4mm} 
\begingroup%
  \makeatletter%
  \providecommand\color[2][]{%
    \errmessage{(Inkscape) Color is used for the text in Inkscape, but the package 'color.sty' is not loaded}%
    \renewcommand\color[2][]{}%
  }%
  \providecommand\transparent[1]{%
    \errmessage{(Inkscape) Transparency is used (non-zero) for the text in Inkscape, but the package 'transparent.sty' is not loaded}%
    \renewcommand\transparent[1]{}%
  }%
  \providecommand\rotatebox[2]{#2}%
  \newcommand*\fsize{\dimexpr\f@size pt\relax}%
  \newcommand*\lineheight[1]{\fontsize{\fsize}{#1\fsize}\selectfont}%
  \ifx\svgwidth\undefined%
    \setlength{\unitlength}{70.86614173bp}%
    \ifx\svgscale\undefined%
      \relax%
    \else%
      \setlength{\unitlength}{\unitlength * \real{\svgscale}}%
    \fi%
  \else%
    \setlength{\unitlength}{\svgwidth}%
  \fi%
  \global\let\svgwidth\undefined%
  \global\let\svgscale\undefined%
  \makeatother%
  \begin{picture}(1,1.04)%
    \lineheight{1}%
    \setlength\tabcolsep{0pt}%
    \put(0,0){\includegraphics[width=\unitlength,page=1]{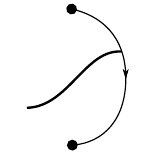}}%
    \put(0.83467506,0.67968369){\color[rgb]{0,0,0}\makebox(0,0)[lt]{\lineheight{1.25}\smash{\begin{tabular}[t]{l}$x_2$\end{tabular}}}}%
    \put(0,0){\includegraphics[width=\unitlength,page=2]{tangles_unoriented_crossing.pdf}}%
    \put(0.84512581,0.28861343){\color[rgb]{0,0,0}\makebox(0,0)[lt]{\lineheight{1.25}\smash{\begin{tabular}[t]{l}$x_1$\end{tabular}}}}%
    \put(0.04558784,0.29500092){\color[rgb]{0,0,0}\makebox(0,0)[lt]{\lineheight{1.25}\smash{\begin{tabular}[t]{l}$y_1$\end{tabular}}}}%
    \put(0.02647766,0.67545753){\color[rgb]{0,0,0}\makebox(0,0)[lt]{\lineheight{1.25}\smash{\begin{tabular}[t]{l}$y_2$\end{tabular}}}}%
    \put(0,0){\includegraphics[width=\unitlength,page=3]{tangles_unoriented_crossing.pdf}}%
    \put(0.94593563,0.47375516){\color[rgb]{0,0,0}\rotatebox{89.156151}{\makebox(0,0)[lt]{\lineheight{1.25}\smash{\begin{tabular}[t]{l}$\succ$\end{tabular}}}}}%
    \put(0.09900697,0.46709855){\color[rgb]{0,0,0}\rotatebox{89.156151}{\makebox(0,0)[lt]{\lineheight{1.25}\smash{\begin{tabular}[t]{l}$\succ$\end{tabular}}}}}%
  \end{picture}%
\endgroup%

\caption{
\\\hspace*{1mm}}
\end{subfigure}      
\hspace{-2mm}
\hfill
\begin{subfigure}[b]{0.17\textwidth}
\hspace{4mm} 
\begingroup%
  \makeatletter%
  \providecommand\color[2][]{%
    \errmessage{(Inkscape) Color is used for the text in Inkscape, but the package 'color.sty' is not loaded}%
    \renewcommand\color[2][]{}%
  }%
  \providecommand\transparent[1]{%
    \errmessage{(Inkscape) Transparency is used (non-zero) for the text in Inkscape, but the package 'transparent.sty' is not loaded}%
    \renewcommand\transparent[1]{}%
  }%
  \providecommand\rotatebox[2]{#2}%
  \newcommand*\fsize{\dimexpr\f@size pt\relax}%
  \newcommand*\lineheight[1]{\fontsize{\fsize}{#1\fsize}\selectfont}%
  \ifx\svgwidth\undefined%
    \setlength{\unitlength}{70.86614173bp}%
    \ifx\svgscale\undefined%
      \relax%
    \else%
      \setlength{\unitlength}{\unitlength * \real{\svgscale}}%
    \fi%
  \else%
    \setlength{\unitlength}{\svgwidth}%
  \fi%
  \global\let\svgwidth\undefined%
  \global\let\svgscale\undefined%
  \makeatother%
  \begin{picture}(1,1.04)%
    \lineheight{1}%
    \setlength\tabcolsep{0pt}%
    \put(0,0){\includegraphics[width=\unitlength,page=1]{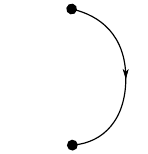}}%
    \put(0.83467506,0.67968369){\color[rgb]{0,0,0}\makebox(0,0)[lt]{\lineheight{1.25}\smash{\begin{tabular}[t]{l}$x_2$\end{tabular}}}}%
    \put(0,0){\includegraphics[width=\unitlength,page=2]{tangles_unoriented_crossing2.pdf}}%
    \put(0.84512581,0.28861343){\color[rgb]{0,0,0}\makebox(0,0)[lt]{\lineheight{1.25}\smash{\begin{tabular}[t]{l}$x_1$\end{tabular}}}}%
    \put(0.04558784,0.29500092){\color[rgb]{0,0,0}\makebox(0,0)[lt]{\lineheight{1.25}\smash{\begin{tabular}[t]{l}$y_1$\end{tabular}}}}%
    \put(0.02647766,0.67545753){\color[rgb]{0,0,0}\makebox(0,0)[lt]{\lineheight{1.25}\smash{\begin{tabular}[t]{l}$y_2$\end{tabular}}}}%
    \put(0,0){\includegraphics[width=\unitlength,page=3]{tangles_unoriented_crossing2.pdf}}%
    \put(0.94594515,0.45571508){\color[rgb]{0,0,0}\rotatebox{89.996297}{\makebox(0,0)[lt]{\lineheight{1.25}\smash{\begin{tabular}[t]{l}$\succ$\end{tabular}}}}}%
    \put(0.10855458,0.45306244){\color[rgb]{0,0,0}\rotatebox{89.996297}{\makebox(0,0)[lt]{\lineheight{1.25}\smash{\begin{tabular}[t]{l}$\succ$\end{tabular}}}}}%
  \end{picture}%
\endgroup%

\caption{
\\\hspace*{1mm}}
\end{subfigure}
\hspace{-2mm}
\hfill
\vspace{-1mm} \\
\caption{Some boundary-ordered tangle diagrams in a biangle (arrows on the boundary arcs indicate the boundary-orientation (Def.\ref{def:surface}))}
\vspace{-2mm}
\label{fig:elementary_tangle_diagrams_in_a_directed_biangle}
\end{figure}

For an (oriented) biangle $B$ viewed as a generalized marked surface as above, choose an internal arc $b$ connecting the two marked points of $B$. Then, cutting $B$ along $b$ yields a unique (up to diffeomorphism) generalized marked surface $(\Sigma',\mathcal{P}')$ as described in Def.\ref{def:cutting}. One easily observes that $(\Sigma',\mathcal{P}')$ is a disjoint union of two biangles. Here is an analog of Prop.\ref{prop:BW_full} for biangles:
\begin{proposition}[quantum trace for biangles; \cite{BW, CL}]
\label{prop:BW_biangles}
Let $A, \omega \in \mathbb{C}^*$ satisfy $A = \omega^{-2}$. Then there is a unique family of algebra homomorphisms
$$
{\rm Tr}^\omega_B \, : \, \mathcal{S}^A_{\rm s}(B) \to \mathbb{C}
$$
defined for all (oriented) biangles $B$, such that

\begin{enumerate}
\item[\rm (1)] (Cutting Property) 
Let $B = (\Sigma,\mathcal{P})$ be a biangle, and $[K,s]\in \mathcal{S}^A_{\rm s}(B)$ be a stated skein for $B$. Let $(\Sigma',\mathcal{P}')$ be the generalized marked surface obtained by cutting $B$ along an internal arc of $B$ connecting the two marked points, as described in Def.\ref{def:cutting}. Then $(\Sigma',\mathcal{P}')$ is disjoint union of two oriented biangles $B_1=(\Sigma_1,\mathcal{P}_1)$ and $B_2=(\Sigma_2,\mathcal{P}_2)$, and the tangle $K \subseteq (\Sigma\setminus\mathcal{P})\times(-1,1)$ yields tangles $K_1 \subseteq (\Sigma_1\setminus\mathcal{P}_1) \times (-1,1)$ and $K_2 \subseteq (\Sigma_2\setminus\mathcal{P}_2) \times (-1,1)$ by this cutting process. Suppose that $K_1$ and $K_2$ are V-tangles. Then one has
$$
{\rm Tr}^\omega_B([K,s]) = \sum_{\mbox{\tiny compatible $s_1,s_2$}} {\rm Tr}^\omega_{B_1}([K_1,s_1]) \, {\rm Tr}^\omega_{B_2}([K_2,s_2]),
$$
where the sum is over all pairs of states $s_1 : \partial K_1 \to \{+,-\}$ and $s_2 : \partial K_2 \to \{+,-\}$ that comprise states $s' : \partial K' \to \{+,-\}$ of $(\Sigma',\mathcal{P}')$ that are compatible with $s : \partial K \to \{+,-\}$ in the sense of Def.\ref{def:compatible_stated_skeins}.

\vs

\item[\rm (2)] (Elementary Cases) For a single biangle $B$, if $[K,s] \in \mathcal{S}^A_{\rm s}(B)$ is a stated skein for $B$ consisting of one component and the  boundary-ordered tangle diagram of $K$ is one of the two cases (I) and (II) in Fig.\ref{fig:elementary_tangle_diagrams_in_a_directed_biangle}:
\begin{enumerate}
\item[\rm (I)] One has
$$
{\rm Tr}^\omega_B([K,s]) = \left\{
\begin{array}{ll}
1 & \mbox{if $s(x)=s(y)$}, \\
0 & \mbox{if $s(x)\neq s(y)$;}
\end{array}
\right.
$$

\item[\rm (II)] One has
$$
{\rm Tr}^\omega_B([K,s]) = \left\{
\begin{array}{ll}
0 & \mbox{if $s(y_1)=s(y_2)$}, \\
-\omega^{-5} & \mbox{if $s(y_1)=-$ and $s(y_2)=+$,} \\
\omega^{-1} & \mbox{if $s(y_1)=+$ and $s(y_2)=-$.}
\end{array}
\right.
$$
\end{enumerate}

\end{enumerate}

\end{proposition}

For any given stated skein for a single biangle $B$, by using skein relations one can resolve it to a linear combination of elementary stated skeins dealt with in Prop.\ref{prop:BW_biangles}(2) and therefore the value under ${\rm Tr}^\omega_B$ can be computed. 
%
%
One caveat is that even when the stated tangle diagram of a stated skein $[K,s]$ over $B$ is without crossing, the value of the biangle quantum trace can still be complicated, instead of being just products of the above elementary cases (I) and (II). For later use, we list some examples:
\begin{lemma}[Biangle quantum trace on remaining elementary tangles]
\label{lem:biangle_quantum_trace_on_remaining_elementary_tangles}
Let $A,\omega$ and $B$ be as in Prop.\ref{prop:BW_biangles}. For a stated skein $[K,s] \in\mathcal{S}^A_{\rm s}(B)$ falling into one of the following cases, we have:
\begin{enumerate}
\item[\rm (1)] If the boundary-ordered tangle diagram of $K$ is as in case (III) of Fig.\ref{fig:elementary_tangle_diagrams_in_a_directed_biangle}, then by \cite[Lem.14]{BW} one has
\begin{align}
\label{eq:BW_biangle_quantum_trace_III}
{\rm Tr}^\omega_B([K,s]) = \left\{
\begin{array}{ll}
0 & \mbox{if $s(x_1)=s(x_2)$}, \\
\omega & \mbox{if $s(x_1)=-$ and $s(x_2)=+$,} \\
-\omega^5 & \mbox{if $s(x_1)=+$ and $s(x_2)=-$.}
\end{array}
\right.
\end{align}

\item[\rm (2)] In case (IV) of Fig.\ref{fig:elementary_tangle_diagrams_in_a_directed_biangle} when the tangle diagram consists of two disjoint parallel lines, if we let
\begin{align}
\label{eq:height_exchange_sign}
\epsilon = \left\{
\begin{array}{ll}
+1 & \mbox{if the vertical orderings are such that $x_1\succ x_2$ and $y_2\succ y_1$,} \\
-1 & \mbox{if $x_2\succ x_1$ and $y_1\succ y_2$,}
\end{array}
\right.
\end{align}
and $s(x_1)=\varepsilon_1$, $s(x_2)=\varepsilon_2$, $s(y_1)=\varepsilon_1'$ and $s(y_2)=\varepsilon_2'$, then one has:
\begin{align}
\label{eq:two_parallel_lines}
{\rm Tr}^\omega_B([K,s]) = \left\{
\begin{array}{ll}
\omega^{2\epsilon} & \mbox{if $\varepsilon_1=\varepsilon_2=\varepsilon_1'=\varepsilon_2'$},  \\
\omega^{2\epsilon} - \omega^{-6\epsilon}  & \mbox{if $\varepsilon_1=\varepsilon_2'=+$ and $\varepsilon_2=\varepsilon_1'=-$}, \\
\omega^{-2\epsilon} & \mbox{if $\varepsilon_1=\varepsilon_1'\neq \varepsilon_2' = \varepsilon_2$}, \\
0 & \mbox{otherwise}.
\end{array}
\right. 
\end{align}

\item[\rm (3)] In case of Fig.\ref{fig:elementary_tangle_diagrams_in_a_directed_biangle}(V), the value of ${\rm Tr}^\omega_B([K,s])$ is given by eq.\eqref{eq:two_parallel_lines} with $\epsilon=-1$ and $s(x_1)=\varepsilon_2$, $s(x_2)=\varepsilon_1$, $s(y_1)=\varepsilon_1'$ and $s(y_2)=\varepsilon_2'$.

\item[\rm (4)] In case of Fig.\ref{fig:elementary_tangle_diagrams_in_a_directed_biangle}(VI), the value of ${\rm Tr}^\omega_B([K,s])$ is given by eq.\eqref{eq:two_parallel_lines} with $\epsilon=+1$ and $s(x_1)=\varepsilon_1$, $s(x_2)=\varepsilon_2$, $s(y_1)=\varepsilon_2'$ and $s(y_2)=\varepsilon_1'$.
\end{enumerate}
\end{lemma}
Genuinely simple basic cases can be conveniently described if one adapts the peculiar picture convention of Bonahon and Wong \cite[\S3.5]{BW} for stated boundary-ordered tangle diagrams in a surface $\frak{S}$, which stipulates that the vertical ordering on each boundary arc $b$ of $\frak{S}$ should match the horizontal ordering with respect to a chosen orientation on $b$. If the stated boundary-ordered tangle diagram in a biangle $B$ consists of disjoint parallel lines under this Bonahon-Wong picture convention where one boundary arc is given the boundary-orientation and the other the opposite of the boundary-orientation, then the value ${\rm Tr}^\omega_B([K,s])$ is either zero or one. Since we do not use the Bonahon-Wong convention, we translate this observation as follows:

\begin{lemma}[biangle factor of a skein whose  diagram consists of (strongly) parallel lines]
\label{lem:parallel_skein_over_biangle}
Let $[K,s]$ be a stated skein over a biangle $B=(\Sigma,\mathcal{P})$. Let $\omega \in \mathbb{C}^*$ and $A = \omega^{-2}$. Let $e$ and $e'$ be the two boundary arcs of $B$, constituting the boundary of $B$. Suppose that the V-tangle $K \subset (\Sigma\setminus\mathcal{P}) \times (-1,1)$ satisfies:
\begin{enumerate}
\item[\rm (1)] the stated (boundary-ordered) tangle diagram of $K$ in $B$ consists of disjoint parallel lines (under our picture convention, not necessarily under Bonahon and Wong's convention);

\item[\rm (2)] the ordering on the segments of $K$ induced by the elevations of the points of $\partial K \cap (e\times (0,1))$ (i.e. the endpoints of the segments of $K$ over $e$) coincides with that induced by the elevations of the points of $\partial K \cap (e' \times (0,1))$.
\end{enumerate}
If, for every segment of $K$, its two endpoints are assigned the same sign by $s$, then ${\rm Tr}^\omega_B([K, s]) = 1$. Otherwise, ${\rm Tr}^\omega_B([K, s])=0$.
\end{lemma}

{\it Proof.} Let $k_1,\ldots,k_r$ be the segments (components) of $K$. For each $j=1,\ldots,r$, let $I_j \subset (-1,1)$ be the image of $k_j$ under the projection $(\Sigma\setminus\mathcal{P})\times (-1,1) \to (-1,1)$ to the second factor, i.e. the collection of elevations of the points of $k_j$. By condition (2), we see that $K$ can be V-isotoped so that $I_1,\ldots,I_r$ are mutually disjoint. So one can rename these segments so that the sequence $k_1,\ldots,k_r$ is arranged in increasing order of elevations. Since ${\rm Tr}^\omega_B : \mathcal{S}^A_{\rm s}(B) \to \mathbb{C}$ is an algebra homomorphism (see Prop.\ref{prop:BW_biangles}), where the product in $\mathcal{S}^A_{\rm s}(B)$ is given by the superposition operation, it follows that
$$
{\rm Tr}^\omega_B([K, s]) = {\rm Tr}^\omega_B([k_1,s_1]) \, {\rm Tr}^\omega_B([k_2,s_2]) \, \cdots \, {\rm Tr}^\omega_B([k_r,s_r]).
$$
Each ${\rm Tr}^\omega_B([k_j,s_j])$ falls under Prop.\ref{prop:BW_biangles}.(2)(I), hence the claim follows.
\qed

\subsection{State-sum formula}
\label{subsec:state-sum_model}

The sought-for state-sum formula for ${\rm Tr}^\omega_\Delta$ requires the notion of ``split" ideal triangulation $\wh{\Delta}$ of an ideal triangulation $\Delta$, where each edge of $\Delta$ is replaced by a biangle. Then the complexity of a skein caused by elevations shall be pushed to biangles by isotopy.

\begin{definition}[\cite{BW}]
\label{def:good_position}
Let $\Delta$ be a triangulation of a generalized marked surface $(\Sigma,\mathcal{P})$, and let $\frak{S} = \Sigma\setminus\mathcal{P}$. Recall that $\Delta$ is a collection of $\mathcal{P}$-arcs in $\Sigma$ satisfying certain conditions. Denote the edges of $\Delta$ by $e_1,e_2,\ldots,e_n$, and the triangles of $\Delta$ by $t_1,t_2,\ldots,t_m$.

\vs

$\bullet$ For each edge $e_i$ of $\Delta$ choose a $\mathcal{P}$-arc $e_i'$ in $\Sigma$ that is $\mathcal{P}$-isotopic to $e_i$, so that no two members of the collection
$$
\wh{\Delta} := \Delta \cup \{e_1',\ldots,e_n'\}
$$
intersect in $\Sigma\setminus\mathcal{P}$. Call this $\wh{\Delta}$ a \ul{\em split $\mathcal{P}$-triangulation} (or, \ul{\em split ideal triangulation}, or \ul{\em split triangulation}) associated to the triangulation $\Delta$. For each $i$, the region bounded by $e_i$ and $e_i'$ is called a \ul{\em biangle} $B_i$. The triangles formed by $\wh{\Delta}$ are in one-to-one correspondence with the triangles of $\Delta$, and we denote them by $\wh{t}_1,\ldots,\wh{t}_m$, correspondingly.

\vs

$\bullet$ A V-tangle $K$ in $\frak{S} \times (-1,1)$ is said to be in a \ul{\em good position} (with respect to $\wh{\Delta}$) if it satisfies (P1) of \S\ref{label:tangles_and_skein_algebra} and all of the following:
\begin{enumerate}
\item[\rm (GP1)] For each constituent edge $e$ of $\wh{\Delta}$, $K$ is transverse to $e \times (-1,1)$. In particular, $K\cap (e\times (-1,1))$ has at most finitely many elements;

\item[\rm (GP2)] For every triangle $\wh{t}_j$ of $\wh{\Delta}$, \, $K \cap (\wh{t}_j \times (-1,1))$ consists of finitely many disjoint arcs, each of which is contained in a constant elevation surface $\wh{t}_j \times *$ and joins two distinct components of $\partial \wh{t}_j \times (-1,1)$, where each component of $\partial \wh{t}_j$ is a side of $\wh{t}_j$ minus the vertices of $\partial \wh{t}_j$;

\item[\rm (GP3)] For every triangle $\wh{t}_j$ of $\wh{\Delta}$, the components of $K \cap (\wh{t}_j \times (-1,1))$ lie at mutually distinct elevations;

\item[\rm (GP4)] None of the crossings of $K$ lie over an edge of $\wh{\Delta}$.
\end{enumerate}

\end{definition}

Note that every triangle of $\wh{\Delta}$ has three distinct sides, even if the corresponding triangle of $\Delta$ is self-folded. In particular, in (GP2) above, for each triangle $\wh{t}_j$ of $\wh{\Delta}$, the number of components of $\partial \wh{t}_j$, hence also the number of components of $\partial \wh{t}_j \times (-1,1)$, is always three. We added (GP4) for convenience.

\begin{lemma}[\cite{BW}]
\label{lem:good_position}
A VH-tangle $K$ in $\frak{S} \times (-1,1)$ can be VH-isotoped to a VH-tangle in a good position.
\end{lemma}

For a V-tangle in a good position, the elevation change occurs only over the biangles. Note that a V-tangle being in a good position does {\em not} guarantee that its tangle diagram on a triangle of $\wh{\Delta}$ has no crossings.

\vs

We find it convenient to define some more terms, both for stating Bonahon and Wong's construction and for later sections.

\begin{definition}
\label{def:segments_and_junctures}
Let $(\Sigma,\mathcal{P})$, $\frak{S}$, $\Delta$, $\wh{\Delta}$, $e_i$, $e_i'$, $t_j$ and $\wh{t}_j$ be as in Def.\ref{def:good_position}; in particular, $i$ runs through $1,\ldots,n$ and $j$ runs through $1,\ldots,m$. Let $K$ be a V-tangle in $\frak{S} \times (-1,1)$ in a good position with respect to $\wh{\Delta}$. Let $E_{\wh{\Delta}} := \bigcup_{i=1}^n (e_i\cup e_i')$. Then $K$ is divided by $E_{\wh{\Delta}}\times (-1,1)$ into \ul{\em (tangle) segments}. Projection on $\Sigma$ of each segment of $K$ is also called a (tangle) segment. The boundary points of a segment are called \ul{\em $\wh{\Delta}$-junctures} of $K$ (or just \ul{\em junctures} of $K$), or \ul{\em endpoints} of that segment.

\vs

A \ul{\em $\wh{\Delta}$-juncture-state} of $K$ is a map $J : \{\mbox{$\wh{\Delta}$-junctures of $K$}\} \to \{+,-\}$ assigning a sign to each juncture.

\end{definition}
The term ``junctures" could have been defined just as elements of $K \cap (E_{\wh{\Delta}} \times (-1,1))$, or their projections in $\Sigma$. The points of $\partial K$ and their projections are also examples of junctures.

\begin{proposition}[state-sum formula of Bonahon-Wong quantum trace; {\cite[\S6]{BW}}]
\label{prop:state-sum_of_BW}
Let $(\Sigma,\mathcal{P})$, $\frak{S}$, $\Delta$, $\wh{\Delta}$, $e_i$, $e_i'$, $B_i$, $t_j$ and $\wh{t}_j$ be as in Def.\ref{def:good_position}, and $K$ be a V-tangle in $\frak{S} \times (-1,1)$ in a good position with respect to $\wh{\Delta}$. Let $s : \partial K \to \{+,-\}$ be a state for the tangle $K$. Let $A, \omega \in \mathbb{C}^*$ satisfy $A = \omega^{-2}$. 

\vs

Let $J : \{\mbox{$\wh{\Delta}$-junctures of $K$}\} \to \{+,-\}$ be a $\wh{\Delta}$-juncture-state of $K$. For each triangle $\wh{t}_j$ of $\wh{\Delta}$, let $K_j := K \cap (\wh{t}_j \times (-1,1))$. Then $[K_j, \, J |_{\partial K_j}]$ is a stated skein for the triangle $\wh{t}_j$, where $J|_{\partial K_j}$ is the restriction of $J$ to $\partial K_j$. Let $k_{j,1},\ldots,k_{j,l_j}$ be the components of $K_j$, i.e. the tangle segments of $K$ over the triangle $\wh{t}_j$, in order of increasing elevation, so that each $[k_{j,\alpha}, J|_{\partial k_{j,\alpha}}]$ is a stated skein for the triangle $\wh{t}_j$, where $J|_{\partial k_{j,\alpha}}$ is the restriction of $J$ to $\partial k_{j,\alpha}$. Then, the element ${\rm Tr}^\omega_{\wh{t}_j}( [k_{j,\alpha}, J|_{\partial k_{j,\alpha}}] )$ of the triangle algebra $\mathcal{T}^\omega_{\wh{t}_j}$ is defined via Prop.\ref{prop:BW_full}(2)(a), and the element ${\rm Tr}^\omega_{\wh{t}_j}([K_j, J|_{\partial K_j}]) $ of $\mathcal{T}^\omega_{\wh{t}_j}$ is given by 
\begin{align}
\label{eq:triangle-factor}
{\rm Tr}^\omega_{\wh{t}_j}([K_j, J|_{\partial K_j}]) = {\rm Tr}^\omega_{\wh{t}_j}( [k_{j,1}, J|_{\partial k_{j,1}}] ) \, {\rm Tr}^\omega_{\wh{t}_j}( [k_{j,2}, J|_{\partial k_{j,2}}] ) \cdots {\rm Tr}^\omega_{\wh{t}_j}( [k_{j,l_j}, J|_{\partial k_{j,l_j}}] ) \in \mathcal{T}^\omega_{\wh{t}_j}.
\end{align}
Via the natural map $\mathcal{T}^\omega_{\wh{t}_j} \to \mathcal{T}^\omega_{t_j}$ induced by the correspondence of the sides of the triangles, this element ${\rm Tr}^\omega_{\wh{t}_j}([K_j, J|_{\partial K_j}])$ can be viewed as an element of $\mathcal{T}^\omega_{t_j}$.

\vs

For each biangle $B_i=(\Sigma_i,\mathcal{P}_i)$ which is bounded by $e_i$ and $e_i'$, let $L_i := K \cap ((\Sigma_i\setminus\mathcal{P}_i)\times (-1,1))$. Then $[L_i, J|_{\partial L_i}]$ is a stated skein for the biangle $B_i$, where $B_i$ is viewed as a generalized marked surface of its own, and $J|_{\partial L_i}$ is the restriction of $J$ to $\partial L_i$. Define the \ul{\em Bonahon-Wong term} for the $\wh{\Delta}$-juncture-state $J$ of the V-tangle $K$ as
\begin{align}
\label{eq:BW_term}
{\rm BW}^\omega_{\wh{\Delta}}(K;J) := \left( \prod_{i=1}^n {\rm Tr}^\omega_{B_i}([L_i, J|_{\partial L_i}]) \right) \left( \bigotimes_{j=1}^m {\rm Tr}^\omega_{\wh{t}_j}([K_j, J|_{\partial K_j}]) \right) ~ \in ~ \bigotimes_{j=1}^m \mathcal{T}^\omega_{t_j},
\end{align}
where the numbers ${\rm Tr}^\omega_{B_i}([L_i, J|_{\partial L_i}]) \in \mathbb{C}$ are given by Prop.\ref{prop:BW_biangles}. Then ${\rm BW}^\omega_{\wh{\Delta}}(K;J) \in \mathcal{Z}^\omega_\Delta \subset \mathcal{T}^\omega_\Delta$.

\vs

Finally, one has
$$
{\rm Tr}^\omega_\Delta ([K,s]) = \sum_{J : J|_{\partial K} = s} {\rm BW}^\omega_{\wh{\Delta}}(K;J) ~ \in ~ \mathcal{Z}^\omega_\Delta ~ \subset ~ \mathcal{T}^\omega_\Delta ~ \subset ~ \bigotimes_{j=1}^m \mathcal{T}^\omega_{t_j},
$$
where the sum is over all $\wh{\Delta}$-juncture-states $J$ that restrict to the given state $s$ at $\partial K$.

\end{proposition}

The above proposition is what one can practically use for actual computation of the values of the Bonahon-Wong quantum trace.

\section{Gabella's quantum holonomy}
\label{sec:Gabella}

\def\cV{\mathcal V}
\def\cP{\mathcal P}
\def\tP{\tilde{\cP}}
\def\tS{\tilde \Sigma_\Delta}
\def\BZ{{\mathbb Z}}
\def\al{\alpha}
\def\sh{\mathrm{s}}
\def\tk{\tilde k}
\def\hD{{\hat \Delta}}
\def\fS{\mathfrak{S}}

\subsection{Branched double cover surface} 
\label{subsec:branched_double_cover} Quantization of the trace-of-monodromy functions (for closed curves) on the Teichm\"uller space, namely the quantum ordering problem mentioned in \S\ref{subsec:deformation_quantization}, is also of interest to physicists, as the coefficients of the monomials of the quantum version of the trace-of-monodromy correspond to the so-called ``framed protected spin characters" in physics. A very interesting solution to this quantum ordering problem is given by Gabella \cite{G}. 
Gabella's construction of quantum holonomy is based on the works of Gaiotto, Moore and Neitzke \cite{GMN13, GMN14} and Galakhov, Longhi and Moore \cite{GLM}. Part of the main ideas of these works are the processes called ``nonabelianization" and ``abelianization", which relate a ${\rm GL}_N$-bundle of a surface $\Sigma$ to an abelian bundle of a certain $N$-fold branched (ramified) covering of $\Sigma$. In our case, $N=2$. We assume $\Sigma$ is connected. We start with the description of a branched double cover of the surface $\Sigma$, and some basic constructions about it.

\vs

Recall that a {\em branched double cover} $\pi:\til{\Sigma}\to \Sigma$, with a finite branching set $\cV \subset \Sigma$, is a continuous map such that the restriction  $\pi':\til{\Sigma} \setminus \pi^{-1}(\cV)\to\Sigma\setminus \cV$ of $\pi$ is an ordinary double covering and the restriction of $\pi$ onto $\pi^{-1}(\cV)$ is bijective. Given the branching set $\cV$, branched double covers $\pi$ of $\Sigma$ are in one-to-one correspondence with ordinary double covers $\pi'$ of $\Sigma\setminus\mathcal{V}$, which are classified by subgroups of index two of the fundamental group of $\Sigma\setminus \cV$, or cohomology classes in $H^1(\Sigma\setminus \cV;\BZ/2\mathbb{Z})$. Here we will identify $\mathbb{Z}/2\mathbb{Z}$ with the multiplicative group $\{+,- \}= \{+1,-1\}$.

\begin{figure}
\centering 
\scalebox{0.70}{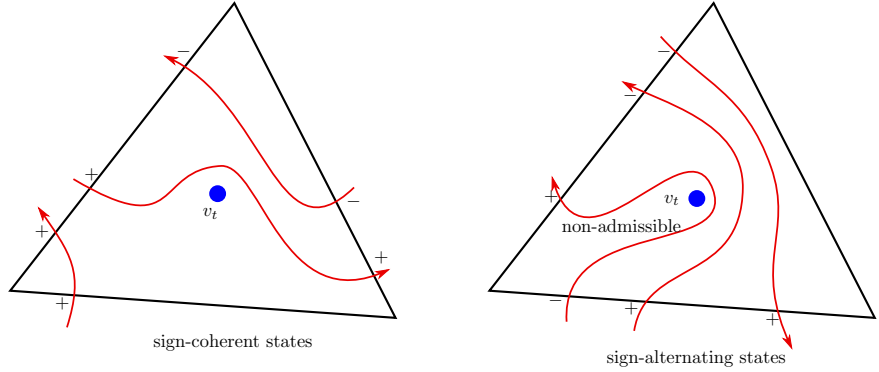}
\vspace{-4mm}
\caption{Gabella lifts in a triangle, and their states}
\vspace{-4mm}
\label{fig:turns}
\end{figure}

\begin{definition}[branched double cover]
\label{def:branched_double_cover}
Let $(\Sigma,\mathcal{P})$ be a triangulable generalized marked surface and $\Delta$ a triangulation of $(\Sigma,\mathcal{P})$. For each triangle $t$ of $\Delta$, choose a point $v_t$ in the interior of $t$. These points are called the \ul{\em branch points}. Denote by $\mathcal{V}$ the set of all branch points. The manifold $\Sigma\setminus \cV$ deformation retracts to the CW-complex consisting of vertices $\cP$ and edges of the triangulation $\Delta$. The 1-cocycle assigning to every edge of this complex the element $-$ of the group $\{+,- \}$ defines a 1-cohomology class in  $H^1(\Sigma\setminus \cV;\BZ/2\mathbb{Z})$, and hence defines a branched double cover
$$
\pi : \til{\Sigma}_\Delta \to \Sigma
$$
with the branching set $\cV$.
\end{definition}
\begin{remark}
The above is a special case of a more general theory of branched $N$-fold cover of $\Sigma$, which corresponds to a ``Seiberg-Witten curve" in physics  \cite{G09, GMN11}. See these references, as well as \cite{G, GMN13, GMN14}, for descriptions using branch cuts.
\end{remark}
A basic observation is that every continuous path $\alpha : [0,1] \to \Sigma\setminus \mathcal{V}$ has two {\em continuous lifts} $[0,1] \to \til{\Sigma}_\Delta\setminus\pi^{-1}(\mathcal{V})$ in the branched double cover. Lifting paths in $\Sigma$ to paths in $\til{\Sigma}_\Delta$ is a crucial part of the nonabelianization process. To keep track of the lifted paths, we need:
\begin{deflem}[states of lifts of marked points]
Let $(\Sigma,\mathcal{P})$, $\Delta$ and $\pi:\til{\Sigma}_\Delta\to\Sigma$ be as in Def.\ref{def:branched_double_cover}. Let $\til{\mathcal{P}} = \pi^{-1}(\mathcal{P})$ be the lifts of the marked points. Define a state function ${\rm s} = {\rm s}_\Delta:\til{\mathcal{P}} \to \{+,-\}$ as follows. First choose a point in $\til{\mathcal{P}}$ and assign value $+$ for the function ${\rm s}$ at this point. There is a unique extension ${\rm s} : \til{\mathcal{P}} \to \{+,-\}$ such that if $p,p'\in \til{\mathcal{P}}$ are connected by a continuous lift of an edge in $\Delta$, then they have different ${\rm s}$ value.  We will say that $p \in \til{\mathcal{P}}$ is on the \ul{\em sheet} ${\rm s}(p) \in \{+,-\}$.
\end{deflem}
The state ${\rm s}: \til{\mathcal{P}} \to \{+,-\}$ depends on the choice of $\Delta$, and there are two possible choices of states ${\rm s}$, for a given $\Delta$; just choose one ${\rm s}$. We now consider lifting paths in $\Sigma$ to paths in $\til{\Sigma}_\Delta$. A crucial example of the branched double cover is the case when the base surface $(\Sigma,\mathcal{P})$ is a triangle.
\begin{definition}
\label{def:Gabella_lift_of_ordered_pair}
Let $t$ be a non-self-folded triangle, viewed as a generalized marked surface $(\Sigma,\mathcal{P})$; that is, $\Sigma$ is diffeomorphic to a closed disc and $\mathcal{P}$ consists of three points on $\partial \Sigma$. Let $\Delta$ be the unique ideal triangulation of $t$, and denote $\til{\Sigma}_\Delta$ by $\til{t}$. Let $x_0$ and $x_1$ be points of $\partial t = \partial \Sigma \setminus \mathcal{P}$ lying in distinct sides of $t$, say $b_0$ and $b_1$ respectively. A \ul{\em Gabella lift of the ordered pair $(x_0,x_1)$} is any proper embedding $\til{\gamma} : [0,1] \to \til{t} \setminus \pi^{-1}(\{v_t\})$ such that $\pi(\til{\gamma}(i))=x_i$ for each $i=0,1$, or its image together with the orientation.
\end{definition}
For a given pair $(x_0,x_1)$, there are two distinct paths in $t\setminus\{v_t\}$ from $x_0$ to $x_1$, up to isotopy in $t\setminus \{v_t\}$. Hence, up to isotopy in $\til{t}$, there are four distinct Gabella lifts of the pair $(x_0,x_1)$, which can be conveniently parametrized by the signs at points $x_0$ and $x_1$ as follows.
\begin{definition}
\label{def:state_of_Gabella_lift_of_ordered_pair}
Let $(\Sigma,\mathcal{P})=t$, $\til{t}$ be as in Def.\ref{def:Gabella_lift_of_ordered_pair}, and $\til{\gamma} : [0,1] \to \til{t}\setminus \pi^{-1}(\{v_t\})$ be a Gabella lift of the pair $(x_0,x_1)$ of points of $\partial t$ lying in distinct sides $b_0$ and $b_1$ of $t$. For each $i=0,1$, among the two continuous lifts in $\til{t}$ of $b_i$, let $\til{b}_i$ be the one that contains $\til{\gamma}(i)$. Let $\til{\gamma}\,'(i)$ be the velocity vector of $\til{\gamma}$ at $\til{\gamma}(i)$. At the point $\til{\gamma}(i) \in \til{b}_i$ choose a tangent vector to $\til{b}_i$ that forms with $\til{\gamma}\,'(i)$ the positive orientation of the surface, and travel along $\til{b}_i$ in the direction of this tangent vector until one reaches a point $\til{p}$ of $\til{\mathcal{P}}$. Define the state of $\til{\gamma}$ at $x_i$ to be $s(x_i):={\rm s}(\til{p})$. The resulting map $s = s_{\til{\gamma}}:\{x_0,x_1\} \to \{+,-\}$ is the \ul{\em state of the Gabella lift $\til{\gamma}$ of the pair $(x_0,x_1)$}.
\end{definition}
See Fig.\ref{fig:turns} for examples. It is easy to observe the following, from the definition.
\begin{lemma}
Given any map $s:\{x_0,x_1\} \to \{+,-\}$, there is a unique, up to isotopy in $\til{t}$, Gabella lift of the pair $(x_0,x_1)$ whose state coincides with $s$. \qed
\end{lemma}
A Gabella lift $\til{\gamma}$ of $(x_0,x_1)$ can be thought of as a lift of a path $\gamma = \pi\circ\til{\gamma}$ in $t$ from $x_0$ to $x_1$. Meanwhile, we should eventually be dealing with the thickenings of $t$ and $\til{t}$. For our purposes, it suffices to consider the following:

\begin{definition}
\label{def:horizontal_triangle_segment}
Let $(\Sigma,\mathcal{P})=t$ and $\til{t}$ be as in Def.\ref{def:Gabella_lift_of_ordered_pair}.

\vs

$\bullet$ A \ul{\em horizontal triangle segment} over $t$ is an oriented VH-tangle $k : [0,1] \to t\times (-1,1)$ in the thickening of $t$ living at a constant elevation $c\in (-1,1)$ (the image of $k$ lies in $t\times \{c\}$), equipped with the upward vertical framing, such that $\partial k \neq {\O}$ and the elements of $\partial k$ lie in distinct components of $\partial t \times (-1,1)$, as in (GP2) of Def.\ref{def:good_position}; in particular, $k(u) = (\gamma(u),c)$ for some path $\gamma: [0,1] \to t$.

\vs

$\bullet$ A \ul{\em Gabella lift of a horizontal triangle segment $k$} is any path $\til{k} : [0,1] \to \til{t} \times (-1,1)$ of the form $\til{k}(u) = (\til{\gamma}(u),c)$, where $\til{\gamma}$ is a Gabella lift of the pair $(\gamma(0),\gamma(1))$.

\vs

$\bullet$ The \ul{\em state of the Gabella lift $\til{k}$ of $k$} is $s=s_{\til{k}}:\{k(0),k(1)\}\to\{+,-\}$ induced by $s_{\til{\gamma}}$ of Def.\ref{def:state_of_Gabella_lift_of_ordered_pair}. 

\vs

$\bullet$ Let $\til{k}$ be a Gabella lift of $k$. Let $e_1$, $e_2$ and $e_3$ be the sides of $t$, appearing clockwise in $\partial t$ in this order. Suppose the endpoints of $\gamma$ lie in $e_i$ and $e_{i+1}$ (with $e_4=e_1$). If the value of the state $s_{\til{\gamma}}$ at the endpoint of $\gamma$ in $e_i$ is $-$ and that at the endpoint in $e_{i+1}$ is $+$, then $\til{k}$ is said to be \ul{\em nonadmissible}. Otherwise, $\til{k}$ is \ul{\em admissible}.
\end{definition}
In Fig.\ref{fig:turns}, the only nonadmissible Gabella lift is the leftmost curve in the picture on the right. Given a horizontal triangle segment $k$ over $t$, there are four Gabella lifts of $k$ up to isotopy in $\til{t}\times(-1,1)$, parametrized by the four possible states $s:\{k(0),k(1)\}\to\{+,-\}$ at the endpoints of $k$. Exactly one of these four Gabella lifts is nonadmissible.

Now we deal with the case of any triangulable generalized marked surface $(\Sigma,\mathcal{P})$. 
\begin{definition}
\label{def:Gabella_lift}
Let $(\Sigma,\mathcal{P})$ be a triangulable generalized marked surface, and $\Delta$ an ideal triangulation of $(\Sigma,\mathcal{P})$. Let $\frak{S} = \Sigma\setminus\mathcal{P}$. Let $\wh{\Delta}$ be a split ideal triangulation of $\Delta$; see Def.\ref{def:good_position}. Let $K$ be an oriented VH-tangle in $\frak{S} \times (-1,1)$ in a good position with respect to $\wh{\Delta}$ (see Def.\ref{def:good_position}), so that the tangle segments of $K$ (see Def.\ref{def:segments_and_junctures}) over each triangle of $\wh{\Delta}$ is a horizontal tangle segment; call these \ul{\em triangle segments of $K$}.

\vs

A \ul{\em Gabella lift} $\til{K}$ (in $\til{\Sigma}_\Delta\times(-1,1)$) of $K$ is the choice of a Gabella lift of each triangle segment of $K$. We say that $\til{K}$ is \ul{\em admissible} if the chosen Gabella lift of each triangle segment of $K$ is admissible. An \ul{\em isotopy} of Gabella lifts of $K$ is a simultaneous isotopy of the Gabella lifts of the triangle segments.
\end{definition}

We now briefly describe physicists' ideas of abelianization and nonabelianization. The original object of our study is the enhanced Teichm\"uller space $\mathcal{X}^+(\Sigma,\mathcal{P})$ of a triangulable generalized marked surface $(\Sigma,\mathcal{P})$, defined in Def.\ref{def:enhanced_Teichmuller_space}. One standard way of studying this space is to view it as the moduli space of ${\rm PSL}(2,\mathbb{R})$-local systems on $\frak{S} = \Sigma\setminus\mathcal{P}$ satisfying some conditions \cite{FG06}, together with certain data at the ``asymptotic boundary points" $\mathcal{P}$. Recall that, for a Lie group ${\rm G}$, a ${\rm G}$-local system on a manifold means a principal ${\rm G}$-bundle on the manifold together with a flat ${\rm G}$-connection on it. The abelianization and nonabelianization processes of Gaiotto, Moore and Neitzke \cite{GMN13} build a correspondence between the moduli space of (certain) ${\rm GL}_N$-local systems on the original surface $\frak{S}$ and the moduli space of (certain) ${\rm GL}_1$-local systems on a branched $N$-fold cover $\til{\frak{S}}$ of $\frak{S}$. The ${\rm GL}_N$ parallel transport map  of a path $\gamma$ in $\frak{S}$ corresponds to the sum of (formal) ${\rm GL}_1$ parallel transports of lifts $\til{\gamma}$ of this path in the branched cover $\til{\frak{S}}$. Only lifts satisfying an admissibility condition contribute to this sum, which corresponds to our admissibility condition when $N=2$. The term associated to $\til{\gamma}$ is written as a product of ${\rm GL}_1$ parallel transports of some elementary paths in $\til{\frak{S}}$, hence it can be viewed as a monomial in some variables. Upon quantization, these monomials should be enhanced to noncommutative monomials, and an appropriate quantum ordering must be chosen for the resulting (Laurent) polynomial to satisfy favorable properties. 

\vs

One answer was suggested by Galakhov, Longhi and Moore \cite{GLM}; namely, consider the Weyl-ordered noncommutative monomial times a power of $q$ to the net sum of signs of self-intersections of the lifted path $\til{\gamma}$. In view of the way they determined the signs, this net sum of signs is really the writhe of a framed path in the thickened $3$-manifold $\til{\frak{S}}\times(-1,1)$ given the ``always going up" elevations that projects to $\til{\gamma}$ in $\til{\frak{S}}$. In particular, a path $\gamma$ is first lifted to a framed path in the 3-dimensional manifold $\frak{S}\times(-1,1)$ with always-going-up elevations, and then the quantum ${\rm GL}_N$ parallel transport is constructed using writhes of lifted paths in $\til{\frak{S}} \times (-1,1)$. So, this method was insufficient to deal with the quantum parallel transport along a {\em closed} path $\gamma$, because it can't be lifted to a closed framed path in $\frak{S} \times (-1,1)$ with always-going-up elevations. Once we choose a starting point of $\gamma$, we lift it to an always-going-up path in $\frak{S}\times(-1,1)$, and then need to add a ``going down" path at the end in order to get a closed path. Gabella's contribution \cite{G} is the consideration of a certain complex-number correction factor associated to this closing-up part, in order to make the final result independent of the choice of starting basepoint. This correction factor, which Gabella refers to as an {\em R-matrix}, is heavily inspired by the biangle factor of Bonahon and Wong \cite{BW} which was reviewed in \S\ref{subsec:biangles} of the present paper.

\vs

In fact, we enhance Gabella's construction and describe it for any {\em general} oriented VH-tangle $K$ in the thickened surface $\frak{S} \times(-1,1)$, not just for a (simple) loop in $\frak{S}$ nor a certain special closed tangle in $\frak{S}\times(-1,1)$ that has always-going-up elevations except at a small part. For now, assume that $K$ is in a good position with respect to some split ideal triangulation $\wh{\Delta}$; later, our main result will allow us to drop this assumption.  Gabella's quantum holonomy along such $K$ in a good position will be given as a sum of certain terms, over all admissible Gabella lifts $\til{K}$ of $K$ up to isotopy. One might want to understand each $\til{K}$ as a framed tangle in $\til{\Sigma}_\Delta \times(-1,1)$ that lifts $K$, but in general, $\til{K}$ cannot be a continuous path in $\til{\Sigma}_\Delta\times(-1,1)$ and must have discontinuities over biangles of $\wh{\Delta}$. However, we still refer to $\til{K}$ as a Gabella lift of $K$ {\em in $\til{\Sigma}_\Delta\times(-1,1)$}, for convenience. 

\vs

In pictures, a Gabella lift $\til{K}$ is drawn in $\frak{S}$ as follows. In each biangle $B$ of $\wh{\Delta}$, just draw the boundary-ordered tangle diagram for $K\cap (B\times (-1,1))$. In each triangle $\wh{t}$ of $\wh{\Delta}$, isotope $\til{K}$ so that the image of $\til{K}$ under the projection $\pi \times {\rm id} : \til{\Sigma}_\Delta\times(-1,1) \to \frak{S}\times(-1,1)$ satisfies (P1) and (P2) of \S\ref{label:tangles_and_skein_algebra} over $\wh{t}$, hence yielding a boundary-ordered tangle diagram in $\wh{t}$. We also indicate the values of the states at the endpoints of each triangle segment over $\wh{t}$. These data forms a \ul{\em diagram of $\til{K}$}. In particular, it is a boundary-ordered tangle diagram in $\frak{S}$ together with some extra data. By a \ul{\em crossing of $\til{K}$} we mean a self-intersection of the image of $\til{K}$ under the projection $\til{\Sigma}_\Delta\times(-1,1) \to \til{\Sigma}_\Delta$, and we indicate this in the diagram of $\til{K}$ drawn in $\frak{S}$ as before, using broken lines. A self-intersection of the diagram of $\til{K}$ that is not a crossing of $\til{K}$ is drawn without broken lines; see Fig.\ref{fig:computation_of_writhe}. 

\subsection{Biangle quantum holonomy via operator invariant}
\label{subsec:biangle_quantum_holonomy}

The strategy of Gabella \cite{G} to construct a quantum holonomy is similar to Bonahon and Wong in the following sense: first define the quantum holonomy for a biangle, and then define the quantum holonomy for a triangulated surface via a state-sum formula. We first focus on biangles. We note that, over a biangle, Gabella \cite[\S5.3--\S5.4]{G} describes the values of the (biangle) quantum holonomy only for simple examples of oriented VH-tangles which he refers to as  ``R-matrix" and ``cup/cap". Gabella mentions that his R-matrix associated to an oriented tangle in a biangle consisting of a single crossing can be interpreted as the R-matrix appearing in the representation theory of the quantum group $\mathcal{U}_q(\frak{gl}_N)$, where we can put $N=2$ for now. This strongly hints that the biangle quantum holonomy of Gabella is related to the Reshetikhin-Turaev invariant of tangles \cite{RT, Turaev16}, associated to a representation of a quantum group. We will eventually confirm this expectation.

\vs

We first review and settle some basic background for constructing operator invariants of tangles in a biangle, by mostly following the contents of Ohtsuki's book \cite{O}. Translating into our language requires us to consider a biangle with fixed labeling of boundary arcs. First, recall from Def.\ref{def:triangulability} and \S\ref{subsec:biangles} that a {\em biangle} $B$ is an example of a generalized marked surface $(\Sigma,\mathcal{P})$. In particular, it is an oriented surface diffeomorphic to a closed disc with two marked points on the boundary, and its boundary $\partial B = \partial \Sigma \setminus \mathcal{P}$ is a disjoint union of two boundary arcs. We will need to work with a following version of biangle.

\begin{figure}[htbp!]
\hfill
\begin{subfigure}[b]{0.17\textwidth}
\hspace*{-3mm} 
\begingroup%
  \makeatletter%
  \providecommand\color[2][]{%
    \errmessage{(Inkscape) Color is used for the text in Inkscape, but the package 'color.sty' is not loaded}%
    \renewcommand\color[2][]{}%
  }%
  \providecommand\transparent[1]{%
    \errmessage{(Inkscape) Transparency is used (non-zero) for the text in Inkscape, but the package 'transparent.sty' is not loaded}%
    \renewcommand\transparent[1]{}%
  }%
  \providecommand\rotatebox[2]{#2}%
  \newcommand*\fsize{\dimexpr\f@size pt\relax}%
  \newcommand*\lineheight[1]{\fontsize{\fsize}{#1\fsize}\selectfont}%
  \ifx\svgwidth\undefined%
    \setlength{\unitlength}{70.86614173bp}%
    \ifx\svgscale\undefined%
      \relax%
    \else%
      \setlength{\unitlength}{\unitlength * \real{\svgscale}}%
    \fi%
  \else%
    \setlength{\unitlength}{\svgwidth}%
  \fi%
  \global\let\svgwidth\undefined%
  \global\let\svgscale\undefined%
  \makeatother%
  \begin{picture}(1,1.04)%
    \lineheight{1}%
    \setlength\tabcolsep{0pt}%
    \put(0,0){\includegraphics[width=\unitlength,page=1]{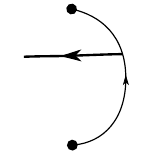}}%
    \put(0.02921547,0.65210342){\color[rgb]{0,0,0}\makebox(0,0)[lt]{\lineheight{1.25}\smash{\begin{tabular}[t]{l}$y$\end{tabular}}}}%
    \put(0,0){\includegraphics[width=\unitlength,page=2]{tangles_identity_forward.pdf}}%
    \put(0.87236481,0.66107025){\color[rgb]{0,0,0}\makebox(0,0)[lt]{\lineheight{1.25}\smash{\begin{tabular}[t]{l}$x$\end{tabular}}}}%
    \put(0.84193916,0.26098465){\color[rgb]{0,0,0}\makebox(0,0)[lt]{\lineheight{1.25}\smash{\begin{tabular}[t]{l}$b_{\rm in}$\end{tabular}}}}%
    \put(0.23333933,0.26718307){\color[rgb]{0,0,0}\makebox(0,0)[lt]{\lineheight{1.25}\smash{\begin{tabular}[t]{l}$b_{\rm out}$\end{tabular}}}}%
  \end{picture}%
\endgroup%

\caption{identity\\\hspace*{4mm}(forward)}
\end{subfigure}
\hfill
\hspace{1mm}
\begin{subfigure}[b]{0.17\textwidth}
\hspace*{-3mm} 
\begingroup%
  \makeatletter%
  \providecommand\color[2][]{%
    \errmessage{(Inkscape) Color is used for the text in Inkscape, but the package 'color.sty' is not loaded}%
    \renewcommand\color[2][]{}%
  }%
  \providecommand\transparent[1]{%
    \errmessage{(Inkscape) Transparency is used (non-zero) for the text in Inkscape, but the package 'transparent.sty' is not loaded}%
    \renewcommand\transparent[1]{}%
  }%
  \providecommand\rotatebox[2]{#2}%
  \newcommand*\fsize{\dimexpr\f@size pt\relax}%
  \newcommand*\lineheight[1]{\fontsize{\fsize}{#1\fsize}\selectfont}%
  \ifx\svgwidth\undefined%
    \setlength{\unitlength}{70.86614173bp}%
    \ifx\svgscale\undefined%
      \relax%
    \else%
      \setlength{\unitlength}{\unitlength * \real{\svgscale}}%
    \fi%
  \else%
    \setlength{\unitlength}{\svgwidth}%
  \fi%
  \global\let\svgwidth\undefined%
  \global\let\svgscale\undefined%
  \makeatother%
  \begin{picture}(1,1.04)%
    \lineheight{1}%
    \setlength\tabcolsep{0pt}%
    \put(0,0){\includegraphics[width=\unitlength,page=1]{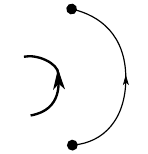}}%
    \put(0.02921547,0.65210342){\color[rgb]{0,0,0}\makebox(0,0)[lt]{\lineheight{1.25}\smash{\begin{tabular}[t]{l}$y_2$\end{tabular}}}}%
    \put(0,0){\includegraphics[width=\unitlength,page=2]{tangles_cup_up.pdf}}%
    \put(0.07276626,0.22511337){\color[rgb]{0,0,0}\makebox(0,0)[lt]{\lineheight{1.25}\smash{\begin{tabular}[t]{l}$y_1$\end{tabular}}}}%
    \put(0.1055255,0.35436626){\color[rgb]{0,0,0}\rotatebox{91.217088}{\makebox(0,0)[lt]{\lineheight{1.25}\smash{\begin{tabular}[t]{l}$\succ$\end{tabular}}}}}%
  \end{picture}%
\endgroup%

\caption{cup (up)\\\hspace*{5mm}($y_1\succ y_2$)}
\end{subfigure}
\hfill
\hspace{-4mm}
\begin{subfigure}[b]{0.17\textwidth}
\hspace*{-3mm} 
\begingroup%
  \makeatletter%
  \providecommand\color[2][]{%
    \errmessage{(Inkscape) Color is used for the text in Inkscape, but the package 'color.sty' is not loaded}%
    \renewcommand\color[2][]{}%
  }%
  \providecommand\transparent[1]{%
    \errmessage{(Inkscape) Transparency is used (non-zero) for the text in Inkscape, but the package 'transparent.sty' is not loaded}%
    \renewcommand\transparent[1]{}%
  }%
  \providecommand\rotatebox[2]{#2}%
  \newcommand*\fsize{\dimexpr\f@size pt\relax}%
  \newcommand*\lineheight[1]{\fontsize{\fsize}{#1\fsize}\selectfont}%
  \ifx\svgwidth\undefined%
    \setlength{\unitlength}{70.86614173bp}%
    \ifx\svgscale\undefined%
      \relax%
    \else%
      \setlength{\unitlength}{\unitlength * \real{\svgscale}}%
    \fi%
  \else%
    \setlength{\unitlength}{\svgwidth}%
  \fi%
  \global\let\svgwidth\undefined%
  \global\let\svgscale\undefined%
  \makeatother%
  \begin{picture}(1,1.04)%
    \lineheight{1}%
    \setlength\tabcolsep{0pt}%
    \put(0,0){\includegraphics[width=\unitlength,page=1]{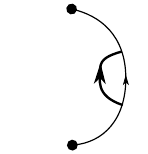}}%
    \put(0.84120651,0.68948082){\color[rgb]{0,0,0}\makebox(0,0)[lt]{\lineheight{1.25}\smash{\begin{tabular}[t]{l}$x_2$\end{tabular}}}}%
    \put(0,0){\includegraphics[width=\unitlength,page=2]{tangles_cap_up.pdf}}%
    \put(0.84512581,0.28861343){\color[rgb]{0,0,0}\makebox(0,0)[lt]{\lineheight{1.25}\smash{\begin{tabular}[t]{l}$x_1$\end{tabular}}}}%
    \put(0.95294832,0.48557454){\color[rgb]{0,0,0}\rotatebox{89.202264}{\makebox(0,0)[lt]{\lineheight{1.25}\smash{\begin{tabular}[t]{l}$\succ$\end{tabular}}}}}%
  \end{picture}%
\endgroup%

\caption{cap (up)\\\hspace*{5mm}($x_1\succ x_2$)}
\end{subfigure}
\hfill
\hspace{-1,5mm}
\begin{subfigure}[b]{0.20\textwidth}
\hspace*{2mm} 
\begingroup%
  \makeatletter%
  \providecommand\color[2][]{%
    \errmessage{(Inkscape) Color is used for the text in Inkscape, but the package 'color.sty' is not loaded}%
    \renewcommand\color[2][]{}%
  }%
  \providecommand\transparent[1]{%
    \errmessage{(Inkscape) Transparency is used (non-zero) for the text in Inkscape, but the package 'transparent.sty' is not loaded}%
    \renewcommand\transparent[1]{}%
  }%
  \providecommand\rotatebox[2]{#2}%
  \newcommand*\fsize{\dimexpr\f@size pt\relax}%
  \newcommand*\lineheight[1]{\fontsize{\fsize}{#1\fsize}\selectfont}%
  \ifx\svgwidth\undefined%
    \setlength{\unitlength}{70.86614173bp}%
    \ifx\svgscale\undefined%
      \relax%
    \else%
      \setlength{\unitlength}{\unitlength * \real{\svgscale}}%
    \fi%
  \else%
    \setlength{\unitlength}{\svgwidth}%
  \fi%
  \global\let\svgwidth\undefined%
  \global\let\svgscale\undefined%
  \makeatother%
  \begin{picture}(1,1.04)%
    \lineheight{1}%
    \setlength\tabcolsep{0pt}%
    \put(0,0){\includegraphics[width=\unitlength,page=1]{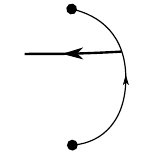}}%
    \put(0.83467506,0.67968369){\color[rgb]{0,0,0}\makebox(0,0)[lt]{\lineheight{1.25}\smash{\begin{tabular}[t]{l}$x_2$\end{tabular}}}}%
    \put(0,0){\includegraphics[width=\unitlength,page=2]{tangles_height_exchange1.pdf}}%
    \put(0.84512581,0.28861343){\color[rgb]{0,0,0}\makebox(0,0)[lt]{\lineheight{1.25}\smash{\begin{tabular}[t]{l}$x_1$\end{tabular}}}}%
    \put(0,0){\includegraphics[width=\unitlength,page=3]{tangles_height_exchange1.pdf}}%
    \put(0.04558784,0.29500092){\color[rgb]{0,0,0}\makebox(0,0)[lt]{\lineheight{1.25}\smash{\begin{tabular}[t]{l}$y_1$\end{tabular}}}}%
    \put(0.02647766,0.67545753){\color[rgb]{0,0,0}\makebox(0,0)[lt]{\lineheight{1.25}\smash{\begin{tabular}[t]{l}$y_2$\end{tabular}}}}%
    \put(0.96829293,0.44942591){\color[rgb]{0,0,0}\rotatebox{89.405289}{\makebox(0,0)[lt]{\lineheight{1.25}\smash{\begin{tabular}[t]{l}$\succ$\end{tabular}}}}}%
    \put(0.04186893,0.53928297){\color[rgb]{0,0,0}\rotatebox{-90.216719}{\makebox(0,0)[lt]{\lineheight{1.25}\smash{\begin{tabular}[t]{l}$\succ$\end{tabular}}}}}%
  \end{picture}%
\endgroup%

\caption{height exchange 1\\\hspace*{3mm} ($x_1\succ x_2$, \, $y_2\succ y_1$)}
\end{subfigure}
\hfill
\begin{subfigure}[b]{0.19\textwidth}
\hspace{4mm} 
\begingroup%
  \makeatletter%
  \providecommand\color[2][]{%
    \errmessage{(Inkscape) Color is used for the text in Inkscape, but the package 'color.sty' is not loaded}%
    \renewcommand\color[2][]{}%
  }%
  \providecommand\transparent[1]{%
    \errmessage{(Inkscape) Transparency is used (non-zero) for the text in Inkscape, but the package 'transparent.sty' is not loaded}%
    \renewcommand\transparent[1]{}%
  }%
  \providecommand\rotatebox[2]{#2}%
  \newcommand*\fsize{\dimexpr\f@size pt\relax}%
  \newcommand*\lineheight[1]{\fontsize{\fsize}{#1\fsize}\selectfont}%
  \ifx\svgwidth\undefined%
    \setlength{\unitlength}{70.86614173bp}%
    \ifx\svgscale\undefined%
      \relax%
    \else%
      \setlength{\unitlength}{\unitlength * \real{\svgscale}}%
    \fi%
  \else%
    \setlength{\unitlength}{\svgwidth}%
  \fi%
  \global\let\svgwidth\undefined%
  \global\let\svgscale\undefined%
  \makeatother%
  \begin{picture}(1,1.04)%
    \lineheight{1}%
    \setlength\tabcolsep{0pt}%
    \put(0,0){\includegraphics[width=\unitlength,page=1]{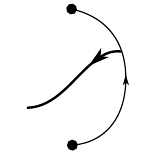}}%
    \put(0.83467506,0.67968369){\color[rgb]{0,0,0}\makebox(0,0)[lt]{\lineheight{1.25}\smash{\begin{tabular}[t]{l}$x_2$\end{tabular}}}}%
    \put(0,0){\includegraphics[width=\unitlength,page=2]{tangles_crossing_positive.pdf}}%
    \put(0.84512581,0.28861343){\color[rgb]{0,0,0}\makebox(0,0)[lt]{\lineheight{1.25}\smash{\begin{tabular}[t]{l}$x_1$\end{tabular}}}}%
    \put(0.04558784,0.29500092){\color[rgb]{0,0,0}\makebox(0,0)[lt]{\lineheight{1.25}\smash{\begin{tabular}[t]{l}$y_1$\end{tabular}}}}%
    \put(0.02647766,0.67545753){\color[rgb]{0,0,0}\makebox(0,0)[lt]{\lineheight{1.25}\smash{\begin{tabular}[t]{l}$y_2$\end{tabular}}}}%
    \put(0,0){\includegraphics[width=\unitlength,page=3]{tangles_crossing_positive.pdf}}%
    \put(0.96355645,0.45902338){\color[rgb]{0,0,0}\rotatebox{90.332903}{\makebox(0,0)[lt]{\lineheight{1.25}\smash{\begin{tabular}[t]{l}$\succ$\end{tabular}}}}}%
    \put(0.09554895,0.46717601){\color[rgb]{0,0,0}\rotatebox{90.332903}{\makebox(0,0)[lt]{\lineheight{1.25}\smash{\begin{tabular}[t]{l}$\succ$\end{tabular}}}}}%
  \end{picture}%
\endgroup%

\caption{positive crossing\\\hspace*{2mm} ($x_1\succ x_2$, \, $y_1\succ y_2$)}
\end{subfigure}      
\hspace{-2mm}
\hfill
\\
\vspace{2mm}
\hfill
\begin{subfigure}[b]{0.17\textwidth}
\hspace*{-3mm} 
\begingroup%
  \makeatletter%
  \providecommand\color[2][]{%
    \errmessage{(Inkscape) Color is used for the text in Inkscape, but the package 'color.sty' is not loaded}%
    \renewcommand\color[2][]{}%
  }%
  \providecommand\transparent[1]{%
    \errmessage{(Inkscape) Transparency is used (non-zero) for the text in Inkscape, but the package 'transparent.sty' is not loaded}%
    \renewcommand\transparent[1]{}%
  }%
  \providecommand\rotatebox[2]{#2}%
  \newcommand*\fsize{\dimexpr\f@size pt\relax}%
  \newcommand*\lineheight[1]{\fontsize{\fsize}{#1\fsize}\selectfont}%
  \ifx\svgwidth\undefined%
    \setlength{\unitlength}{70.86614173bp}%
    \ifx\svgscale\undefined%
      \relax%
    \else%
      \setlength{\unitlength}{\unitlength * \real{\svgscale}}%
    \fi%
  \else%
    \setlength{\unitlength}{\svgwidth}%
  \fi%
  \global\let\svgwidth\undefined%
  \global\let\svgscale\undefined%
  \makeatother%
  \begin{picture}(1,1.04)%
    \lineheight{1}%
    \setlength\tabcolsep{0pt}%
    \put(0,0){\includegraphics[width=\unitlength,page=1]{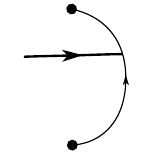}}%
    \put(0.02921547,0.65210342){\color[rgb]{0,0,0}\makebox(0,0)[lt]{\lineheight{1.25}\smash{\begin{tabular}[t]{l}$y$\end{tabular}}}}%
    \put(0,0){\includegraphics[width=\unitlength,page=2]{tangles_identity_backward.pdf}}%
    \put(0.87236481,0.66107025){\color[rgb]{0,0,0}\makebox(0,0)[lt]{\lineheight{1.25}\smash{\begin{tabular}[t]{l}$x$\end{tabular}}}}%
    \put(0.84193916,0.26098465){\color[rgb]{0,0,0}\makebox(0,0)[lt]{\lineheight{1.25}\smash{\begin{tabular}[t]{l}$b_{\rm in}$\end{tabular}}}}%
    \put(0.23333933,0.26718307){\color[rgb]{0,0,0}\makebox(0,0)[lt]{\lineheight{1.25}\smash{\begin{tabular}[t]{l}$b_{\rm out}$\end{tabular}}}}%
  \end{picture}%
\endgroup%

\caption{identity\\\hspace*{5mm}(backward)}
\end{subfigure}
\hfill
\hspace{1mm}
\begin{subfigure}[b]{0.17\textwidth}
\hspace*{-1mm} 
\begingroup%
  \makeatletter%
  \providecommand\color[2][]{%
    \errmessage{(Inkscape) Color is used for the text in Inkscape, but the package 'color.sty' is not loaded}%
    \renewcommand\color[2][]{}%
  }%
  \providecommand\transparent[1]{%
    \errmessage{(Inkscape) Transparency is used (non-zero) for the text in Inkscape, but the package 'transparent.sty' is not loaded}%
    \renewcommand\transparent[1]{}%
  }%
  \providecommand\rotatebox[2]{#2}%
  \newcommand*\fsize{\dimexpr\f@size pt\relax}%
  \newcommand*\lineheight[1]{\fontsize{\fsize}{#1\fsize}\selectfont}%
  \ifx\svgwidth\undefined%
    \setlength{\unitlength}{70.86614173bp}%
    \ifx\svgscale\undefined%
      \relax%
    \else%
      \setlength{\unitlength}{\unitlength * \real{\svgscale}}%
    \fi%
  \else%
    \setlength{\unitlength}{\svgwidth}%
  \fi%
  \global\let\svgwidth\undefined%
  \global\let\svgscale\undefined%
  \makeatother%
  \begin{picture}(1,1.04)%
    \lineheight{1}%
    \setlength\tabcolsep{0pt}%
    \put(0,0){\includegraphics[width=\unitlength,page=1]{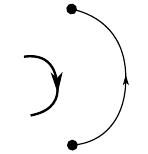}}%
    \put(0.02921547,0.65210342){\color[rgb]{0,0,0}\makebox(0,0)[lt]{\lineheight{1.25}\smash{\begin{tabular}[t]{l}$y_2$\end{tabular}}}}%
    \put(0,0){\includegraphics[width=\unitlength,page=2]{tangles_cup_down.pdf}}%
    \put(0.07276626,0.22511337){\color[rgb]{0,0,0}\makebox(0,0)[lt]{\lineheight{1.25}\smash{\begin{tabular}[t]{l}$y_1$\end{tabular}}}}%
    \put(0.1055255,0.35436626){\color[rgb]{0,0,0}\rotatebox{91.217088}{\makebox(0,0)[lt]{\lineheight{1.25}\smash{\begin{tabular}[t]{l}$\succ$\end{tabular}}}}}%
  \end{picture}%
\endgroup%

\caption{cup (down)\\\hspace*{5mm}($y_1\succ y_2$)}
\end{subfigure}
\hfill
\hspace{-2mm}
\begin{subfigure}[b]{0.17\textwidth}
\hspace*{0mm} 
\begingroup%
  \makeatletter%
  \providecommand\color[2][]{%
    \errmessage{(Inkscape) Color is used for the text in Inkscape, but the package 'color.sty' is not loaded}%
    \renewcommand\color[2][]{}%
  }%
  \providecommand\transparent[1]{%
    \errmessage{(Inkscape) Transparency is used (non-zero) for the text in Inkscape, but the package 'transparent.sty' is not loaded}%
    \renewcommand\transparent[1]{}%
  }%
  \providecommand\rotatebox[2]{#2}%
  \newcommand*\fsize{\dimexpr\f@size pt\relax}%
  \newcommand*\lineheight[1]{\fontsize{\fsize}{#1\fsize}\selectfont}%
  \ifx\svgwidth\undefined%
    \setlength{\unitlength}{70.86614173bp}%
    \ifx\svgscale\undefined%
      \relax%
    \else%
      \setlength{\unitlength}{\unitlength * \real{\svgscale}}%
    \fi%
  \else%
    \setlength{\unitlength}{\svgwidth}%
  \fi%
  \global\let\svgwidth\undefined%
  \global\let\svgscale\undefined%
  \makeatother%
  \begin{picture}(1,1.04)%
    \lineheight{1}%
    \setlength\tabcolsep{0pt}%
    \put(0,0){\includegraphics[width=\unitlength,page=1]{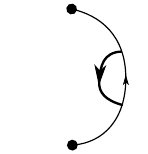}}%
    \put(0.84120651,0.68948082){\color[rgb]{0,0,0}\makebox(0,0)[lt]{\lineheight{1.25}\smash{\begin{tabular}[t]{l}$x_2$\end{tabular}}}}%
    \put(0,0){\includegraphics[width=\unitlength,page=2]{tangles_cap_down.pdf}}%
    \put(0.84512581,0.28861343){\color[rgb]{0,0,0}\makebox(0,0)[lt]{\lineheight{1.25}\smash{\begin{tabular}[t]{l}$x_1$\end{tabular}}}}%
    \put(0.95294832,0.48557454){\color[rgb]{0,0,0}\rotatebox{89.202264}{\makebox(0,0)[lt]{\lineheight{1.25}\smash{\begin{tabular}[t]{l}$\succ$\end{tabular}}}}}%
  \end{picture}%
\endgroup%

\caption{cap (down)\\\hspace*{5mm}($x_1\succ x_2$)}
\end{subfigure}
\hfill
\hspace{1mm}
\begin{subfigure}[b]{0.20\textwidth}
\hspace*{5mm} 
\begingroup%
  \makeatletter%
  \providecommand\color[2][]{%
    \errmessage{(Inkscape) Color is used for the text in Inkscape, but the package 'color.sty' is not loaded}%
    \renewcommand\color[2][]{}%
  }%
  \providecommand\transparent[1]{%
    \errmessage{(Inkscape) Transparency is used (non-zero) for the text in Inkscape, but the package 'transparent.sty' is not loaded}%
    \renewcommand\transparent[1]{}%
  }%
  \providecommand\rotatebox[2]{#2}%
  \newcommand*\fsize{\dimexpr\f@size pt\relax}%
  \newcommand*\lineheight[1]{\fontsize{\fsize}{#1\fsize}\selectfont}%
  \ifx\svgwidth\undefined%
    \setlength{\unitlength}{70.86614173bp}%
    \ifx\svgscale\undefined%
      \relax%
    \else%
      \setlength{\unitlength}{\unitlength * \real{\svgscale}}%
    \fi%
  \else%
    \setlength{\unitlength}{\svgwidth}%
  \fi%
  \global\let\svgwidth\undefined%
  \global\let\svgscale\undefined%
  \makeatother%
  \begin{picture}(1,1.04)%
    \lineheight{1}%
    \setlength\tabcolsep{0pt}%
    \put(0,0){\includegraphics[width=\unitlength,page=1]{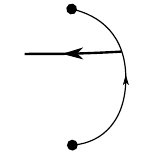}}%
    \put(0.83467506,0.67968369){\color[rgb]{0,0,0}\makebox(0,0)[lt]{\lineheight{1.25}\smash{\begin{tabular}[t]{l}$x_2$\end{tabular}}}}%
    \put(0,0){\includegraphics[width=\unitlength,page=2]{tangles_height_exchange2.pdf}}%
    \put(0.84512581,0.28861343){\color[rgb]{0,0,0}\makebox(0,0)[lt]{\lineheight{1.25}\smash{\begin{tabular}[t]{l}$x_1$\end{tabular}}}}%
    \put(0,0){\includegraphics[width=\unitlength,page=3]{tangles_height_exchange2.pdf}}%
    \put(0.04558784,0.29500092){\color[rgb]{0,0,0}\makebox(0,0)[lt]{\lineheight{1.25}\smash{\begin{tabular}[t]{l}$y_1$\end{tabular}}}}%
    \put(0.02647766,0.67545753){\color[rgb]{0,0,0}\makebox(0,0)[lt]{\lineheight{1.25}\smash{\begin{tabular}[t]{l}$y_2$\end{tabular}}}}%
    \put(0.96829293,0.44942591){\color[rgb]{0,0,0}\rotatebox{89.405289}{\makebox(0,0)[lt]{\lineheight{1.25}\smash{\begin{tabular}[t]{l}$\prec$\end{tabular}}}}}%
    \put(0.04186893,0.53928297){\color[rgb]{0,0,0}\rotatebox{-90.216719}{\makebox(0,0)[lt]{\lineheight{1.25}\smash{\begin{tabular}[t]{l}$\prec$\end{tabular}}}}}%
  \end{picture}%
\endgroup%

\caption{height exchange 2\\\hspace*{3mm} ($x_2\succ x_1$, \, $y_1\succ y_2$)}
\end{subfigure}
\hfill
\hspace{2mm}
\begin{subfigure}[b]{0.19\textwidth}
\hspace{4mm} 
\begingroup%
  \makeatletter%
  \providecommand\color[2][]{%
    \errmessage{(Inkscape) Color is used for the text in Inkscape, but the package 'color.sty' is not loaded}%
    \renewcommand\color[2][]{}%
  }%
  \providecommand\transparent[1]{%
    \errmessage{(Inkscape) Transparency is used (non-zero) for the text in Inkscape, but the package 'transparent.sty' is not loaded}%
    \renewcommand\transparent[1]{}%
  }%
  \providecommand\rotatebox[2]{#2}%
  \newcommand*\fsize{\dimexpr\f@size pt\relax}%
  \newcommand*\lineheight[1]{\fontsize{\fsize}{#1\fsize}\selectfont}%
  \ifx\svgwidth\undefined%
    \setlength{\unitlength}{70.86614173bp}%
    \ifx\svgscale\undefined%
      \relax%
    \else%
      \setlength{\unitlength}{\unitlength * \real{\svgscale}}%
    \fi%
  \else%
    \setlength{\unitlength}{\svgwidth}%
  \fi%
  \global\let\svgwidth\undefined%
  \global\let\svgscale\undefined%
  \makeatother%
  \begin{picture}(1,1.04)%
    \lineheight{1}%
    \setlength\tabcolsep{0pt}%
    \put(0,0){\includegraphics[width=\unitlength,page=1]{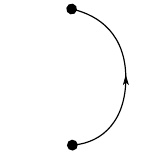}}%
    \put(0.83467506,0.67968369){\color[rgb]{0,0,0}\makebox(0,0)[lt]{\lineheight{1.25}\smash{\begin{tabular}[t]{l}$x_2$\end{tabular}}}}%
    \put(0,0){\includegraphics[width=\unitlength,page=2]{tangles_crossing_negative.pdf}}%
    \put(0.84512581,0.28861343){\color[rgb]{0,0,0}\makebox(0,0)[lt]{\lineheight{1.25}\smash{\begin{tabular}[t]{l}$x_1$\end{tabular}}}}%
    \put(0.04558784,0.29500092){\color[rgb]{0,0,0}\makebox(0,0)[lt]{\lineheight{1.25}\smash{\begin{tabular}[t]{l}$y_1$\end{tabular}}}}%
    \put(0.02647766,0.67545753){\color[rgb]{0,0,0}\makebox(0,0)[lt]{\lineheight{1.25}\smash{\begin{tabular}[t]{l}$y_2$\end{tabular}}}}%
    \put(0,0){\includegraphics[width=\unitlength,page=3]{tangles_crossing_negative.pdf}}%
    \put(0.96355645,0.45902338){\color[rgb]{0,0,0}\rotatebox{90.332903}{\makebox(0,0)[lt]{\lineheight{1.25}\smash{\begin{tabular}[t]{l}$\succ$\end{tabular}}}}}%
    \put(0.09554895,0.46717601){\color[rgb]{0,0,0}\rotatebox{90.332903}{\makebox(0,0)[lt]{\lineheight{1.25}\smash{\begin{tabular}[t]{l}$\succ$\end{tabular}}}}}%
  \end{picture}%
\endgroup%

\caption{negative crossing\\\hspace*{2mm} ($x_1\succ x_2$, \, $y_1\succ y_2$)}
\end{subfigure}
\hspace{-2mm}
\hfill
\vspace{-1mm} \\
\caption{Elementary boundary-ordered oriented tangle diagrams in a directed biangle (orientations on boundary arcs $b_{\rm out}$ and $b_{\rm in}$ are as in Def.\ref{def:directed_biangle})}
\vspace{-2mm}
\label{fig:elementary_bo_oriented_tangle_diagrams_in_a_directed_biangle}
\end{figure}

\begin{definition}[directed biangle]
\label{def:directed_biangle}
$\bullet$ A \ul{\em direction} of a biangle $B$ is a bijective map
\begin{align}
{\rm dir} = {\rm dir} _B : \{\mbox{the two boundary arcs of $B$}\} \to \{{\rm in}, {\rm out}\}.
\end{align}
The boundary arc mapping to ${\rm in}$ is called the \ul{\em inward} boundary arc $b_{\rm in}$, and the one mapping to ${\rm out}$ the \ul{\em outward} boundary arc $b_{\rm out}$.

\vs

$\bullet$ A \ul{\em directed biangle} $\vec{B}$ is a pair $(B,{\rm dir})$ consisting of a biangle $B$ and a direction ${\rm dir}$ of $B$.

\vs

$\bullet$ In a directed biangle $\vec{B}=(B,{\rm dir})$, the outward boundary arc $b_{\rm out}$ is given the boundary-orientation (see Def.\ref{def:surface}) on $b_{\rm out}$ coming from the surface orientation on $B$. The inward boundary arc $b_{\rm in}$ is given the orientation opposite to the boundary-orientation on $b_{\rm in}$. These orientations are depicted in Fig.\ref{fig:elementary_bo_oriented_tangle_diagrams_in_a_directed_biangle}.
\end{definition}
Any two directed biangles can be identified by a diffeomorphism preserving the direction. Note that once a direction is chosen on a biangle, there is one of the two marked points of $B$ that both the oriented boundary arcs $b_{\rm in}$ and $b_{\rm out}$ point toward. So the choice of a direction is equivalent to the choice of such a distinguished marked point of $B$. To get Ohtsuki's pictures, one should rotate our pictures by 90 degrees clockwise.

\vs

\begin{definition}
\label{def:D_B}
$\bullet$ Let $\mathcal{D}(\vec{B})$ be the set of all equivalence classes $[D]$ of boundary-ordered oriented tangle diagrams $D$ (see Def.\ref{def:tangle_diagram}) in $\vec{B}$, where two boundary-ordered oriented tangle diagrams are defined to be \ul{\em equivalent} if one can be obtained from the other by a sequence of isotopies of boundary-ordered oriented tangle diagrams and moves (M1) in \S\ref{label:tangles_and_skein_algebra}, i.e. the framed Reidemeister moves I, II and III.

\vs

$\bullet$ Let $\mathcal{D}_{\rm s}(\vec{B})$ be the set of all equivalence classes $[D,s]$ of stated boundary-ordered oriented tangle diagrams $(D,s)$, where the equivalence is defined analogously.
\end{definition}
By Prop.\ref{prop:tangle_diagrams_and_tangles}(4) we see that $\mathcal{D}(\vec{B})$ is in one-to-one correspondence with the set of all VH-isotopy classes of VH-tangles in the thickening of a directed biangle. As any two directed biangles $\vec{B}$ and $\vec{B}'$ are identified by a diffeomorphism preserving the directions, the sets $\mathcal{D}(\vec{B})$ and $\mathcal{D}(\vec{B}')$ can be naturally identified. We now study some basic operations on $\mathcal{D}(\vec{B})$.

\begin{definition}
Let $\vec{B}$ be a directed biangle, with the underlying biangle denoted by $B = (\Sigma,\mathcal{P})$. Identify $\frak{S} = \Sigma\setminus\mathcal{P}$ with $[0,1] \times \mathbb{R}$ by choosing a diffeomorphism such that $b_{\rm in}$ maps to $\{1\}\times \mathbb{R}$ and $b_{\rm out}$ to $\{0\} \times \mathbb{R}$, where the orientations of $b_{\rm in}$ and $b_{\rm out}$ match the usual increasing orientation of $\mathbb{R}$. The $\mathbb{R}$-coordinate is called a \ul{\em horizontal} coordinate.

\vs

We say that a collection $D_1,\ldots,D_n$ of tangle diagrams in $\vec{B}$ are \ul{\em horizontally disjoint} if the images of them under the second projection $[0,1] \times\mathbb{R} \to \mathbb{R}$ are mutually disjoint. We say $D_i$ is \ul{\em horizontally higher} than $D_j$ if the $\mathbb{R}$-coordinates of points of $D_i$ are bigger than those of $D_j$.

\vs

We say that a collection $D_1,\ldots,D_n$ of boundary-ordered tangle diagrams in $\vec{B}$ are \ul{\em vertically disjoint} if for each pair $i,j$ of distinct indices in $\{1,\ldots,n\}$, either $D_i$ is \ul{\em vertically higher} than $D_j$, meaning on each of the two boundary arcs $b$ any element of $\partial_b D_i$ has higher vertical ordering than any element of $\partial_b D_j$, or $D_j$ is vertically higher than $D_i$.

\vs

A collection of boundary-ordered tangle diagrams in $\vec{B}$ is said to be \ul{\em completely disjoint} if it is horizontally disjoint and vertically disjoint.
\end{definition}

For a tangle diagram $D$ in a directed biangle $\vec{B}$, we write
\begin{align}
\label{eq:partial_in_D_partial_out_D}
\partial_{\rm in}D = \partial_{b_{\rm in}}D \quad\mbox{and}\quad \partial_{\rm out}D = \partial_{b_{\rm out}}D,
\end{align}
where $b_{\rm in}$ and $b_{\rm out}$ are inward and outward boundary arcs of $\vec{B}$ respectively.

\begin{figure}
\centering 
\begin{subfigure}[b]{0.4\textwidth}
\hspace{13mm} \scalebox{1.0}{
\begingroup%
  \makeatletter%
  \providecommand\color[2][]{%
    \errmessage{(Inkscape) Color is used for the text in Inkscape, but the package 'color.sty' is not loaded}%
    \renewcommand\color[2][]{}%
  }%
  \providecommand\transparent[1]{%
    \errmessage{(Inkscape) Transparency is used (non-zero) for the text in Inkscape, but the package 'transparent.sty' is not loaded}%
    \renewcommand\transparent[1]{}%
  }%
  \providecommand\rotatebox[2]{#2}%
  \newcommand*\fsize{\dimexpr\f@size pt\relax}%
  \newcommand*\lineheight[1]{\fontsize{\fsize}{#1\fsize}\selectfont}%
  \ifx\svgwidth\undefined%
    \setlength{\unitlength}{70.86614173bp}%
    \ifx\svgscale\undefined%
      \relax%
    \else%
      \setlength{\unitlength}{\unitlength * \real{\svgscale}}%
    \fi%
  \else%
    \setlength{\unitlength}{\svgwidth}%
  \fi%
  \global\let\svgwidth\undefined%
  \global\let\svgscale\undefined%
  \makeatother%
  \begin{picture}(1,1.04)%
    \lineheight{1}%
    \setlength\tabcolsep{0pt}%
    \put(0,0){\includegraphics[width=\unitlength,page=1]{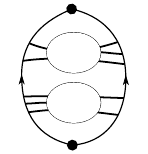}}%
    \put(0.42195251,0.65734695){\color[rgb]{0,0,0}\makebox(0,0)[lt]{\lineheight{1.25}\smash{\begin{tabular}[t]{l}$D_1$\end{tabular}}}}%
    \put(0.41569336,0.31735359){\color[rgb]{0,0,0}\makebox(0,0)[lt]{\lineheight{1.25}\smash{\begin{tabular}[t]{l}$D_2$\end{tabular}}}}%
  \end{picture}%
\endgroup%
}
\vspace{-2mm}
\caption*{tensor product $[D_1]\otimes[D_2]$}
\end{subfigure}
~
\begin{subfigure}[b]{0.4\textwidth}
\hspace{15mm}  \scalebox{1.0}{
\begingroup%
  \makeatletter%
  \providecommand\color[2][]{%
    \errmessage{(Inkscape) Color is used for the text in Inkscape, but the package 'color.sty' is not loaded}%
    \renewcommand\color[2][]{}%
  }%
  \providecommand\transparent[1]{%
    \errmessage{(Inkscape) Transparency is used (non-zero) for the text in Inkscape, but the package 'transparent.sty' is not loaded}%
    \renewcommand\transparent[1]{}%
  }%
  \providecommand\rotatebox[2]{#2}%
  \newcommand*\fsize{\dimexpr\f@size pt\relax}%
  \newcommand*\lineheight[1]{\fontsize{\fsize}{#1\fsize}\selectfont}%
  \ifx\svgwidth\undefined%
    \setlength{\unitlength}{85.03937008bp}%
    \ifx\svgscale\undefined%
      \relax%
    \else%
      \setlength{\unitlength}{\unitlength * \real{\svgscale}}%
    \fi%
  \else%
    \setlength{\unitlength}{\svgwidth}%
  \fi%
  \global\let\svgwidth\undefined%
  \global\let\svgscale\undefined%
  \makeatother%
  \begin{picture}(1,0.86666667)%
    \lineheight{1}%
    \setlength\tabcolsep{0pt}%
    \put(0,0){\includegraphics[width=\unitlength,page=1]{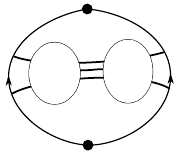}}%
    \put(0.22289921,0.43821195){\color[rgb]{0,0,0}\makebox(0,0)[lt]{\lineheight{1.25}\smash{\begin{tabular}[t]{l}$D_1$\end{tabular}}}}%
    \put(0.65073598,0.44373016){\color[rgb]{0,0,0}\makebox(0,0)[lt]{\lineheight{1.25}\smash{\begin{tabular}[t]{l}$D_2$\end{tabular}}}}%
    \put(0,0){\includegraphics[width=\unitlength,page=2]{tangles_composition.pdf}}%
  \end{picture}%
\endgroup%
}
\vspace{-2mm}
\caption*{composition $[D_1] \circ [D_2]$}
\end{subfigure}
\vspace{-2mm}
\caption{Tensor product and composition of tangle diagrams in a directed biangle}
\vspace{-3mm}
\label{fig:tensor_product_and_composition_of_tangle_diagrams_in_a_directed_biangle}
\end{figure}

\begin{definition}[tensor product and composition of oriented tangle diagrams in a directed biangle]
\label{def:tensor_product_and_composition_of_oriented_tangle_diagrams_in_a_directed_biangle}

Let $\vec{B}$ be a directed biangle, and let $[D_1],[D_2] \in \mathcal{D}(\vec{B})$.

$\bullet$ Let $D_1'$ and $D_2'$ be boundary-ordered oriented tangle diagrams in $\vec{B}$ that are isotopic to $D_1$ and $D_2$ respectively, such that $D_1'$ and $D_2'$ are completely disjoint and $D_1'$ is horizontally higher than $D_2'$. Let $D$ be the union of $D_1'$ and $D_2'$, with $D_1'$ being set to be vertically higher than $D_2'$; see Fig.\ref{fig:tensor_product_and_composition_of_tangle_diagrams_in_a_directed_biangle}. Define the \ul{\em tensor product} of the equivalence classes of boundary-ordered oriented tangle diagrams $[D_1]$ and $[D_2]$ in $\vec{B}$ as $[D_1] \otimes [D_2] := [D]$.

\vs

$\bullet$ Suppose that $\partial_{\rm in} D_1$ and $\partial_{\rm out} D_2$ have the same cardinality, that the unique bijection $\partial_{\rm in} D_1 \to \partial_{\rm out} D_2$ that preserves the horizontal orderings also preserves the vertical orderings, and that this bijection is compatible with the orientations on $D_1$ and $D_2$ in the sense that it sends sinks to sources and vice versa; we then say that $[D_1]$ is \ul{\em composable} with $[D_2]$. Now, say that $D_1$ is a boundary-ordered oriented tangle diagram in a directed biangle $\vec{B}_1$, and $D_2$ is one in $\vec{B}_2$. Let $\vec{B}'$ be the directed biangle obtained by gluing $\vec{B}_1$ and $\vec{B}_2$ along the inward boundary arc of $\vec{B}_1$ and the outward boundary arc of $\vec{B}_2$ respecting their orientations, so that the gluing map restricts to the above bijection $\partial_{\rm in} D_1 \to \partial_{\rm out}D_2$. Let $D$ be the boundary-ordered oriented tangle diagram in $\vec{B}'$ obtained as the union of the images of $D_1$ and $D_2$ under the gluing map; see Fig.\ref{fig:tensor_product_and_composition_of_tangle_diagrams_in_a_directed_biangle}. Define the \ul{\em composition} of the equivalence classes of boundary-ordered oriented tangle diagrams $[D_1]$ and $[D_2]$ in a directed biangle as $[D_1] \circ [D_2] := [D]$.
\end{definition}

It is well known that any oriented tangle diagram in a directed biangle $\vec{B}$ can be obtained from certain elementary ones by applying a sequence of tensor products and compositions; we formulate the result in terms of boundary-ordered oriented tangle diagrams.

\begin{definition}[elementary tangle diagrams in a directed biangle \cite{O}]
\label{def:elementary_tangle_diagrams_in_a_directed_biangle}
$\bullet$ A boundary-ordered oriented tangle diagram $D$ in a directed biangle $\vec{B}$ is said to be \ul{\em elementary} if it is one of the ten cases in Fig.\ref{fig:elementary_bo_oriented_tangle_diagrams_in_a_directed_biangle}. We say the corresponding equivalence class $[D]$ is elementary.
\end{definition}
\begin{lemma}
Any $[D] \in \mathcal{D}(\vec{B})$ can be obtained from elementary classes in $\mathcal{D}(\vec{B})$ by applying a sequence of tensor products and compositions.
\end{lemma}

We are now ready to formulate Gabella's biangle quantum holonomy as matrix elements of an operator invariant. From now on let
\begin{align}
\label{eq:V}
V := \mbox{the $2$-dimensional $\mathbb{C}$ vector space with the ordered basis $\{\xi_+, \xi_-\}$}.
\end{align}
We note that, in Gabella's notation \cite{G}, the basis vectors $\xi_+$ and $\xi_-$ correspond to the symbols $s_1$ and $s_2$. For $[D]\in \mathcal{D}(\vec{B})$ we will consider $V^{\otimes |\partial_{\rm in} D|}$ and $V^{\otimes |\partial_{\rm out} D|}$, where the tensor factors are thought to be ordered according to the (increasing) horizontal orderings on $\partial_{\rm in}D$ and $\partial_{\rm out} D$; for example, in case of Fig.\ref{fig:elementary_tangle_diagrams_in_a_directed_biangle}(IV), we have $V^{\otimes|\partial_{\rm in} D|} = V\otimes V$, with the first factor corresponding to the endpoint $x_1$ and the second factor to $x_2$. We set $V^{\otimes0} := \mathbb{C}$.
\begin{proposition}[Gabella's biangle quantum holonomy packaged as an operator invariant]
\label{prop:G_D}
There exists a unique map $G$ assigning to each $[D] \in \mathcal{D}(\vec{B})$ a linear operator
$$
G([D]) : V^{\otimes|\partial_{\rm in} D|} \to V^{\otimes |\partial_{\rm out}D|}
$$
that satisfies
$$
G([D_1]\circ [D_2]) = G([D_1])\circ G([D_2]), \qquad
G([D_1]\otimes [D_2]) = G([D_1]) \otimes G([D_2]),
$$
and whose values at elementary classes $[D] \in \mathcal{D}(\vec{B})$ are:
\begin{enumerate}
\item[\rm (GB1)] (identity) if $D$ is as in Fig.\ref{fig:elementary_bo_oriented_tangle_diagrams_in_a_directed_biangle}(I) or (VI), then $G([D]) = {\rm id} : V \to V$;

\item[\rm (GB2)] (cup) if $D$ is as in Fig.\ref{fig:elementary_bo_oriented_tangle_diagrams_in_a_directed_biangle}(II) or (VII), then $G([D]) : \mathbb{C} \to V\otimes V$ sends $1 \in \mathbb{C}$ to $\xi_+ \otimes \xi_- - \omega^4 \xi_- \otimes \xi_+$;

\item[\rm (GB3)] (cap) if $D$ is as in Fig.\ref{fig:elementary_bo_oriented_tangle_diagrams_in_a_directed_biangle}(III) or (VIII), then $G([D]):V\otimes V \to \mathbb{C}$ sends $\xi_+ \otimes \xi_-$ to $-\omega^{-4} $ and $\xi_- \otimes \xi_+$ to $1$ while sending other basis vectors to zero;

\item[\rm (GB4)] (height exchange) if $D$ is as in Fig.\ref{fig:elementary_bo_oriented_tangle_diagrams_in_a_directed_biangle}(IV), then the map $G([D]) :V \otimes V \to V\otimes V$ is given on the basis vectors as:
\begin{align}
\label{eq:GB4}
G([D]) : \xi_i \otimes \xi_j \mapsto \left\{
\begin{array}{ll}
\xi_+ \otimes \xi_- + (\omega^4 - \omega^{-4})\xi_-\otimes \xi_+, & \mbox{if $i=+$ and $j=-$,} \\	
\xi_i \otimes \xi_j, & \mbox{otherwise,}
\end{array}
\right.
\end{align}
while if $D$ is as in Fig.\ref{fig:elementary_bo_oriented_tangle_diagrams_in_a_directed_biangle}(IX), then the map $G([D]):V\otimes V\to V\otimes V$ is given as the inverse of the map in eq.\eqref{eq:GB4};

\item[\rm (GB5)] (positive crossing) if $D$ is as in Fig.\ref{fig:elementary_bo_oriented_tangle_diagrams_in_a_directed_biangle}(V), then $G([D]) : V\otimes V \to V\otimes V$ is given on the basis vectors as:
\begin{align}
\label{eq:GB5}
G([D]) : \xi_i \otimes \xi_j \mapsto \left\{
\begin{array}{ll}
\omega^{-4} \xi_i \otimes \xi_j, & \mbox{if $i= j$,} \\
\xi_- \otimes \xi_+, & \mbox{if $i=+$ and $j=-$,} \\	
\xi_+ \otimes \xi_- + (\omega^{-4}-\omega^4) \xi_- \otimes \xi_+, & \mbox{if $i=-$ and $j=+$;}
\end{array}
\right.
\end{align}

\item[\rm (GB6)] (negative crossing) if $D$ is as in Fig.\ref{fig:elementary_bo_oriented_tangle_diagrams_in_a_directed_biangle}(X), then $G([D]) : V\otimes V \to V\otimes V$ is given as the inverse of the map in eq.\eqref{eq:GB5}.
\end{enumerate}

\end{proposition}
What will be actually used in the construction of Gabella's quantum holonomy over the entire triangulated surface are matrix elements.

\begin{definition}[Basis of $V^{\otimes M}$ and bilinear form]
\label{def:basis_of_tensor_space_and_bilinear_form}
For each ordered sequence $\vec{\varepsilon} = \{ \varepsilon_1,\ldots,\varepsilon_M\}$ in $\{+,-\}$, define the basis vector of $V^{\otimes M}$ as
$$
\xi^{\vec{\varepsilon}} := \xi_{\varepsilon_1} \otimes \xi_{\varepsilon_2} \otimes \cdots \otimes \xi_{\varepsilon_M} \in V\otimes V\otimes \cdots \otimes V = V^{\otimes M}.
$$
Define a symmetric bilinear form $\langle \cdot,\cdot\rangle$ on $V^{\otimes M}$ as 
$\langle \xi^{\vec{\varepsilon}}, \xi^{\vec{\varepsilon}\,'} \rangle = \delta_{\vec{\varepsilon},\vec{\varepsilon}\,'}$, the Kronecker delta.
\end{definition}

\begin{definition}[Gabella's biangle quantum holonomy as matrix elements of $G$]
\label{def:Gabella_biangle_quantum_holonomy_as_matrix_elements}
Let $\vec{B}$ be a directed biangle, with the underlying biangle being $B = (\Sigma,\mathcal{P})$, with $\frak{S} = \Sigma\setminus\mathcal{P}$. Let $\omega\in \mathbb{C}^*$. Let $[K,s]$ be the VH-isotopy class of a stated oriented VH-tangle in $\frak{S}\times(-1,1)$.

\vs

Denote by $[D,s] \in \mathcal{D}_{\rm s}(\vec{B})$ the stated boundary-ordered oriented tangle diagram of $[K,s]$ and let \\$\partial_{\rm in}(D) = \{x_1,x_2,\ldots,x_{|\partial_{\rm in}D|}\}$ and $\partial_{\rm out}(D) = \{y_1,\ldots,y_{|\partial_{\rm out} D|}\}$, arranged in the respective (increasing) horizontal orderings. Define the \ul{\em inward/outward basis vector for $s$} as
\begin{align}
\label{eq:inward_and_outward_basis_vectors}
\xi_{\rm in}^s := \xi_{s(x_1)} \otimes \xi_{s(x_2)} \otimes \cdots \otimes \xi_{s(x_{|\partial_{\rm in}D|})} \in  V^{\otimes|\partial_{\rm in}D|}, \quad
\xi_{\rm out}^s := \xi_{s(y_1)} \otimes \cdots \otimes \xi_{s(y_{|\partial_{\rm out}D|)}} \in V^{\otimes|\partial_{\rm out}D|},
\end{align}
where in case $\partial_{\rm in} D= {\O}$ or $\partial_{\rm out} D= {\O}$, we let $\xi^s_{\rm in}=1$ or $\xi^s_{\rm out}=1$ respectively. Let $G([D,s])\in \mathbb{C}$ be the matrix element of the linear map $G([D]) : V^{\otimes|\partial_{\rm in}D|} \to V^{\otimes|\partial_{\rm out}D|}$ of Prop.\ref{prop:G_D} for these basis vectors, i.e.
$$
G([D,s]) := \langle \xi^s_{\rm out}, \, G([D]) \xi^s_{\rm in}  \rangle
$$
Define the \ul{\em Gabella biangle quantum holonomy} as the map ${\rm TrHol}^\omega_{\vec{B}}$ assigning to each $[K,s]$ the value:
$$
{\rm TrHol}_{\vec{B}}^\omega([K,s]) := G([D,s]) \in \mathbb{C}.
$$
\end{definition}
For example, if $D$ is as in Prop.\ref{prop:G_D}(GB2) with $y_1\succ y_2$, $s(y_1)=-$ and $s(y_2)=+$, then ${\rm TrHol}_{\vec{B}}^\omega([K,s]) = -\omega^4$, while if $D$ is as in Prop.\ref{prop:G_D}(GB4) with $s(x_1)=+$, $s(x_2)=-$, $s(y_1)=-$ and  $s(y_2)=+$, then ${\rm TrHol}_{\vec{B}}^\omega([K,s]) = \omega^4-\omega^{-4}$.

\vs

The following useful observations can be checked in a straightforward manner; we will provide a proof of it at the end of \S\ref{subsec:Gabella_biangle_factors_as_twisting_of_BW_biangle_factors}, using the relationship with the Bonahon-Wong biangle quantum trace.
\begin{lemma}
\label{eq:G_D_useful_properties}
The Gabella biangle quantum holonomy satisfies:
\label{lem:simple_properties_of_G_D}
\begin{enumerate}
\item[\rm (1)] (direction independence) ${\rm TrHol}_{\vec{B}}^\omega([K,s])$ is independent of the choice of the direction on $B$, so we can define
$$
{\rm TrHol}^\omega_B([K,s]) := {\rm TrHol}^\omega_{\vec{B}}([K,s]).
$$

\item[\rm (2)] (charge conservation) A state $s$ of a tangle diagram $D$ in $\vec{B}$ or of a VH-tangle in the thickening of $\vec{B}$ is said to be \ul{\em charge-preserving} if and only if the net sum of signs assigned by $s$ to elements of $\partial_{\rm in} D$  equals that for $\partial_{\rm out}D$ (where $\pm$ are viewed as $\pm 1$). We have ${\rm TrHol}_{\vec{B}}^\omega([K,s])=0$ for any state $s$ that is not charge-preserving.
\end{enumerate}
\end{lemma}
Above is a rigorous treatment via operator invariants of what Gabella meant in \cite{G}. In particular, one can verify by inspection that the values of ${\rm TrHol}_B^\omega$ at the elementary stated oriented VH-tangles defined as above coincide with Gabella's values in \cite[\S5]{G} except for the cups and caps. The values at cups and caps originally assigned by Gabella in \cite[eq.(5.7)]{G} should be corrected to our values given by Def.\ref{def:Gabella_biangle_quantum_holonomy_as_matrix_elements} and Prop.\ref{prop:G_D}(GB2)--(GB3), in order to guarantee the well-definedness (i.e. isotopy invariance) of Gabella's quantum holonomy for the entire surface which will be constructed in the next subsection. We will give a full justification of our construction of ${\rm TrHol}_B^\omega$, including the proof of Prop.\ref{prop:G_D}, later in the next section.

\subsection{Quantum holonomy}

We now describe Gabella's quantum holonomy \cite{G} associated to an oriented VH-tangle $K$ in $\frak{S} \times (-1,1)$. We first define a Gabella lift of a stated oriented VH-tangle $(K,s)$, and enumerate their isotopy classes by the collection of states for the Gabella lifts of triangle segments, which we package as a $\wh{\Delta}$-juncture-state of $K$.
\begin{definition}
\label{def:Gabella_lift_of_K_s}
Let $(\Sigma,\mathcal{P})$, $\frak{S}$, $\Delta$ and $\wh{\Delta}$ be as in Def.\ref{def:Gabella_lift}. Let $K$ be an oriented VH-tangle in $\frak{S}\times(-1,1)$ in a good position with respect to $\wh{\Delta}$. Let $\til{K}$ be a Gabella lift in $\til{\Sigma}_\Delta\times(-1,1)$ of $K$. The $\wh{\Delta}$-juncture-state $J_{\til{K}}$ of $K$ whose restriction to the endpoints of each triangle segment $k$ of $K$ coincides with the state $s_{\til{k}}$ of the Gabella lift $\til{k}$ of $k$ determined by $\til{K}$ is the \ul{\em $\wh{\Delta}$-juncture-state of $K$ associated to Gabella lift $\til{K}$}.

\vs

If $s$ is a state of $K$, by a \ul{\em Gabella lift of the stated oriented VH-tangle $(K,s)$} we mean a Gabella lift $\til{K}$ of $K$ whose associated $\wh{\Delta}$-juncture-state $J_{\til{K}}$ restricts to $s$ at $\partial K$.
\end{definition}

\begin{lemma}
Let $(K,s)$ be a stated oriented VH-tangle in $\frak{S}\times(-1,1)$ in a good position with respect to $\wh{\Delta}$. The correspondence $\til{K} \mapsto J_{\til{K}}$ yields a bijection between the set of all isotopy classes of Gabella lifts of $K$ and the set of all $\wh{\Delta}$-juncture-states of $K$. Gabella lifts of $(K,s)$ correspond to $\wh{\Delta}$-juncture-states of $K$ restricting to $s$ at $\partial K$. 
\end{lemma}
\begin{definition}
For a $\wh{\Delta}$-juncture-state $J$ of $K$, denote by $\til{K}^J$ a Gabella lift of $K$ whose associated $\wh{\Delta}$-juncture-state coincides with $J$.

\vs

A $\wh{\Delta}$-juncture-state $J$ of $K$ is said to be \ul{\em admissible} if $\til{K}^J$ is an admissible Gabella lift of $K$ (see Def.\ref{def:Gabella_lift}), that is, if none of the triangle segments of $K$ fall into the case of Fig.\ref{fig:elementary_skeins_in_the_triangle}(a) with $J(x_1)=-$ and $J(x_2)=+$. 
\end{definition}

\begin{remark}
Notice that the state-sum formula of the Bonahon-Wong quantum trace, as seen in \S\ref{subsec:state-sum_model}, was written as a sum over all $\wh{\Delta}$-juncture-states of $K$. One may observe that the terms corresponding to nonadmissible $\wh{\Delta}$-juncture-states are zero.
\end{remark}

We are now ready to state our version of  the construction of the sought-for Gabella quantum holonomy.

\begin{definition}[enhanced and modified version of Gabella quantum holonomy]
\label{def:Gabella_quantum_holonomy}
Let $(\Sigma,\mathcal{P})$, $\frak{S}$, $\Delta$ and $\wh{\Delta}$ be as in Def.\ref{def:Gabella_lift}. Let $\omega\in \mathbb{C}^*$, and let $q = \omega^4$. Let $(K,s)$ be a stated oriented VH-tangle in $\frak{S} \times (-1,1)$. 

\vs

Isotope $(K,s)$ into a stated oriented VH-tangle $(K',s')$ in $\frak{S}\times(-1,1)$ in a good position with respect to $\wh{\Delta}$ (see Lem.\ref{lem:good_position}) through a VH-isotopy.

\vs

For each Gabella lift $\til{K}^J$ in $\til{\Sigma}_\Delta \times (-1,1)$ of $(K',s')$ associated to a $\wh{\Delta}$-juncture-state $J : \{\mbox{$\wh{\Delta}$-junctures of $K'$}\} \to \{+,-\}$ of $K'$ that restricts to $s'$ at $\partial K'$, we define a monomial $\wh{Z}_{\til{K}^J}$ in the (square-root) generators $\wh{Z}_1,\ldots,\wh{Z}_n$ of the Chekhov-Fock algebra $\mathcal{T}^\omega_\Delta$ (see Def.\ref{def:CF_balanced_square-root_algebra}) and a coefficient $\ol{\ul{\Omega}}(\til{K}^J;\omega) \in \mathbb{Z}[\omega,\omega^{-1}]$ as follows:
\begin{enumerate}
\item[\rm (G1)] (the monomial part) For each edge $e$ of the original triangulation $\Delta$, denote again by $e$ one of any of the two edges $e$ and $e'$ in $\wh{\Delta}$ corresponding to $e$. Let
$$
b_e(J) \in \mathbb{Z}
$$
be the sum of signs of all $\wh{\Delta}$-junctures of $K'$ over $e$ assigned by the $\wh{\Delta}$-juncture-state $J$, where $+$ and $-$ are thought of as $1$ and $-1$ respectively. Now let
$$
\wh{Z}_{\til{K}^{\hspace{-0,3mm}J}} := \left[\prod_{e\in \Delta} \wh{Z}_e^{\,b_e\hspace{-0,3mm}(J)}\right]_{\rm Weyl} \quad \in \quad \mathcal{T}^\omega_\Delta,
$$
where $[\sim]_{\rm Weyl}$ is the Weyl-ordered monomial defined in Def.\ref{def:Weyl-ordered_product}.

\vs

\begin{figure}
\centering 
\hspace{5mm} \scalebox{0.8}{
\begingroup%
  \makeatletter%
  \providecommand\color[2][]{%
    \errmessage{(Inkscape) Color is used for the text in Inkscape, but the package 'color.sty' is not loaded}%
    \renewcommand\color[2][]{}%
  }%
  \providecommand\transparent[1]{%
    \errmessage{(Inkscape) Transparency is used (non-zero) for the text in Inkscape, but the package 'transparent.sty' is not loaded}%
    \renewcommand\transparent[1]{}%
  }%
  \providecommand\rotatebox[2]{#2}%
  \newcommand*\fsize{\dimexpr\f@size pt\relax}%
  \newcommand*\lineheight[1]{\fontsize{\fsize}{#1\fsize}\selectfont}%
  \ifx\svgwidth\undefined%
    \setlength{\unitlength}{226.77165354bp}%
    \ifx\svgscale\undefined%
      \relax%
    \else%
      \setlength{\unitlength}{\unitlength * \real{\svgscale}}%
    \fi%
  \else%
    \setlength{\unitlength}{\svgwidth}%
  \fi%
  \global\let\svgwidth\undefined%
  \global\let\svgscale\undefined%
  \makeatother%
  \begin{picture}(1,0.5)%
    \lineheight{1}%
    \setlength\tabcolsep{0pt}%
    \put(0,0){\includegraphics[width=\unitlength,page=1]{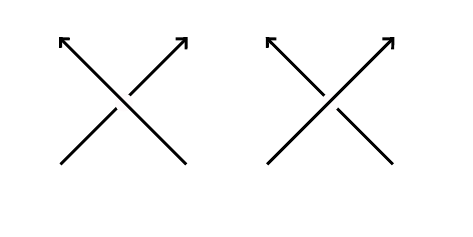}}%
    \put(0.23977864,0.08764329){\color[rgb]{0,0,0}\makebox(0,0)[lt]{\lineheight{1.25}\smash{\begin{tabular}[t]{l}$(-)$\end{tabular}}}}%
    \put(0.6796484,0.08864274){\color[rgb]{0,0,0}\makebox(0,0)[lt]{\lineheight{1.25}\smash{\begin{tabular}[t]{l}$(+)$\end{tabular}}}}%
  \end{picture}%
\endgroup%
}
\vspace{-7mm}
\caption{Sign of a crossing}
\vspace{-5mm}
\label{fig:sign_of_crossing}
\end{figure}

\item[\rm (G2)] (the $q$-power coefficient part) To each of the crossings of $\til{K}^J$ (as defined at the end of \S\ref{subsec:branched_double_cover}), associate the sign $+1$ if it is of type $(+)$ in Fig.\ref{fig:sign_of_crossing}, and $-1$ if it is of type $(-)$ in Fig.\ref{fig:sign_of_crossing}. Let 
$$
{\rm wr}_{\wh{\Delta}}(\til{K}^J) \in \mathbb{Z}
$$
be the net sum of all signs of crossings over all triangles of $\wh{\Delta}$, i.e. the usual \ul{\em writhe} of the tangle $\til{K}^J$ considered only over triangles of $\wh{\Delta}$.

\vs

\item[\rm (G3)] (the biangle factor; ``R-matrix" and cup/cap) For each $i=1,\ldots,n$, let $L_i := K' \cap ((\Sigma_i\setminus\mathcal{P}_i)\times(-1,1))$ be the part of $K'$ over the biangle $B_i =(\Sigma_i,\mathcal{P}_i)$ of $\wh{\Delta}$ corresponding to the edge $e_i$ of $\Delta$. Then $(L_i, J|_{\partial L_i})$ is a stated oriented VH-tangle in $(\Sigma_i\setminus\mathcal{P}_i)\times(-1,1)$. Consider its value under the biangle quantum holonomy map ${\rm TrHol}^\omega_{B_i}([L_i,J|_{\partial L_i}]) \in \mathbb{C}$.
\end{enumerate}
Let
\begin{align}
\label{eq:Gabella_coefficient}
\ol{\ul{\Omega}}(\til{K}^J;\omega) := \textstyle (\prod_{i=1}^n {\rm TrHol}^\omega_{B_i}([L_i,J|_{\partial L_i}]))  \cdot q^{-{\rm wr}_{\wh{\Delta}}(\til{K}^J)} \quad \in \quad \mathbb{Z}[\omega,\omega^{-1}].
\end{align}

\vs

Define the (enhanced) \ul{\em Gabella quantum holonomy} for the stated oriented VH-tangle $(K,s)$ in $\frak{S}\times(-1,1)$ by the formula
\begin{align}
\label{def:TrHol_sum_formula}
{\rm TrHol}^\omega_\Delta(K,s) := \sum_{\til{K}^{\hspace{-0,3mm}J}} \ol{\ul{\Omega}}(\til{K}^J;\omega) \, \wh{Z}_{\til{K}^{\hspace{-0,3mm}J}} \quad \in \quad \mathcal{T}^\omega_\Delta
\end{align}
where the sum is over all isotopy classes of {\em admissible} Gabella lifts $\til{K}^J$ of $(K',s')$ (see Def.\ref{def:Gabella_lift}), i.e. over all admissible $\wh{\Delta}$-juncture-states $J$ of $K'$ that restrict to $s'$ at $\partial K'$.
\end{definition}
One might wonder why we can choose any one of $e$ and $e'$ in (G1) above. We claim that for any $\wh{\Delta}$-juncture-state $J$ of $K'$ for which the coefficient $\ol{\ul{\Omega}}(\til{K}^J;\omega)$ is nonzero, we have $b_e(J) = b_{e'}(J)$; this is an easy consequence of the charge conservation property (see Lem.\ref{lem:simple_properties_of_G_D}(2)) of the biangle quantum holonomy. The ideas of (G1) and (G2) already appeared in \cite{GLM} for a connected nonclosed tangle that is given the always-going-up elevations, and as mentioned already, the main contribution of Gabella \cite{G} is the biangle factor (G3), which enables one to define the quantum holonomy for closed tangles. Meanwhile, as mentioned in \S\ref{subsec:biangle_quantum_holonomy}, the biangle quantum holonomy was only partially dealt with in \cite{G}, and the independence of the value of ${\rm TrHol}^\omega_\Delta(K,s)$ (even in case $K={\O}$) on the choice made in the construction, namely the choice of an isotopy transformation into a tangle $K'$ in a good position, was only partially proved in \cite{G}. This independence was made as a conjecture in the earlier version of the present paper \cite{KSv1}, and will be obtained as a consequence of the main result of this paper, in a fully general case of a stated oriented VH-tangle.

\vs

In the meantime, one actual difference between Gabella's original construction \cite{G} and ours is the monomial part (G1); Gabella uses different normalization. Namely, denote by $|b_e|(K)$ be the number of all junctures of $K$ on $e$. Then, in place of our monomial $\wh{Z}_{\til{K}^{\hspace{-0,3mm}J}}$, Gabella uses
\begin{align}
\label{eq:Gabella_normalization}
\wh{X}_{\til{K}^{\hspace{-0,3mm}J}} := \left[\prod_{e\in \Delta} \wh{X}_e^{\,\frac{1}{2}(b_e\hspace{-0,3mm}(J)+|b_e|(K))}\right]_{\rm Weyl} \quad \in \quad \mathcal{T}^q_\Delta.
\end{align}
Gabella's choice causes the final result to enjoy many properties, and has the advantage of avoiding the square-root variables (since $\frac{1}{2}(b_e\hspace{-0,3mm}(J)+|b_e|(K))\in\mathbb{Z}$) and making the lowest term $1$. Later, we will discuss why we had to modify as we did.

\vs

As already mentioned, Gabella's original construction is for a simple loop in the surface $\frak{S}$.
\begin{definition}
\label{def:Gabella_K_gamma}
Let $\Delta$ be a triangulation of a triangulable generalized marked surface $(\Sigma,\mathcal{P})$, and let $\gamma$ be an oriented  simple loop in $\frak{S}=\Sigma\setminus\mathcal{P}$.

\vs

Let $K_\gamma$ be the oriented tangle in $\frak{S}\times(-1,1)$ obtained as a constant-elevation lift of $\gamma$, as in Def.\ref{def:K_gamma}. Define the Gabella quantum holonomy of the oriented simple loop $\gamma$ in the surface $\frak{S}$ as
\begin{align}
\label{eq:Gabella_solution}
{\rm TrHol}^\omega_\Delta(\gamma) := {\rm TrHol}^\omega_\Delta(K_\gamma,{\O}) \in \mathcal{T}^\omega_\Delta.
\end{align}
\end{definition}
Some of the important properties of ${\rm TrHol}^\omega_\Delta(\gamma)$ proved by Gabella in \cite{G}  read as follows, when translated to our normalization:
\begin{proposition}[properties of Gabella quantum holonomy \cite{G}]
\label{prop:properties_of_Gabella}
Let $(\Sigma,\mathcal{P})$, $\frak{S}$, $\Delta$, $\gamma$ and $K=K_\gamma$ be as in Def.\ref{def:Gabella_K_gamma}. Suppose further that $\gamma$ is not a contractible loop in $\frak{S}$.
\begin{enumerate}
\item[\rm (1)] When $\omega=1$, ${\rm TrHol}^1_\Delta(\gamma)$ recovers the (absolute value of the) usual trace-of-holonomy function $\mathbb{I}(\gamma)$ on $\mathcal{X}^+(\Sigma,\mathcal{P})$ along $\gamma$;

\item[\rm (2)] In terms of the partial ordering on the monomials $\wh{Z}_{\til{K}^{\hspace{-0,3mm}J}} \in \mathcal{T}^\omega_\Delta$ appearing in the right hand side of the sum formula eq.\eqref{def:TrHol_sum_formula} for ${\rm TrHol}^\omega_\Delta(\gamma) = {\rm TrHol}^\omega_\Delta(K,{\O})$ induced by the powers of the generators $\wh{Z}_1,\ldots,\wh{Z}_n$, there is a unique Gabella lift $\til{K}^{\hspace{-0,3mm}J}$ giving the highest term, for which the power $a_e$ of each generator $\wh{Z}_e$ in the monomial $\wh{Z}_{\til{K}^{\hspace{-0,3mm}J}}$, associated to the edge $e$ of $\Delta$, equals the intersection number $a_e(\gamma)$ defined in Lem.\ref{lem:a_e} (two times the Fock-Goncharov tropical $\mathcal{A}$-coordinate for the lamination $\gamma$), and the coefficient $\ol{\ul{\Omega}}(\til{K}^{\hspace{-0,3mm}J};\omega)$ is $1$;

\item[\rm (3)] For each $\til{K}^{\hspace{-0,3mm}J}$ appearing in eq.\eqref{def:TrHol_sum_formula}, the coefficient $\ol{\ul{\Omega}}(\til{K}^{\hspace{-0,3mm}J};\omega)$ belongs to $\mathbb{Z}_{\ge 0}[q,q^{-1}]$, i.e. is a {\em positive} integral Laurent polynomial in $q=\omega^4$; 

\item[\rm (4)] For each $\til{K}^{\hspace{-0,3mm}J}$ appearing in eq.\eqref{def:TrHol_sum_formula}, the coefficient $\ol{\ul{\Omega}}(\til{K}^{\hspace{-0,3mm}J};\omega)$ is $*$-invariant, in the sense that $\ol{\ul{\Omega}}(\til{K}^{\hspace{-0,3mm}J};\omega) = \ol{\ul{\Omega}}(\til{K}^{\hspace{-0,3mm}J};\omega^{-1})$, or that it is invariant under the exchange $q\leftrightarrow q^{-1}$.
\end{enumerate}
\end{proposition}
One remark is that the proof given in \cite[\S6.4]{G} of the positivity property of part (3) is not quite sufficient; we note that this positivity will follow from the main result of our paper along with \cite{CKKO}.

\section{Equality of the two constructions}

\subsection{Statement of the main theorem} For a stated oriented VH-tangle $(K,s)$ in $\frak{S}\times(-1,1)$, we investigated the Bonahon-Wong quantum trace ${\rm Tr}^\omega_{\Delta}([K,s]) \in \mathcal{T}^\omega_\Delta$ (by viewing $(K,s)$ as a stated oriented V-tangle) and Gabella's quantum holonomy ${\rm TrHol}^\omega_\Delta([K,s])\in \mathcal{T}^\omega_\Delta$. They have several properties in common, which naturally leads to the question of whether they are equal, as mentioned in \cite{G} for the case when $K$ is closed and connected. This equality is our main theorem. In order to formulate the result in full generality, incorporating possibly nonclosed tangles, we first need to introduce:

\begin{definition}[signed order correction amount]
\label{def:signed_order_correction_amount}
Let $(\Sigma,\mathcal{P})$ be a generalized marked surface, not necessarily triangulable. Let $\frak{S} = \Sigma\setminus\mathcal{P}$. Let $e$ be a $\mathcal{P}$-arc (see Def.\ref{def:P-link}) of $(\Sigma,\mathcal{P})$, with an orientation chosen on $e$. 

\vs

$\bullet$ Let $Z$ be an ordered finite subset of $e$, and $s_0:Z\to \{+,-\}$ be a state of $Z$. Define the \ul{\em signed order} \ul{\em correction amount} of the stated ordered subset $(Z,s_0)$ on the oriented arc $e$ as the integer
\begin{align*}
\mathcal{C}(e;Z,s_0) := \sum_{x,y\in Z \, | \, x<y} {\rm sgn}(e;\overrightarrow{xy}) \, s_0(x) \, s_0(y) ~\in~\mathbb{Z},
\end{align*}
where each sign $\pm$ is regarded as the number $\pm 1$, and
\begin{align*}
{\rm sgn}(e;\overrightarrow{xy}) := \left\{
\begin{array}{rl}
1, & \mbox{if the direction from $x$ to $y$ matches the given orientation on $e$,} \\
-1, & \mbox{otherwise.}
\end{array}
\right.
\end{align*}
When $Z={\O}$, we set $\mathcal{C}(e;{\O},{\O})=0$.

\vs

$\bullet$ When $Z$ is a finite subset of the thickening $e \times (-1,1)$ of $e$ equipped with a state $s_0:Z\to\{+,-\}$, such that the restrictions to $Z$ of the projection maps $e\times (-1,1) \to e$ and $e\times (-1,1) \to (-1,1)$ are injective, we define $\mathcal{C}(e;Z,s_0)$ using the projection of $(Z,s_0)$ onto $e$, with the ordering on the projection of $Z$ coming from the elevations of the elements of $Z$, i.e. the vertical ordering. 

\vs

$\bullet$ When $e$ is a boundary arc and no orientation is specified, we use the boundary-orientation (see Def.\ref{def:surface}) coming from the surface orientation on $\Sigma$.

\vs

$\bullet$ If $(K,s)$ is a stated VH-tangle in $\frak{S} \times (-1,1)$, the \ul{\em boundary signed order correction amount} of $(K,s)$ is defined as
$$
\partial \mathcal{C}_{(\Sigma,\mathcal{P})}(K,s) := \sum_{\mbox{\tiny boundary arcs $b$ of $(\Sigma,\mathcal{P})$}} \mathcal{C}(b; \partial K \cap (b\times(-1,1)), \, s|_{\partial K \cap (b\times(-1,1))})
$$
where we always use the boundary-orientation on $b$'s even if $(\Sigma,\mathcal{P})$ is a directed biangle.
\end{definition}

\vs

\begin{theorem}[main theorem]
\label{thm:main}
Let $(\Sigma,\mathcal{P})$ be a triangulable generalized marked surface, and $\Delta$ a triangulation of $(\Sigma,\mathcal{P})$. Let $\frak{S} = \Sigma\setminus\mathcal{P}$. Let $\omega \in \mathbb{C}^*$, $q = \omega^4$ and $A = \omega^{-2}$. Let $(K,s)$ be a stated oriented VH-tangle in $\frak{S} \times (-1,1)$\, satisfying (P1) and (P2) of \S\ref{label:tangles_and_skein_algebra}. Then the Bonahon-Wong quantum trace ${\rm Tr}^\omega_\Delta([K,s]) \in \mathcal{Z}^\omega_\Delta \subset \mathcal{T}^\omega_\Delta$, constructed in \S\ref{sec:BW} for the stated skein $[K,s] \in \mathcal{S}^A_{\rm s}(\Sigma,\mathcal{P})$, is related to the enhanced Gabella quantum holonomy ${\rm TrHol}^\omega_\Delta(K,s) \in \mathcal{T}^\omega_\Delta$ constructed in 
\S\ref{sec:Gabella} as
\begin{align}
\label{eq:main_thm}
{\rm Tr}^\omega_\Delta([K,s]) = \omega^{2{\rm wr}(K)} \, \omega^{\partial \mathcal{C}_{(\Sigma,\mathcal{P})}(K,s)} \, {\rm TrHol}^\omega_\Delta(K,s),
\end{align}
where the writhe ${\rm wr}(K)\in\mathbb{Z}$ of the VH-tangle $K$ is defined as
$$
{\rm wr}(K) := \#\mbox{(crossings of type $(+)$ in Fig.\ref{fig:sign_of_crossing})} - \#\mbox{(crossings of type $(-)$ in Fig.\ref{fig:sign_of_crossing})},
$$
counted over $\Sigma$ (not over the branched double cover $\til{\Sigma}_\Delta$).

\vs

In particular, this equality holds for a constant-elevation lift $K = K_\gamma$ of an oriented simple loop $\gamma$ in $\frak{S}$ (see Def.\ref{def:K_gamma}), where ${\rm wr}(K)=0$, $\partial K={\O}$, $s={\O}$ and $\partial \mathcal{C}_{(\Sigma,\mathcal{P})}(K,{\O})=0$. Hence, Allegretti and Kim's solution $\wh{\mathbb{I}}^\omega_\Delta(\gamma)$ in eq.\eqref{eq:AK_solution} to the quantum ordering problem for the classical function $\mathbb{I}(\gamma)$ addressed in \S\ref{subsec:deformation_quantization} coincides with Gabella's solution ${\rm TrHol}^\omega_\Delta(\gamma)$ in eq.\eqref{eq:Gabella_solution}.
\end{theorem}
This theorem has some immediate consequences on the Gabella quantum holonomy ${\rm TrHol}^\omega_\Delta(K,s)$ following from the corresponding properties of the Bonahon-Wong quantum trace ${\rm Tr}^\omega_\Delta([K,s])$, and some obvious observations about ${\rm wr}(K)$ and $\partial \mathcal{C}_{(\Sigma,\mathcal{P})}(K,s)$. Namely, we obtain a proof of Prop.\ref{prop:properties_of_Gabella}, especially the positivity in part (3), and the following corollaries:
\begin{corollary}[isotopy invariance]
\label{cor:isotopy_invariance}
The enhanced Gabella quantum holonomy defined in Def.\ref{def:Gabella_quantum_holonomy}, ${\rm TrHol}^\omega_\Delta(K,s)$, depends only on the VH-isotopy class of a stated oriented VH-tangle $(K,s)$.
\end{corollary}

\begin{corollary}[mutation compatibility]
\label{corollary:Gabella_mutation_compatibility}
Let $(\Sigma,\mathcal{P})$ and $\Delta$ be as in Def.\ref{def:good_position}. Let $\frak{S} = \Sigma\setminus\mathcal{P}$. Let $\omega\in \mathbb{C}^*$ and $q = \omega^4$. Let $\Delta'$ be another triangulation of $(\Sigma,\mathcal{P})$.

\vs

Let $(K,s)$ be an oriented stated VH-tangle in $\frak{S}\times(-1,1)$. Then we have
$$
{\rm TrHol}^\omega_\Delta(K,s) = \Theta^\omega_{\Delta\Delta'} ( {\rm TrHol}^\omega_{\Delta'}(K,s) ),
$$
where $\Theta^\omega_{\Delta\Delta'}$ is the square-root quantum coordinate change map as in Prop.\ref{prop:square-root_skew-field_isomorphisms}. In particular, \\${\rm TrHol}^\omega_\Delta(K,s) \in \mathcal{Z}^\omega_\Delta$ and  ${\rm TrHol}^\omega_{\Delta'}(K,s) \in \mathcal{Z}^\omega_{\Delta'}$ (see Def.\ref{def:CF_balanced_square-root_algebra}). 

\vs

As a consequence, we also obtain ${\rm TrHol}^\omega_\Delta(K,s) \in {\bf L}^\omega_\Delta$, where ${\bf L}^\omega_\Delta$ is defined as in eq.\eqref{eq:bf_L_omega}.
\end{corollary}

The main theorem and these corollaries hold for the Gabella quantum holonomy ${\rm TrHol}^\omega_\Delta(K,s)$ defined using Prop.\ref{prop:G_D} and the normalization Def.\ref{def:Gabella_quantum_holonomy}(G1) as in our present paper, but not for Gabella's original construction \cite{G} which uses different values for Prop.\ref{prop:G_D}(GB2)--(GB3) and the normalization in eq.\eqref{eq:Gabella_normalization}.

\vs

We note that, in the earlier version \cite{KSv1} of the present paper, the main theorem was proved only for the case when the tangle $K$ is closed, connected, and has no crossing, because the biangle quantum holonomy was only partially dealt with. In the present version, through the contribution of the newly participating second author, we now have a complete treatment of the biangle quantum holonomy, hence the main theorem in full generality. The new idea added in this version also made the proofs shorter and more direct.

\vs

Let's now start proving the main theorem. Let $(K,s)$ be any stated oriented VH-tangle in $\frak{S}\times (-1,1)$. Isotope it to a stated oriented VH-tangle $(K',s')$ in a good position with respect to a split ideal triangulation $\wh{\Delta}$ of $\Delta$ (see Def.\ref{def:good_position}), through a VH-isotopy. Then both sides of sought-for eq.\eqref{eq:main_thm} are expressed as sums over admissible $\wh{\Delta}$-juncture-states $J$ of $K'$ restricting to $s'$ at $\partial K'$; see Prop.\ref{prop:state-sum_of_BW} ad Def.\ref{def:Gabella_quantum_holonomy}. It is enough to show that the summands of the two sides for each $J$ coincide, i.e. to show the \ul{\em term-by-term equality}
\begin{align}
\label{eq:term_equality}
{\rm BW}^\omega_{\wh{\Delta}}(K';J) = \omega^{2{\rm wr}(K)} \, \omega^{\partial \mathcal{C}_{(\Sigma,\mathcal{P})}(K,s)} \, \ol{\ul{\Omega}}(\til{K}^J;\omega) \, \wh{Z}_{\til{K}^{\hspace{-0,3mm}J}}
\end{align}
for each $J$. The rest of the present section is devoted to a proof of the term-by-term equality, eq.\eqref{eq:term_equality}; for convenience, from now on we may identify $K$ with $K'$:
$$
K = K'
$$
Throughout the entire section, we will reserve the symbols $K$, $K'$ and $s$ as such.

\subsection{Bonahon-Wong triangle factors}

We further break down the Bonahon-Wong summand \\ ${\rm BW}^\omega_{\wh{\Delta}}(K';J)={\rm BW}^\omega_{\wh{\Delta}}(K;J)$ in the left hand side of the sought-for eq.\eqref{eq:term_equality} into the product of the biangle factors ${\rm Tr}^\omega_{B_i}([L_i, J|_{\partial L_i}])$ and the triangle factors ${\rm Tr}^\omega_{\wh{t}_j}([K_j, J|_{\partial K_j}])$, as in its very definition in eq.\eqref{eq:BW_term}. In the present subsection we focus on the triangle factors.

\vs

Let $(\Sigma,\mathcal{P})$, $\frak{S}$, $\Delta$, $\wh{\Delta}$, $e_i$, $e_i'$, $B_i$, $t_j$ and $\wh{t}_j$ be as in Def.\ref{def:good_position}. Let $\omega\in\mathbb{C}^*$ and $q = \omega^4$. As done in Prop.\ref{prop:state-sum_of_BW}, we let $K_j := K \cap (\wh{t}_j \times(-1,1))$ be the part of $K$ over the triangle $\wh{t}_j$ of $\wh{\Delta}$, and let $k_{j,1}, \ldots, k_{j,l_j}$ the components of $K_j$ (the triangle segments of $K$ over $\wh{t}_j$) in order of increasing elevation. The triangle factor for this triangle $\wh{t}_j$ is as written in eq.\eqref{eq:triangle-factor}:
\begin{align}
\nonumber
{\rm Tr}^\omega_{\wh{t}_j}([K_j, J|_{\partial K_j}]) = {\rm Tr}^\omega_{\wh{t}_j}( [k_{j,1}, J|_{\partial k_{j,1}}] ) \, {\rm Tr}^\omega_{\wh{t}_j}( [k_{j,2}, J|_{\partial k_{j,2}}] ) \cdots {\rm Tr}^\omega_{\wh{t}_j}( [k_{j,l_j}, J|_{\partial k_{j,l_j}}] ) \in \mathcal{T}^\omega_{\wh{t}_j}.
\end{align}
We will now rewrite the right hand side as the Weyl-ordered monomial times some integer power of $\omega$, which we refer to as the {\em deviation}. 
\begin{definition}
\label{def:deviation_of_BW_traingle_factor}
The \ul{\em deviation} of the Bonahon-Wong triangle factor for triangle $\wh{t}_j$ from the Weyl-ordering, which is associated to the stated V-tangle $(K_j,J|_{\partial K_j})$ in the thickening of $\wh{t}_j$, is the unique integer ${\rm dev}_{\wh{t}_j}(K_j, J|_{\partial K_j})$ such that
$$
{\rm Tr}^\omega_{\wh{t}_j}([K_j, J|_{\partial K_j}]) = \omega^{{\rm dev}_{\wh{t}_j}(K_j,J|_{\partial K_j})} \, \left[ {\rm Tr}^\omega_{\wh{t}_j}( [k_{j,1}, J|_{\partial k_{j,1}}] ) \cdots {\rm Tr}^\omega_{\wh{t}_j}( [k_{j,l_j}, J|_{\partial k_{j,l_j}}] ) \right]_{\rm Weyl}.
$$
In the right hand side we define $[\sim]_{\rm Weyl}$ as $[\wh{Z}_1^{a_1} \cdots \wh{Z}_n^{a_n}]_{\rm Weyl}$, when $\sim \, = \omega^m \wh{Z}_1^{a_1}\cdots \wh{Z}_n^{a_n}$ for some $m,a_1,\ldots,a_n \in \mathbb{Z}$.
\end{definition}
If $J$ is not admissible both sides of the above equation are zero, so the deviation would not be uniquely determined, but we only deal with admissible $J$. We could apply the definition of deviation to more general stated V-tangles satisfying certain conditions, e.g. (GP2)--(GP3) of Def.\ref{def:good_position}. For example, for the stated V-tangle consisting of two of the components of $(K_j,J|_{\partial K_j})$, whenever $1\le r<u\le l_j$ we have
$$
\hspace{-4mm} {\rm Tr}^\omega_{\wh{t}_j} ([k_{j,r}, \, J|_{\partial k_{j,r}}]) \, {\rm Tr}^\omega_{\wh{t}_j}([k_{j,u}, J|_{\partial k_{j,u}}]) = \omega^{{\rm dev}_{\wh{t}_j}(k_{j,r}\cup k_{j,u}, J|_{\partial k_{j,r} \cup \partial k_{j,u}})} \left[ {\rm Tr}^\omega_{\wh{t}_j} ([k_{j,r}, \, J|_{\partial k_{j,r}}]) \, {\rm Tr}^\omega_{\wh{t}_j}([k_{j,u}, J|_{\partial k_{j,u}}]) \right]_{\rm Weyl}.
$$
It is straightforward to prove the following well-known observation, which enables us to express the deviation for a triangle as a sum of deviations over pairs of tangle segments in the triangle:
\begin{align}
\label{eq:dev_as_sum}
{\rm dev}_{\wh{t}_j} (K_j, J|_{\partial K_j}) = \sum_{1\le r<u\le l_j} {\rm dev}_{\wh{t}_j}(k_{j,r}\cup k_{j,u}, \, J|_{\partial k_{j,r}\cup \partial k_{j,u}}).
\end{align}

We establish a computational lemma for each summand of the right-hand-side of eq.\eqref{eq:dev_as_sum}.

\begin{figure}
\centering 
\begin{subfigure}[b]{0.24\textwidth}
\hspace*{2mm} 
\scalebox{0.8}{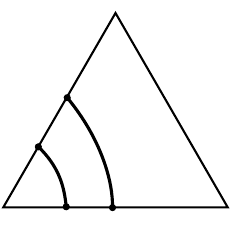}
\vspace{-3mm}
\caption*{Case (1)}
\label{fig:pair_of_loop_segments_case1}
\end{subfigure}
~
\begin{subfigure}[b]{0.24\textwidth}
\hspace*{2mm}  \scalebox{0.8}{
\begingroup%
  \makeatletter%
  \providecommand\color[2][]{%
    \errmessage{(Inkscape) Color is used for the text in Inkscape, but the package 'color.sty' is not loaded}%
    \renewcommand\color[2][]{}%
  }%
  \providecommand\transparent[1]{%
    \errmessage{(Inkscape) Transparency is used (non-zero) for the text in Inkscape, but the package 'transparent.sty' is not loaded}%
    \renewcommand\transparent[1]{}%
  }%
  \providecommand\rotatebox[2]{#2}%
  \newcommand*\fsize{\dimexpr\f@size pt\relax}%
  \newcommand*\lineheight[1]{\fontsize{\fsize}{#1\fsize}\selectfont}%
  \ifx\svgwidth\undefined%
    \setlength{\unitlength}{113.38582677bp}%
    \ifx\svgscale\undefined%
      \relax%
    \else%
      \setlength{\unitlength}{\unitlength * \real{\svgscale}}%
    \fi%
  \else%
    \setlength{\unitlength}{\svgwidth}%
  \fi%
  \global\let\svgwidth\undefined%
  \global\let\svgscale\undefined%
  \makeatother%
  \begin{picture}(1,1)%
    \lineheight{1}%
    \setlength\tabcolsep{0pt}%
    \put(0,0){\includegraphics[width=\unitlength,page=1]{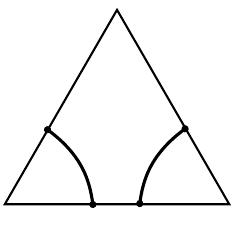}}%
    \put(0.32921321,0.34395084){\color[rgb]{0,0,0}\makebox(0,0)[lt]{\lineheight{1.25}\smash{\begin{tabular}[t]{l}$k_2$\end{tabular}}}}%
    \put(0.57070012,0.34069834){\color[rgb]{0,0,0}\makebox(0,0)[lt]{\lineheight{1.25}\smash{\begin{tabular}[t]{l}$k_1$\end{tabular}}}}%
    \put(0.56797829,0.06288586){\color[rgb]{0,0,0}\makebox(0,0)[lt]{\lineheight{1.25}\smash{\begin{tabular}[t]{l}$x_2$\end{tabular}}}}%
    \put(0.79948826,0.47960458){\color[rgb]{0,0,0}\makebox(0,0)[lt]{\lineheight{1.25}\smash{\begin{tabular}[t]{l}$x_1$\end{tabular}}}}%
    \put(0.33646817,0.06288576){\color[rgb]{0,0,0}\makebox(0,0)[lt]{\lineheight{1.25}\smash{\begin{tabular}[t]{l}$x_3$\end{tabular}}}}%
    \put(0.1049583,0.47960463){\color[rgb]{0,0,0}\makebox(0,0)[lt]{\lineheight{1.25}\smash{\begin{tabular}[t]{l}$x_4$\end{tabular}}}}%
    \put(0,0){\includegraphics[width=\unitlength,page=2]{commutatorcase2_new2.pdf}}%
    \put(0.47023706,0.05906229){\color[rgb]{0,0,0}\makebox(0,0)[lt]{\lineheight{1.25}\smash{\begin{tabular}[t]{l}$\succ$\end{tabular}}}}%
  \end{picture}%
\endgroup%
}
\vspace{-3mm}
\caption*{Case (2)}
\label{fig:pair_of_loop_segments_case2}
\end{subfigure}
~
\begin{subfigure}[b]{0.24\textwidth}
\hspace*{2mm} \scalebox{0.8}{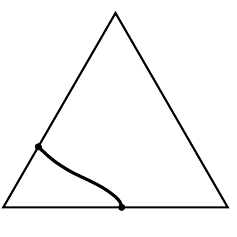}
\vspace{-3mm}
\caption*{Case (3)}
\label{fig:pair_of_loop_segments_case3}
\end{subfigure}
~
\begin{subfigure}[b]{0.24\textwidth}
\hspace*{2mm} \scalebox{0.8}{
\begingroup%
  \makeatletter%
  \providecommand\color[2][]{%
    \errmessage{(Inkscape) Color is used for the text in Inkscape, but the package 'color.sty' is not loaded}%
    \renewcommand\color[2][]{}%
  }%
  \providecommand\transparent[1]{%
    \errmessage{(Inkscape) Transparency is used (non-zero) for the text in Inkscape, but the package 'transparent.sty' is not loaded}%
    \renewcommand\transparent[1]{}%
  }%
  \providecommand\rotatebox[2]{#2}%
  \newcommand*\fsize{\dimexpr\f@size pt\relax}%
  \newcommand*\lineheight[1]{\fontsize{\fsize}{#1\fsize}\selectfont}%
  \ifx\svgwidth\undefined%
    \setlength{\unitlength}{113.38582677bp}%
    \ifx\svgscale\undefined%
      \relax%
    \else%
      \setlength{\unitlength}{\unitlength * \real{\svgscale}}%
    \fi%
  \else%
    \setlength{\unitlength}{\svgwidth}%
  \fi%
  \global\let\svgwidth\undefined%
  \global\let\svgscale\undefined%
  \makeatother%
  \begin{picture}(1,1)%
    \lineheight{1}%
    \setlength\tabcolsep{0pt}%
    \put(0,0){\includegraphics[width=\unitlength,page=1]{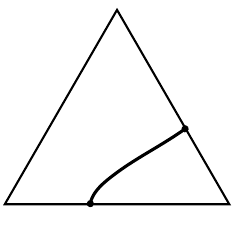}}%
    \put(0.31749642,0.38975286){\color[rgb]{0,0,0}\makebox(0,0)[lt]{\lineheight{1.25}\smash{\begin{tabular}[t]{l}$k_2$\end{tabular}}}}%
    \put(0.56004853,0.38863068){\color[rgb]{0,0,0}\makebox(0,0)[lt]{\lineheight{1.25}\smash{\begin{tabular}[t]{l}$k_1$\end{tabular}}}}%
    \put(0.58120747,0.06288586){\color[rgb]{0,0,0}\makebox(0,0)[lt]{\lineheight{1.25}\smash{\begin{tabular}[t]{l}$x_3$\end{tabular}}}}%
    \put(0.79948826,0.47960458){\color[rgb]{0,0,0}\makebox(0,0)[lt]{\lineheight{1.25}\smash{\begin{tabular}[t]{l}$x_1$\end{tabular}}}}%
    \put(0.33646817,0.06288576){\color[rgb]{0,0,0}\makebox(0,0)[lt]{\lineheight{1.25}\smash{\begin{tabular}[t]{l}$x_2$\end{tabular}}}}%
    \put(0.1049583,0.47960463){\color[rgb]{0,0,0}\makebox(0,0)[lt]{\lineheight{1.25}\smash{\begin{tabular}[t]{l}$x_4$\end{tabular}}}}%
    \put(0,0){\includegraphics[width=\unitlength,page=2]{commutatorcase4_new2.pdf}}%
    \put(0.45372891,0.06560011){\color[rgb]{0,0,0}\makebox(0,0)[lt]{\lineheight{1.25}\smash{\begin{tabular}[t]{l}$\prec$\end{tabular}}}}%
  \end{picture}%
\endgroup%
}
\vspace{-3mm}
\caption*{Case (4)}
\label{fig:pair_of_loop_segments_case4}
\end{subfigure}
\vspace{-3mm}
\caption{Pairs of tangle segments over a triangle (arrows in boundary are clockwise)}
\vspace{-4mm}
\label{fig:pair_of_loop_segments_case}
\end{figure}

\begin{lemma}[deviation for a pair of stated tangle segments in a triangle]
\label{lem:computation_of_commutator_of_log_quantum_traces}
If $(k_1\cup k_2, s_0)$ is a stated V-tangle in the thickening of a non-self-folded triangle $\wh{t}$ consisting of two tangle segments $k_1$ and $k_2$ as in the cases (1) or (2) of Fig.\ref{fig:pair_of_loop_segments_case}, and if we write $\varepsilon_1=s_0(x_1)$, $\varepsilon_2=s_0(x_2)$, $\varepsilon_3=s_0(x_3)$ and $\varepsilon_4 = s_0(x_4)$, then we have
\begin{align*}
{\rm dev}_{\wh{t}}(k_1 \cup k_2, s_0) = \left\{
\begin{array}{ll}
\varepsilon_1 \varepsilon_4 - \varepsilon_2 \varepsilon_3, & \mbox{in case (1)}, \\
\varepsilon_1 \varepsilon_3 - \varepsilon_1 \varepsilon_4 + \varepsilon_2 \varepsilon_4 & \mbox{in case (2)}.
\end{array}
\right. 
\end{align*}
\end{lemma}

The proof follows from straightforward computation. 

\subsection{Gabella triangle factors}

Now we turn to the triangle factors of the term for Gabella quantum holonomy, the right hand side of eq.\eqref{eq:term_equality}. Recall from \eqref{eq:Gabella_coefficient} that 
$\ol{\ul{\Omega}}(\til{K}^J;\omega)$ is the product of the biangle quantum holonomy factors for biangles and a power of $q$ to the (negative) writhe of the Gabella lift $\til{K}^J$ counted over triangles of $\wh{\Delta}$. One must be careful that this writhe is considered over the branched double cover surface $\wh{\Sigma}_\Delta$, and not over $\Sigma$. Writhe is a net sum of signs of all crossings, where each crossing is formed by two segments living over a same triangle of $\wh{\Delta}$. Thus we can write this writhe ${\rm wr}_{\wh{\Delta}}(\til{K}^J)$ as a double sum, where the outer sum is over triangles of $\wh{\Delta}$ and the inner sum is over pairs of segments living over each triangle. For triangle $\wh{t}_j$ of $\wh{\Delta}$, denote by $\til{K}^J_j$ the part of $\til{K}^J$ over $\wh{t}_j$; more precisely, the subset of $\til{K}^J$ of points whose images under the composition $\til{\Sigma}_\Delta \times(-1,1) \to \Sigma \times(-1,1) \to \Sigma$ lie in $\wh{t}_j$. Define
${\rm wr}_{\wh{t}_j}(\til{K}^J_j)$ to be the writhe of $\til{K}^J_j$, i.e. the net sum of signs of the crossings (see Def.\ref{def:Gabella_quantum_holonomy}(G2)) of $\til{K}^J_j$. Then
$$
{\rm wr}_{\wh{\Delta}}(\til{K}^J) = \sum_{j=1}^m {\rm wr}_{\wh{t}_j}(\til{K}^J_j)
$$
Note that $\til{K}^J_j$ is a Gabella lift of the stated oriented VH-tangle $K_j$ living in the thickening of the triangle $\wh{t}_j$. Denote the components of $\til{K}^J_j$ by $\til{k}_{j,1}$, \ldots, $\til{k}_{j,l_j}$. Then one can observe that these are Gabella lifts of the corresponding components of $K_j$ which are stated oriented tangle segments $(k_{j,1}, J|_{\partial k_{j,1}})$, \ldots, $(k_{j,l_j}, J|_{\partial k_{j,l_j}})$. That is, $\til{k}_{j,r}$ is a Gabella lift of the horizontal triangle segment $k_{j,r}$ with the state coinciding with $J|_{\partial k_{j,r}}$; see Def.\ref{def:horizontal_triangle_segment}. One observes
\begin{align}
\label{eq:wr_as_sum}
{\rm wr}_{\wh{t}_j}(\til{K}^J_j) = \sum_{1\le r<u\le l_j} {\rm wr}_{\wh{t}_j}(\til{k}_{j,r} \cup \til{k}_{j,u}), 
\end{align}
thus we indeed expressed ${\rm wr}_{\wh{\Delta}}(\til{K}^J)$ as a double sum. Now we establish a computational lemma for the innermost summand.
\begin{lemma}[writhe for a pair of stated oriented tangle segments in a triangle]
\label{lem:computation_of_writhe_for_a_pair_of_stated_loop_segments}
Let $\wh{t}$ be a non-self-folded triangle, viewed as a generalized marked surface $(\Sigma',\mathcal{P}')$, with the unique ideal triangulation $\Delta'$. Let $\frak{S}' = \Sigma'\setminus\mathcal{P}'$. Let $(k_1\cup k_2,s_0)$ be a stated oriented VH-tangle in $\frak{S}'\times(-1,1)$ as in the case (1) or (2) of Fig.\ref{fig:pair_of_loop_segments_case}, with {\em arbitrary} orientations on $k_1$ and $k_2$. Let $\til{k}_1 \cup \til{k}_2$ be a Gabella lift in $\til{\Sigma}'_{\Delta'} \times (-1,1)$ of $(k_1\cup k_2,s_0)$ (in the sense of Def.\ref{def:Gabella_lift}), where each of $\til{k}_1$ and $\til{k}_2$ is a Gabella lift of $(k_1,s_0|_{\partial k_1})$ and $(k_2,s_0|_{\partial k_2})$ respectively (in the sense of Def.\ref{def:horizontal_triangle_segment} with the prescribed states).

\vs

If we write $\varepsilon_1=s_0(x_1)$, $\varepsilon_2=s_0(x_2)$, $\varepsilon_3=s_0(x_3)$ and $\varepsilon_4 = s_0(x_4)$, then the writhe of the Gabella lift $\til{k}_1\cup \til{k}_2$ is given by: 
$$
{\rm wr}_{\wh{t}} \, (\til{k}_1 \cup \til{k}_2) = \left\{ 
{\renewcommand{\arraystretch}{1.2} \begin{array}{ll}
- \frac{1}{4}(\varepsilon_1 - \varepsilon_2)(\varepsilon_3+\varepsilon_4), & \mbox{in case (1)}, \\
- \frac{1}{4}(\varepsilon_1 - \varepsilon_2)(\varepsilon_3-\varepsilon_4), & \mbox{in case (2)},
\end{array}}
\right.
$$
where each sign $\pm$ is understood as the number $\pm 1$. 
\end{lemma}

\begin{remark}
This is another occurrence of the phenomenon of orientation insensitivity.
\end{remark}

\begin{figure}
\centering 
\hspace{5mm} \scalebox{0.8}{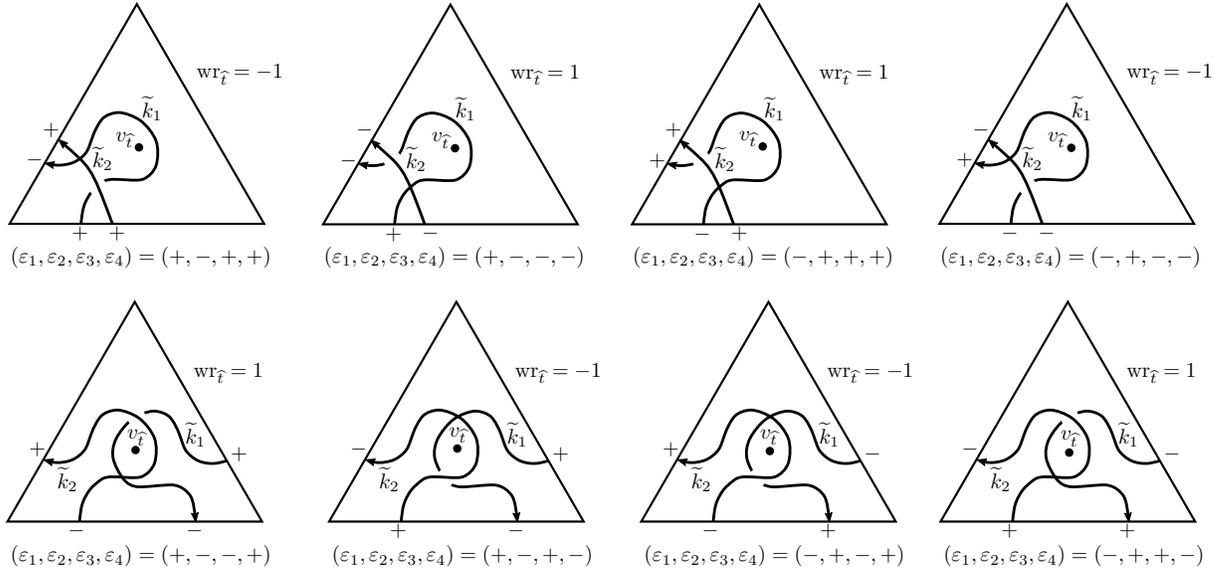}
\caption{Diagrams for a Gabella lift $\til{k}_1 \cup \til{k}_2$ in $\til{\Sigma}_\Delta$ of a union of two stated oriented VH-tangles in a triangle $\wh{t}$, for computation of writhe}
\vspace{-3mm}
\label{fig:computation_of_writhe}
\end{figure}

{\it Proof.} Note that the Gabella lift $\til{k}_1 \cup \til{k}_2$ is uniquely determined up to VH-isotopy  in this situation; any Gabella lift gives the same answer for the writhe, because writhe is well-defined up to (VH-)isotopy. Consider the case (1) of Fig.\ref{fig:pair_of_loop_segments_case}, where both segments $k_1$ and $k_2$ live in the same corner; with respect to the vertex of this corner, $k_1$ is located ``inner" and $k_2$ ``outer". If $\varepsilon_1=\varepsilon_2$ then the Gabella lift $\til{k}_1$ of the inner segment bounds a corner region of the triangle not containing the branch point (see Fig.\ref{fig:turns}), hence can be isotoped so that it doesn't have a crossing with $\til{k}_2$. If $\varepsilon_3\neq \varepsilon_4$ the Gabella lift $\til{k}_2$ of the outer segment bounds a corner region of the triangle containing the branch point (see Fig.\ref{fig:turns}), hence can be isotoped so that it doesn't have a crossing with $\til{k}_1$. So, in these cases, we have ${\rm wr}_{\wh{t}}\,(\til{k}_1\cup \til{k}_2)=0$. Meanwhile, note that $-\frac{1}{4}(\varepsilon_1 - \varepsilon_2)(\varepsilon_3+\varepsilon_4)$ equals zero in these cases, so we get the sought-for equality. Now assume that both $\varepsilon_1\neq \varepsilon_2$ and $\varepsilon_3 = \varepsilon_4$ hold; there are four such possibilities for signs. The upper half of Fig.\ref{fig:computation_of_writhe} presents the diagram  of $\til{k}_1\cup \til{k}_2$ for each of these four possibilities under a particular choice of orientations of $k_1$ and $k_2$; one can easily verify the equality in these cases. For each case in the upper half of Fig.\ref{fig:computation_of_writhe}, note that the projections of $\til{k}_1$ and $\til{k}_2$ in $\wh{t}$ meet at two points in $\Sigma$, only one of which is a crossing point over the branched double cover $\til{\Sigma}_\Delta$ of $\til{k}_1\cup \til{k}_2$. If one changes the orientation of either one of $k_1$ and $k_2$, say $k_i$, then the projections of $\til{k}_1$ and $\til{k}_2$ in $\wh{t}$ are same as before, with the orientation of $\til{k}_i$ reversed. Then the diagram of each of $\til{k}_1$ and $\til{k}_2$ stays the same (though the diagram of $\til{k}_1 \cup \til{k}_2$ is changed), with the orientation of $\til{k}_i$ reversed, and the sheet numbers of points of $\til{k}_i$ changed from before. So, out of the two intersections of the projections of $\til{k}_1$ and $\til{k}_2$ in $\wh{t}$, the one that used to be the crossing of $\til{k}_1\cup \til{k}_2$ is not a crossing anymore, and the remaining one now becomes the crossing point in $\til{\Sigma}_\Delta$. However, one easily verifies that the sign of the crossing is same as before. Similar argument holds whenever one changes the orientation of both $k_1$ and $k_2$.

\vs

Now consider the case (2) of Fig.\ref{fig:pair_of_loop_segments_case}. If $\varepsilon_1=\varepsilon_2$ or $\varepsilon_3=\varepsilon_4$ holds, it is easy to see that both sides of the equation in the statement are zero. When both $\varepsilon_1\neq \varepsilon_2$ and $\varepsilon_3\neq \varepsilon_4$ hold, one can check the equality case by case, for each of the four possibilities for the signs; see the lower half of Fig.\ref{fig:computation_of_writhe}, drawn for particular choices of orientations. As in the previous case, one easily observes that the equality also holds for the other choices of orientations. \qed

\vs

Using the computational lemmas developed so far, we now arrive at the relationship between the triangle factors of the Bonahon-Wong term and the Gabella term. In fact, they are not exactly equal, but equal only up to powers of $\omega$ to the signed order correction amounts defined in Def.\ref{def:signed_order_correction_amount}, which is due to the fact that the Bonahon-Wong quantum trace is defined for V-tangles up to V-isotopy, but writhe (and hence Gabella quantum holonomy) is defined for VH-tangles up to VH-isotopy. The precise statement is as follows, and is the main technical lemma of the proof of our main theorem.

\begin{lemma}[equality of the triangle factors]
\label{lem:equality_of_the_triangle_factors}
One has
$$
{\rm dev}_{\wh{t}_j} (K_j, J|_{\partial K_j}) = - 4 {\rm wr}_{\wh{t}_j}(\til{K}^J_j) + 2{\rm wr}(K_j) + \partial \mathcal{C}_{\wh{t}_j}(K_j, J|_{\partial K_j}) 
$$
\end{lemma}

{\it Proof.} We first rewrite the boundary signed order correction amount as the sum over pairs of tangle segments $k_{j,1}$, \ldots, $k_{j,l_j}$ of $K_j$:
\begin{align}
\label{eq:boundary_signed_order_correction_amount_as_sum_over_pairs}
\partial \mathcal{C}_{\wh{t}_j} (K_j, J|_{\partial K_j}) = \sum_{1\le r<u\le l_j} \partial \mathcal{C}_{\wh{t}_j}( k_{j,r}\cup k_{j,u}, \, J|_{\partial k_{j,r}\cup\partial k_{j,u}}).
\end{align}
This is easily proved by looking at the contributions at each boundary arc of $\wh{t}_j$.

\vs

Let $1\le r < u \le l_j$. In view of the sum expression for the deviation obtained in eq.\eqref{eq:dev_as_sum}, and of
$$
{\rm wr}(K_j) = \sum_{1\le r<u\le l_j} {\rm wr}(k_{j,r} \cup k_{j,u}),
$$
we now investigate the number
\begin{align*}
(*)_{j;r,u} := {\rm dev}_{\wh{t}_j}(k_{j,r} \cup k_{j,u}, J|_{\partial k_{j,r} \cup \partial k_{j,u}}) + 4 {\rm wr}_{\wh{t}_j}(\til{k}_{j,r} \cup \til{k}_{j,u}) - 2{\rm wr}(k_{j,r} \cup k_{j,u}), 
\end{align*}
and will show that it equals the $(r,u)$-th summand in the right hand side of eq.\eqref{eq:boundary_signed_order_correction_amount_as_sum_over_pairs}. We apply Lem.\ref{lem:computation_of_commutator_of_log_quantum_traces} to $(k_{j,r} \cup k_{j,u}, J_{\partial k_{j,r}\cup \partial k_{j,u}})$, for which we also apply Lem.\ref{lem:computation_of_writhe_for_a_pair_of_stated_loop_segments} with the Gabella lift $\til{k}_{j,r} \cup \til{k}_{j,u}$. 

\vs

If $k_{j,r}$ and $k_{j,u}$ are like $k_1$ and $k_2$ of the case (1) of Fig.\ref{fig:pair_of_loop_segments_case}, we have
\begin{align*}
(*)_{j;r,u} & = (\underbrace{J(x_1) J(x_4) - J(x_2) J(x_3)}_{\mbox{\tiny Lem.\ref{lem:computation_of_commutator_of_log_quantum_traces}}})  - \underbrace{(J(x_1) - J(x_2))(J(x_3) + J(x_4))}_{\mbox{\tiny Lem.\ref{lem:computation_of_writhe_for_a_pair_of_stated_loop_segments}}}  - 2 \underbrace{ {\rm wr}(k_{j,r}\cup k_{j,u}) }_{=0}  \\
& = ( - J(x_1) J(x_3) + J(x_2) J(x_4) )  \\
& = {\rm sgn}(e; \overrightarrow{x_1x_3}) J(x_1) J(x_3) + {\rm sgn}(e';\overrightarrow{x_2x_4}) J(x_2) J(x_4),
\end{align*}
where $e$ is the side of $\wh{t}_j$ containing $x_1$ and $x_3$, while $e'$ is the side containing $x_2$ and $x_4$; note $x_1<x_3$ and $x_2<x_4$ in the ordering coming from elevations, because $k_{j,r}$ has lower elevation than $k_{j,u}$. Therefore
\begin{align}
\label{eq:star_j_r_u_equality}
(*)_{j,r;u} = \partial \mathcal{C}_{\wh{t}_j}(k_{j,r}\cup k_{j,u}, \, J_{\partial k_{j,r}\cup\partial k_{j,u}})
\end{align}
holds in this case, as desired.

\vs

Suppose now that $k_{j,r}$ and $k_{j,u}$ are like $k_1$ and $k_2$ of the case (3) of Fig.\ref{fig:pair_of_loop_segments_case}. As tangles, this case can be obtained from the case (1) by a V-isotopy that slides the endpoints $x_1$ of $k_1$ and $x_3$ of $k_2$. After this sliding (or boundary exchange move (M2) in \S\ref{label:tangles_and_skein_algebra}), note that ${\rm dev}_{\wh{t}_j}(k_{j,r} \cup k_{j,u}, J|_{\partial k_{j,r} \cup \partial k_{j,u}})$ does not change, and the boundary signed order correction amount $\partial \mathcal{C}_{\wh{t}_j}(k_{j,r}\cup k_{j,u}, \, J|_{\partial k_{j,r}\cup\partial k_{j,u}})$ increases by two if $J(x_1)=J(x_3)$ and decreases by two if $J(x_1)\neq J(x_3)$. Suppose that the orientations on $k_1$ and $k_2$ are parallel, e.g. going from $x_1$ to $x_2$ and $x_3$ to $x_4$. After sliding, $-2{\rm wr}(k_{j,r}\cup k_{j,u})$ decreases by two, while $4 {\rm wr}_{\wh{t}_j}(\til{k}_{j,r} \cup \til{k}_{j,u})$ increases by four if $J(x_1)=J(x_3)$ (because the parts of the two segments near the boundary edge are living in the same sheet) and stays the same if $J(x_1) \neq J(x_3)$ (because the parts of the two segments near the boundary edge are living in different sheets). Thus eq.\eqref{eq:star_j_r_u_equality} holds. Now suppose that the orientations of $k_1$ and $k_2$ are not parallel; e.g. $k_2$ is going from $x_4$ to $x_3$. After the sliding move from the case (1), $-2{\rm wr}(k_{j,r}\cup k_{j,u})$ increases by two, while $4 {\rm wr}_{\wh{t}_j}(\til{k}_{j,r} \cup \til{k}_{j,u})$ stays the same if $J(x_1) = J(x_3)$ and decreases by four if $J(x_1)\neq J(x_3)$. Thus eq.\eqref{eq:star_j_r_u_equality} still holds.

\vs

If $k_{j,r}$ and $k_{j,u}$ are like $k_1$ and $k_2$ of the case (2) of Fig.\ref{fig:pair_of_loop_segments_case}, we have
\begin{align*}
\hspace{-1,5mm} (*)_{j;r,u} & = (\underbrace{J(x_1) J(x_3) - J(x_1) J(x_4) + J(x_2)J(x_4)}_{\mbox{\tiny Lem.\ref{lem:computation_of_commutator_of_log_quantum_traces}}}) 
- \underbrace{(J(x_1) - J(x_2))(J(x_3) - J(x_4))}_{\mbox{\tiny Lem.\ref{lem:computation_of_writhe_for_a_pair_of_stated_loop_segments}}} ) - 2 \underbrace{ {\rm wr}(k_{j,r}\cup k_{j,u}) }_{=0}  \\
& = J(x_2) J(x_3) 
= {\rm sgn}(e;\overrightarrow{x_2x_3}) J(x_2) J(x_3),
\end{align*}
where $e$ is the side of $\wh{t}_j$ containing $x_2$ and $x_3$; note $x_2<x_3$ in the ordering coming from elevations, because $k_{j,r}$ has lower elevation than $k_{j,u}$. Thus eq.\eqref{eq:star_j_r_u_equality} holds also in this case. 

\vs

Suppose now that $k_{j,r}$ and $k_{j,u}$ are like $k_1$ and $k_2$ of the case (4) of Fig.\ref{fig:pair_of_loop_segments_case}. As tangles, this case can be obtained from the case (2) by a V-isotopy that slides the endpoints $x_2$ of $k_1$ and $x_3$ of $k_2$. After this sliding (or boundary exchange move), note that ${\rm dev}_{\wh{t}_j}(k_{j,r} \cup k_{j,u}, J|_{\partial k_{j,r} \cup \partial k_{j,u}})$ does not change, and the boundary signed order correction amount $\partial \mathcal{C}_{\wh{t}_j}(k_{j,r}\cup k_{j,u}, \, J|_{\partial k_{j,r}\cup\partial k_{j,u}})$ decreases by two if $J(x_2)=J(x_3)$ and increases by two if $J(x_2)\neq J(x_3)$. Suppose that the orientations on $k_1$ and $k_2$ are parallel, e.g. going from $x_2$ to $x_1$ and $x_3$ to $x_4$. After sliding, $-2{\rm wr}(k_{j,r}\cup k_{j,u})$ increases by two, while $4 {\rm wr}_{\wh{t}_j}(\til{k}_{j,r} \cup \til{k}_{j,u})$ decreases by four if $J(x_2)=J(x_3)$ and stays the same if $J(x_2) \neq J(x_3)$. Thus eq.\eqref{eq:star_j_r_u_equality} holds. Now suppose that the orientations of $k_1$ and $k_2$ are not parallel, e.g. $k_1$ is going from $x_1$ to $x_2$. After the sliding move from the case (2), $-2{\rm wr}(k_{j,r}\cup k_{j,u})$ decreases by two, while $4 {\rm wr}_{\wh{t}_j}(\til{k}_{j,r} \cup \til{k}_{j,u})$ stays the same if $J(x_2) = J(x_3)$ and increases by four if $J(x_2)\neq J(x_3)$. Thus eq.\eqref{eq:star_j_r_u_equality} still holds. 

\vs

When $k_{j,r}$ and $k_{j,u}$ play the roles of $k_2$ and $k_1$ of Fig.\ref{fig:pair_of_loop_segments_case} in each of the cases we dealt with so far, eq.\eqref{eq:star_j_r_u_equality} still holds because of the skew-symmetry of each of the three terms of $(*)_{j;r,u}$, as well as that of the boundary signed order correction amount $\partial \mathcal{C}_{\wh{t}_j}(k_{j,r}\cup k_{j,u}, \, J|_{\partial k_{j,r}\cup\partial k_{j,u}})$.

\vs

Now, taking the sum of the equality eq.\eqref{eq:star_j_r_u_equality} over all pairs $(r,u)$ with $1\le r<u\le l_j$, we get the desired result, in view of eq.\eqref{eq:dev_as_sum}, eq.\eqref{eq:wr_as_sum}, and eq.\eqref{eq:boundary_signed_order_correction_amount_as_sum_over_pairs}. \qed

\vs

The boundary signed order correction amount will also appear in the relationship between the biangle factors of the Bonahon-Wong term and the Gabella term, canceling the ones appearing for the triangles.

\subsection{Bonahon-Wong biangle factors via the $\mathcal{U}_q(\frak{sl}_2)$ Reshetikhin-Turaev operator invariant}
\label{subsec:BW_biangle_factors_as_RT} Now we move on to the biangle factors. Our approach is via the operator invariants of tangles in a biangle, as it was in \S\ref{subsec:biangle_quantum_holonomy} when investigating Gabella's biangle quantum holonomy. We follow notations of \S\ref{subsec:biangle_quantum_holonomy}, some of which we recall now. Let $\vec{B}=(B,{\rm dir})$ be a directed biangle (see Def.\ref{def:directed_biangle}) whose boundary arcs are denoted by $b_{\rm in}$ and $b_{\rm out}$, and let $\mathcal{D}(\vec{B})$ be the set of all equivalence classes $[D]$ of boundary-ordered oriented tangle diagrams (see Def.\ref{def:D_B}). Let the vector space $V$ be given as in eq.\eqref{eq:V}, with the ordered basis $\{\xi_+,\xi_-\}$. We claimed in Prop.\ref{prop:G_D} that there exists a unique assignment $G$ that assigns to each $[D] \in \mathcal{D}(\vec{B})$ a linear operator $G([D]) : V^{\otimes|\partial_{\rm in} D|} \to V^{\otimes |\partial_{\rm out}D|}$ satisfying certain properties, where $\partial_{\rm in} D = \partial_{b_{\rm in}}D$ and $\partial_{\rm out} D = \partial_{b_{\rm out}}D$. Then the Gabella biangle quantum holonomy ${\rm TrHol}^\omega_{\vec{B}}([k,s_0])$ for a stated oriented VH-tangle $(k,s_0)$ in the thickening of $\vec{B}$ was defined as a matrix element of the operator $G([D])$ (see Def.\ref{def:Gabella_biangle_quantum_holonomy_as_matrix_elements}), where $D$ is the boundary-ordered oriented tangle diagram of $k$. We still have to prove Prop.\ref{prop:G_D}, and also to find a relationship between the Gabella biangle quantum holonomy and the Bonahon-Wong biangle quantum trace dealt with in \S\ref{subsec:biangles}. We accomplish these two tasks simultaneously.

\vs

The main ingredient of the present section is the Reshetikhin-Turaev operator invariant \cite{RT, Turaev16} for oriented tangles in a thickened biangle, associated to the irreducible $2$-dimensional representation of the quantum group $\mathcal{U}_q(\frak{sl}_2)$. While suggesting the readers consult the excellent monograph by Ohtsuki \cite{O} for details on the theory of operator invariants, we just state and use the known result here.
\begin{proposition}[Reshetikhin-Turaev operator invariant for 2-dimensional irreducible representation of $\mathcal{U}_q(\frak{sl}_2)$ {\cite{O, RT, Turaev16}}]
\label{prop:F_D}
There exists a unique map $F$ assigning to each $[D] \in \mathcal{D}(\vec{B})$ a linear operator
$$
F([D]) : V^{\otimes|\partial_{\rm in} D|} \to V^{\otimes |\partial_{\rm out}D|}
$$
that respects the composition and tensor product (see Def.\ref{def:tensor_product_and_composition_of_oriented_tangle_diagrams_in_a_directed_biangle}), i.e.
$$
F([D_1]\circ [D_2]) = F([D_1])\circ F([D_2]) \quad\mbox{and}\quad
F([D_1]\otimes [D_2]) = F([D_1]) \otimes F([D_2]),
$$
and whose values at elementary classes $[D] \in \mathcal{D}(\vec{B})$, defined in Def.\ref{def:elementary_tangle_diagrams_in_a_directed_biangle} using Fig.\ref{fig:elementary_bo_oriented_tangle_diagrams_in_a_directed_biangle}, are:
\begin{enumerate}
\item[\rm (RT1)] (identity) if $D$ is as in Fig.\ref{fig:elementary_bo_oriented_tangle_diagrams_in_a_directed_biangle}(I) or (VI), then $F([D]) = {\rm id} : V \to V$;

\item[\rm (RT2)] (cup) if $D$ is as in Fig.\ref{fig:elementary_bo_oriented_tangle_diagrams_in_a_directed_biangle}(II) or (VII), then $F([D]) : \mathbb{C} \to V\otimes V$ sends $1 \in \mathbb{C}$ to $\omega\xi_+ \otimes \xi_- - \omega^5 \xi_- \otimes \xi_+$; 

\item[\rm (RT3)] (cap) if $D$ is as in Fig.\ref{fig:elementary_bo_oriented_tangle_diagrams_in_a_directed_biangle}(III) or (VIII), then $F([D]):V\otimes V \to \mathbb{C}$ sends $\xi_+ \otimes \xi_-$ to $-\omega^{-5}$ and $\xi_- \otimes \xi_+$ to $\omega^{-1}$, while sending other basis vectors to zero;

\item[\rm (RT4)] (height exchange) if $D$ is as in Fig.\ref{fig:elementary_bo_oriented_tangle_diagrams_in_a_directed_biangle}(IV), then the map $F([D]) :V \otimes V \to V\otimes V$ is given on the basis vectors as:
\begin{align}
\label{eq:RT4}
F([D]) : \xi_i \otimes \xi_j \mapsto \left\{
\begin{array}{ll}
\omega^2 \, \xi_i \otimes \xi_j, & \mbox{if $i= j$,} \\
\omega^{-2} ( \xi_+ \otimes \xi_- + (\omega^4-\omega^{-4})\xi_-\otimes \xi_+), & \mbox{if $i=+$ and $j=-$,} \\	
\omega^{-2} \xi_- \otimes \xi_+, & \mbox{if $i=-$ and $j=+$,}
\end{array}
\right.
\end{align}
while if $D$ is as in Fig.\ref{fig:elementary_bo_oriented_tangle_diagrams_in_a_directed_biangle}(IX), then the map $F([D]):V\otimes V\to V\otimes V$ is given by the inverse of the map in eq.\eqref{eq:RT4};

\item[\rm (RT5)] (positive crossing) if $D$ is as in Fig.\ref{fig:elementary_bo_oriented_tangle_diagrams_in_a_directed_biangle}(V), then $F([D]) : V\otimes V \to V\otimes V$ is given 
by the composition $(\mbox{inverse of the map in eq.\eqref{eq:RT4}}) \circ P$, where $P:V\otimes V \to V\otimes V$ is the position exchange map
\begin{align}
\label{eq:P}
P(\xi_i\otimes \xi_j) = \xi_j\otimes \xi_i, \qquad \forall i,j\in\{+,-\}.
\end{align}

\item[\rm (RT6)] (negative crossing) if $D$ is as in Fig.\ref{fig:elementary_bo_oriented_tangle_diagrams_in_a_directed_biangle}(X), then $F([D]) : V\otimes V \to V\otimes V$ is given by the inverse of the map $V\otimes V\to V\otimes V$ for case (RT5).
\end{enumerate}
\end{proposition}

\begin{lemma}[orientation insensitivity {\cite[\S3]{O}}]
$F([D])$ does not depend on the orientation of components of $D$.
\end{lemma}

To make the situation clear, we emphasize again that the validity of the above proposition and lemma is well known, whose proofs can be found in the original paper \cite{RT} or the books \cite{O, Turaev16}. One remark is that in these original references the Reshetikhin-Turaev operator invariants are formulated for H-tangles, and the condition (RT4) in the above Prop.\ref{prop:F_D} is to adapt such constructions to the setting of VH-tangles.

\vs

Our strategy to prove Prop.\ref{prop:G_D} is to find a relationship between $G([D])$ and $F([D])$, then to use Prop.\ref{prop:F_D}. However, a proof of Prop.\ref{prop:G_D} for $G([D])$ is not the only reason we recalled the Reshetikhin-Turaev operator invariant $F([D])$. In fact, we observe that the matrix elements of $F([D])$ are precisely the Bonahon-Wong biangle quantum trace. Note that this fact was not mentioned in the original work \cite{BW}; although it was already observed in \cite{CL}, we present a proof here, for completeness.
\begin{proposition}[Bonahon-Wong's biangle quantum trace as matrix elements of $F$ {\cite[Thm.5.2]{CL}}]
\label{prop:Bonahon-Wong_biangle_quantum_trace_as_matrix_elements}
Let $\vec{B}$ be a directed biangle, with the underlying biangle being $B =(\Sigma',\mathcal{P}')$, with $\frak{S}' = \Sigma'\setminus\mathcal{P}'$. Let $\omega \in \mathbb{C}^*$. Let $[k,s_0]$ be the VH-isotopy class of a stated VH-tangle in $\frak{S}'\times(-1,1)$. 

\vs

Denote by $[D,s_0] \in \mathcal{D}_{\rm s}(\vec{B})$ the equivalence class of stated boundary-ordered tangle diagrams for $[k,s_0]$ (with slight abuse of notation for $s_0$). Define the inward and outward basis vectors $\xi^{s_0}_{\rm in} \in V^{\otimes|\partial_{\rm in} D|}$ and $\xi^{s_0}_{\rm out} \in V^{\otimes|\partial_{\rm out}D|}$ as in Def.\ref{def:Gabella_biangle_quantum_holonomy_as_matrix_elements}, and let $F([D,s_0]) \in \mathbb{C}$ be the matrix element of the linear map $F([D]) : V^{\otimes|\partial_{\rm in} D|} \to V^{\otimes|\partial_{\rm out}D|}$ of Prop.\ref{prop:F_D} for these basis vectors, i.e.
$$
F([D,s_0]) := \langle \xi^{s_0}_{\rm out}, \, F([D]) \xi^{s_0}_{\rm in}\rangle
$$
where $\langle \cdot, \cdot \rangle$ is as defined in Def.\ref{def:basis_of_tensor_space_and_bilinear_form}. Then $F([D,s_0])$ coincides with the Bonahon-Wong biangle quantum trace ${\rm Tr}^\omega_B([k,s_0])$ of Prop.\ref{prop:BW_biangles}:
$$
{\rm Tr}^\omega_B ([k,s_0]) = F([D,s_0]).
$$
\end{proposition}

{\it Proof.} For each $[D] \in \mathcal{D}(\vec{B})$ we define $\wh{F}([D]) : V^{\otimes|\partial_{\rm in D}|} \to V^{\otimes|\partial_{\rm out}D|}$ as the unique linear map whose matrix element of each state $s_0$ of $D$ coincides with the number ${\rm Tr}^\omega_B ([k,s_0])$, where $(k,s_0)$ is a stated VH-tangle in $\frak{S}'\times(-1,1)$ whose projection to $\vec{B}$ yields a stated boundary-ordered tangle diagram equivalent to $(D,s_0)$ (by a slight abuse of notation for $s_0$):
$$
\langle \xi^{s_0}_{\rm out}, \, \wh{F}([D]) \xi^{s_0}_{\rm in} \rangle = {\rm Tr}^\omega_B ([k,s_0]).
$$
By Prop.\ref{prop:tangle_diagrams_and_tangles}(4), any such $(k,s_0)$ are VH-isotopic to each other. Since the Bonahon-Wong biangle quantum trace is defined on a V-isotopy class of stated tangles, we see that $\wh{F}([D])$ is well-defined. We will now show that $\wh{F}([D])$ respects the tensor product and composition operations in $\mathcal{D}(\vec{B})$, and that it has same values as $F([D])$ for elementary classes $[D]$.

\vs

First, let $[D_1],[D_2]\in\mathcal{D}(\vec{B})$ and consider the tensor product $[D_1]\otimes[D_2]$ defined as $[D_1' \cup D_2']$ with $D_1'$ and $D_2'$ as in Def.\ref{def:tensor_product_and_composition_of_oriented_tangle_diagrams_in_a_directed_biangle}; we shall show that $\wh{F}([D_1]\otimes[D_2]) = \wh{F}([D_1]) \otimes \wh{F}([D_2])$. Let $k'_1 \cup k'_2$ be a VH-tangle in $\frak{S}'\times(-1,1)$ whose projection gives a diagram equivalent to $D_1'\cup D_2'$, where $k'_1$ corresponds to $D_1'$ and $k'_2$ to $D_2'$. For any state $s_0$ of $D_1' \cup D_2'$, we can write it as $s_0 = s_1 \cup s_2$ for states $s_1$ and $s_2$ of $D_1'$ and $D_2'$, so that the inward basis vector (eq.\eqref{eq:inward_and_outward_basis_vectors}) for $s_0$ is $\xi^{s_0}_{\rm in} = \xi^{s_1}_{\rm in} \otimes \xi^{s_2}_{\rm in} \in V^{\otimes|\partial_{\rm in} D_1'|} \otimes V^{\otimes|\partial_{\rm in} D_2'|} \equiv V^{\otimes|\partial_{\rm in}(D_1' \cup D_2')|}$. By construction in Def.\ref{def:tensor_product_and_composition_of_oriented_tangle_diagrams_in_a_directed_biangle}, $D_1$ is vertically higher than $D_2$, hence $k_1'$ has higher elevation than $k_2'$. Hence the stated skein $[k_1'\cup k_2', s_0] = [k_1'\cup k_2',s_1\cup s_2] \in \mathcal{S}^A_{\rm s}(B)$ equals the product $[k_2',s_2] [k_1',s_1]$ of stated skeins. Using the fact that the Bonahon-Wong biangle quantum trace map ${\rm Tr}^\omega_B$ is an algebra homomorphism, we obtain
\begin{align*}
&  \langle \xi^{s_0}_{\rm out}, \, \wh{F}([D_1]\otimes[D_2])\xi^{s_0}_{\rm in} \rangle = \langle \xi^{s_0}_{\rm out}, \, \wh{F}([D_1'\cup D_2'])\xi^{s_0}_{\rm in} \rangle
= {\rm Tr}^\omega_B ([k_1'\cup k_2', \, s_0])  = {\rm Tr}^\omega_B ([k_2',  s_2][k_1',s_1])  \\
&  = {\rm Tr}^\omega_B ([k_2',s_2]) \, {\rm Tr}^\omega_B([k_1',s_1])
= \langle \xi^{s_2}_{\rm out}, \wh{F}([D_2]) \xi^{s_2}_{\rm in}\rangle \, \langle \xi^{s_1}_{\rm out}, \wh{F}([D_1]) \xi^{s_1}_{\rm in}\rangle  = \langle \xi^{s_1}_{\rm out}, \wh{F}([D_1]) \xi^{s_1}_{\rm in}\rangle \, \langle \xi^{s_2}_{\rm out}, \wh{F}([D_2]) \xi^{s_2}_{\rm in}\rangle \\
& = \langle \xi^{s_1}_{\rm out} \otimes \xi^{s_2}_{\rm out}, \, \wh{F}([D_1])\xi^{s_1}_{\rm in} \otimes \wh{F}([D_2])\xi^{s_2}_{\rm in} \rangle
 = \langle \xi^{s_0}_{\rm out}, \, (\wh{F}([D_1])\otimes\wh{F}([D_2]))(\xi^{s_0}_{\rm in}) \rangle
\end{align*}
hence $\wh{F}([D_1]\otimes[D_2]) = \wh{F}([D_1])\otimes\wh{F}([D_2])$ as desired.

\vs

Now let $[D_1], [D_2] \in \mathcal{D}(\vec{B})$, and suppose that $[D_1]$ is composable with $[D_2]$ as in Def.\ref{def:tensor_product_and_composition_of_oriented_tangle_diagrams_in_a_directed_biangle}; we shall show $\wh{F}([D_1] \circ [D_2]) = \wh{F}([D_1]) \circ \wh{F}([D_2])$. Write $[D] = [D_1] \circ [D_2]$ (as in Def.\ref{def:tensor_product_and_composition_of_oriented_tangle_diagrams_in_a_directed_biangle}), and let $k$ be a VH-tangle in $\frak{S}'\times(-1,1)$ whose projection gives a diagram equivalent to $D$, so that cutting $k$ (see Def.\ref{def:cutting}) along an internal arc of $\vec{B}$ yields VH-tangles $k_1$ and $k_2$ in biangles $\vec{B}_1$ and $\vec{B}_2$ giving diagrams equivalent to $D_1$ and $D_2$ respectively. Then we have a natural identifications $\partial_{\rm out} k = \partial_{\rm out}k_1$, $\partial_{\rm in} k_1=\partial_{\rm out}k_2$ and $\partial_{\rm in} k = \partial_{\rm in} k_2$. For each state $s_0$ of $k$, by the cutting property (see Prop.\ref{prop:BW_biangles}(1)) we have
$$
{\rm Tr}^\omega_B([k,s_0]) = \sum_{\mbox{\tiny compatible $s_1,s_2$}} {\rm Tr}^\omega_{B_1}([k_1,s_1]) \, {\rm Tr}^\omega_{B_2}([k_2,s_2]),
$$
where the sum is over all states $s_1$ and $s_2$ of $k_1$ and $k_2$ such that $s_1|_{\partial_{\rm out} k_1} = s_0|_{\partial_{\rm out} k}$, $s_1|_{\partial_{\rm in} k_1} = s_2|_{\partial_{\rm out} k_2}$ and $s_2|_{\partial_{\rm in} k_2} = s_0|_{\partial_{\rm in} k}$. Note
\begin{align*}
& \langle \xi^{s_0}_{\rm out}, \,\wh{F}([D_1] \circ [D_2]) \xi^{s_0}_{\rm in} \rangle = \langle \xi^{s_0}_{\rm out}, \, \wh{F}([D]) \xi^{s_0}_{\rm in} \rangle = {\rm Tr}^\omega_B([k,s_0]) \\
& = \sum_{\mbox{\tiny compatible $s_1,s_2$}} {\rm Tr}^\omega_{B_1}([k_1,s_1]) \, {\rm Tr}^\omega_{B_2}([k_2,s_2]) = \sum_{\mbox{\tiny compatible $s_1,s_2$}} \langle \xi^{s_1}_{\rm out}, \wh{F}([D_1])\xi^{s_1}_{\rm in}\rangle \, \langle \xi^{s_2}_{\rm out}, \wh{F}([D_2])\xi^{s_2}_{\rm in} \rangle.
\end{align*}
On the other hand, note
\begin{align*}
& \langle \xi^{s_0}_{\rm out}, (\wh{F}([D_1]) \circ \wh{F}([D_2])) \xi^{s_0}_{\rm in}\rangle 
= \langle \xi^{s_0}_{\rm out}, \, \wh{F}([D_1])( \wh{F}([D_2])\xi^{s_0}_{\rm in}) \rangle \\
& = \langle \xi^{s_0}_{\rm out}, \, \wh{F}([D_1])( {\textstyle \sum}_{\vec{\varepsilon}} \, \langle \xi^{\vec{\varepsilon}}, \wh{F}([D_2]) \xi^{s_0}_{\rm in} \rangle\, \xi^{\vec{\varepsilon}}) \rangle
= {\textstyle \sum}_{\vec{\varepsilon}} \, \langle \xi^{\vec{\varepsilon}}, \wh{F}([D_2]) \xi^{s_0}_{\rm in} \rangle \, \langle \xi^{s_0}_{\rm out}, \wh{F}([D_1]) \xi^{\vec{\varepsilon}} \rangle,
\end{align*}
where the sum is over all sign sequences $\vec{\varepsilon} \in \{+,-\}^{|\partial_{\rm out}D_2|}$. For each $\vec{\varepsilon}$, denote by $s_2$ the state of $k_2$ such that $s_2|_{\partial_{\rm in} k_2} = s_0|_{\partial_{\rm in} k}$ and $s_2|_{\partial_{\rm out} k_2} = \vec{\varepsilon}$ (using the horizontal ordering of $\partial_{\rm out} k_2$), and denote by $s_1$ the state of $k_1$ such that $s_1|_{\partial_{\rm in} k_1} = \vec{\varepsilon}$ (using the horizontal ordering of $\partial_{\rm in} k_1$) and $s_1|_{\partial_{\rm out} k_1} = s_0|_{\partial_{\rm out} k}$. Then this sum over $\vec{\varepsilon}$ can be easily seen to be equal to the above sum over compatible $s_1$ and $s_2$. So we obtain $\wh{F}([D_1]\circ [D_2]) = \wh{F}([D_1]) \circ \wh{F}([D_2])$ as desired.

\vs

It remains to check the values on the elementary classes. Prop.\ref{prop:F_D}(RT1) is easily seen to match Prop.\ref{prop:BW_biangles}(2)(I). For  Prop.\ref{prop:F_D}(RT2) the values of $\wh{F}([D])$ can be read from eq.\eqref{eq:BW_biangle_quantum_trace_III} of Lem.\ref{lem:biangle_quantum_trace_on_remaining_elementary_tangles}(1), with $x_1$ and $x_2$ playing the roles of $y_2$ and $y_1$, and match the values of $F([D])$ by inspection. For Prop.\ref{prop:F_D}(RT3) the values of $\wh{F}([D])$ can be read from Prop.\ref{prop:BW_biangles}(2)(II), with $y_1$ and $y_2$ playing the roles of $x_2$ and $x_1$, and match the values of $F([D])$. For Prop.\ref{prop:F_D}(RT4), the values of $F([D])$ match the values of $\wh{F}([D])$ as in eq.\eqref{eq:two_parallel_lines} of Lem.\ref{lem:biangle_quantum_trace_on_remaining_elementary_tangles}(2) by inspection. The remaining cases (RT5) and (RT6) can be similarly compared with the items (3) and (4) of Lem.\ref{lem:biangle_quantum_trace_on_remaining_elementary_tangles}, respectively. \qed

\vs

\begin{lemma}[charge conservation of $F$]
\label{lem:charge_conservation_of_F}
If $s$ is a state of a tangle diagram $D$ in $\vec{B}$ that is not charge-preserving (defined as in Lem.\ref{lem:simple_properties_of_G_D}), then the matrix element $\langle \xi^s_{\rm out}, \, F([D])\xi^s_{\rm in} \rangle$ is zero.
\end{lemma}
This lemma is well known in the operator invariant theory \cite{O}, and can also be directly proved in the style of the above proof of Prop.\ref{prop:Bonahon-Wong_biangle_quantum_trace_as_matrix_elements}, that is, investigate the behavior under the tensor product and composition in $\mathcal{D}(\vec{B})$, and then elementary classes $[D]$.

\subsection{Gabella biangle factors as twisting of Bonahon-Wong biangle factors}
\label{subsec:Gabella_biangle_factors_as_twisting_of_BW_biangle_factors}

We continue from the last subsection, to compare Bonahon and Wong's biangle quantum trace and Gabella's biangle quantum holonomy. We will see that they are not exactly equal, but differ by the signed order correction amount defined in Def.\ref{def:signed_order_correction_amount}; more precisely, we need an operator version of it.

\begin{definition}[signed order correction operator]
\label{def:signed_order_correction_operator}
Let $b$ be a boundary arc of a directed biangle $\vec{B}$, given the boundary-orientation; see Def.\ref{def:surface}. Let $Z$ be an ordered finite set in the interior of $b$. Define the \ul{\em signed order correction operator} of the ordered subset $Z$ on the boundary arc $b$ as the linear operator
$$
\mathcal{C}(b;Z) : V^{\otimes |Z|} \to V^{\otimes |Z|},
$$
given as follows. Let $Z = \{x_1,\ldots,x_{|Z|}\}$, written in the increasing order with respect to the horizontal ordering (see Def.\ref{def:horizontal_ordering}) on $Z$ coming from the orientation on $b$. For each state $s_0:Z\to \{+,-\}$, denote by $\xi^{s_0}$ the corresponding basis vector of $V^{\otimes|Z|}$
$$
\xi^{s_0} := \xi_{s_0(x_1)} \otimes \xi_{s_0(x_2)} \otimes \cdots \otimes \xi_{s_0(x_{|Z|})} ~\in~ V^{\otimes|Z|}.
$$
Then the map $\mathcal{C}(b;Z)$ is defined by its values on these basis vectors given as
$$
\mathcal{C}(b;Z) \xi^{s_0} := \omega^{\mathcal{C}(b;Z,s_0)} \, \xi^{s_0}
$$
for all states $s_0$ of $Z$, where $\mathcal{C}(b;Z,s_0)$ is as defined in Def.\ref{def:signed_order_correction_amount}.
\end{definition}
So $\mathcal{C}(b;Z)$ is diagonal and invertible, with the inverse map given by $\mathcal{C}(b;Z)^{-1} \xi^{s_0} = \omega^{-\mathcal{C}(b;Z,s_0)} \xi^{s_0}$. 

The following lemma is also one of the crucial technical tools of the present paper.
\begin{lemma}[relationship between the operator invariants {$G([D])$} and {$F([D])$}]
\label{lem:relationship_between_the_operator_invariants_G_and_F}
For any $[D] \in \mathcal{D}(\vec{B})$,
$$
F([D]) = \omega^{2{\rm wr}(D)} \, \mathcal{C}(b_{\rm out}; \partial_{\rm out} D) \circ G([D]) \circ \mathcal{C}(b_{\rm in}; \partial_{\rm in} D),
$$
where the orderings on $\partial_{\rm out} D$ and $\partial_{\rm in} D$ (see eq.\eqref{eq:partial_in_D_partial_out_D}) are the vertical orderings. 
\end{lemma}
Keep in mind that in Lem.\ref{lem:relationship_between_the_operator_invariants_G_and_F}, $b_{\rm out}$ and $b_{\rm in}$ are given the boundary-orientations (see Def.\ref{def:surface}), instead of the orientations described in Def.\ref{def:directed_biangle}.

\vs

{\it Proof.} We define
\begin{align*}
\wh{G}([D]):= \omega^{-2{\rm wr}(D)} \cdot \mathcal{C}(b_{\rm out}; \partial_{\rm out} D)^{-1} \circ F([D]) \circ \mathcal{C}(b_{\rm in}; \partial_{\rm in} D)^{-1}.
\end{align*}
The strategy is to show that $\wh{G}([D])$ respects the composition and tensor product operations of $\mathcal{D}(\vec{B})$, and has same values as $G([D])$ on the elementary classes of $\mathcal{D}(\vec{B})$. First, let $[D_1],[D_2]\in\mathcal{D}(\vec{B})$ and consider the tensor product $[D_1]\otimes[D_2]$ defined as $[D_1' \cup D_2']$ with $D_1'$ and $D_2'$ as in Def.\ref{def:tensor_product_and_composition_of_oriented_tangle_diagrams_in_a_directed_biangle}; we shall show that
$$
\wh{G}([D_1]\otimes[D_2]) = \wh{G}([D_1]) \otimes \wh{G}([D_2]).
$$
Since $D_1'$ and $D_2'$ are disjoint, we have ${\rm wr}(D_1'\cup D_2')={\rm wr}(D_1')+{\rm wr}(D_2')$. For any state $s_0$ of $D_1' \cup D_2'$, we can write it as $s_0 = s_1 \cup s_2$ for states $s_1$ and $s_2$ of $D_1'$ and $D_2'$, so that the inward basis vector (eq.\eqref{eq:inward_and_outward_basis_vectors}) for $s_0$ is given by the tensor product $\xi^{s_0}_{\rm in} = \xi^{s_1}_{\rm in} \otimes \xi^{s_2}_{\rm in} \in V^{\otimes|\partial_{\rm in} D_1'|} \otimes V^{\otimes|\partial_{\rm in} D_2'|} \equiv V^{\otimes|\partial_{\rm in}(D_1' \cup D_2')|}$. Note
$$
\mathcal{C}(b_{\rm in};  \partial_{\rm in} (D_1'\cup D_2'))^{-1} (\xi^{s_0}_{\rm in}) = \omega^{-\mathcal{C}(b_{\rm in} ; (\partial_{\rm in} D_1')\cup(\partial_{\rm in} D_2'), s_{1;{\rm in}} \cup s_{2;{\rm in}})} \, \xi^{s_0}_{\rm in},
$$
where $s_{1;{\rm in}} = s_1|_{\partial_{\rm in} D_1'}$ and $s_{2;{\rm in}} = s_2|_{\partial_{\rm in} D_2'}$. Observe that
\begin{align*}
\mathcal{C}(b_{\rm in}; (\partial_{\rm in} D_1') \cup (\partial_{\rm in} D_2'), s_{1;{\rm in}} \cup s_{2;{\rm in}}) & = \mathcal{C}(b_{\rm in} ; \partial_{\rm in} D_1', s_{1;{\rm in}}) + \mathcal{C}(b_{\rm in}; \partial_{\rm in} D_1', s_{2;{\rm in}}) \\
& \quad + \underbrace{ \sum_{x\in \partial_{\rm in} D_2', \, y\in \partial_{\rm in} D_1'} {\rm sgn}(b_{\rm in}; \overrightarrow{xy}) s_0(x) s_0(y)}.
\end{align*}
Because $\partial_{\rm in} D_2'$ is vertically and horizontally lower than $\partial_{\rm in} D_1'$, and the boundary-orientation on $b_{\rm in}$ is going toward the horizontally decreasing direction, we have ${\rm sgn}(b_{\rm in}; \overrightarrow{xy})=-1$ for $x\in \partial_{\rm in} D_2'$ and $y\in \partial_{\rm in} D_1'$, and hence the underbraced part equals $- (\sum_{x\in \partial_{\rm in} D_2'} s_0(x)) (\sum_{y\in \partial_{\rm in}D_1'} s_0(y))$. Likewise,
$$
\mathcal{C}(b_{\rm out}; \partial_{\rm out} (D_1'\cup D_2'))^{-1} (\xi^{s_0}_{\rm out}) = \omega^{-\mathcal{C}(b_{\rm out} ; (\partial_{\rm out} D_1')\cup(\partial_{\rm out} D_2'), s_{1;{\rm out}} \cup s_{2;{\rm out}})} \, \xi^{s_0}_{\rm out},
$$
where $s_{1;{\rm out}} = s_1|_{\partial_{\rm out} D_1'}$ and $s_{2;{\rm out}} = s_2|_{\partial_{\rm out} D_2'}$. Observe that
\begin{align*}
\mathcal{C}(b_{\rm out}; (\partial_{\rm out} D_1') \cup (\partial_{\rm out} D_2'), s_{1;{\rm out}} \cup s_{2;{\rm out}}) & = \mathcal{C}(b_{\rm out} ; \partial_{\rm out} D_1', s_{1;{\rm out}}) + \mathcal{C}(b_{\rm out}; \partial_{\rm out} D_1', s_{2;{\rm out}}) \\
& \quad + \underbrace{ \sum_{x\in \partial_{\rm out} D_2', \, y\in \partial_{\rm out} D_1'} {\rm sgn}(b_{\rm out}; \overrightarrow{xy}) s_0(x) s_0(y)}.
\end{align*}
Because $\partial_{\rm out}D_2'$ is vertically and horizontally lower than $\partial_{\rm out}D_1'$, and the boundary-orientation on $b_{\rm out}$ is going toward the horizontally increasing direction, we have ${\rm sgn}(b_{\rm out}; \overrightarrow{xy})=1$ for $x\in \partial_{\rm out} D_2'$ and $y\in \partial_{\rm out} D_1'$, and hence the underbraced part equals $(\sum_{x\in \partial_{\rm out}D_2'} s_0(x))(\sum_{y\in \partial_{\rm out}D_1'}s_0(y))$. Combining these, together with the fact
$$
F([D_1] \otimes [D_2])(\xi^{s_0}_{\rm in}) = F([D_1])(\xi^{s_1}_{\rm in}) \otimes F([D_2])(\xi^{s_2}_{\rm in}) = (F([D_1]) \otimes F([D_2]))(\xi^{s_0}_{\rm in}),
$$
one obtains
\begin{align*}
\langle \xi^{s_0}_{\rm out}, \, \wh{G}([D_1]\otimes [D_2])(\xi^{s_0}_{\rm in}) \rangle & = \omega^{-(*)_{s_0}} \, \langle \xi^{s_1}_{\rm out} \otimes \xi^{s_2}_{\rm out}, \, \wh{G}([D_1]) (\xi^{s_1}_{\rm in}) \otimes \wh{G}([D_2])(\xi^{s_2}_{\rm in}) \rangle \\
& = \omega^{-(*)_{s_0}} \, \langle \xi^{s_0}_{\rm out}, \, 
(\wh{G}([D_1])\otimes \wh{G}([D_2]))(\xi^{s_0}_{\rm in}) \rangle,
\end{align*}
where $(*)_{s_0}$ is the sum of the underbraced terms:
$$
(*)_{s_0} = - (\underset{x\in\partial_{\rm in} D_2'}{\textstyle \sum} s_0(x))(\underset{y\in \partial_{\rm in}D_1'}{\textstyle \sum} s_0(y)) + (\underset{x\in\partial_{\rm out} D_2'}{\textstyle \sum} s_0(x))(\underset{y\in \partial_{\rm out}D_1'}{\textstyle \sum} s_0(y)).
$$
Notice that $s_0$ is charge-preserving (as defined in Lem.\ref{lem:simple_properties_of_G_D}) if both $s_1$ and $s_2$ are charge-preserving. In this case, we have $(*)_{s_0}=0$, and therefore
$$
\langle \xi^{s_0}_{\rm out}, \wh{G}([D_1]\otimes[D_2])(\xi^{s_0}_{\rm in})\rangle = \langle \xi^{s_0}_{\rm out}, (\wh{G}([D_1])\otimes \wh{G}([D_2]))(\xi^{s_0}_{\rm in})\rangle
$$
holds. When $s_1$ or $s_2$ is not charge-preserving, both sides of this equality can be easily shown to be zero, using Lem.\ref{lem:charge_conservation_of_F}. Hence we get
$$
\wh{G}([D_1]\otimes[D_2])(\xi^{s_0}_{\rm in}) = (\wh{G}([D_1])\otimes \wh{G}([D_2]))(\xi^{s_0}_{\rm in}).
$$

\vs

Now let $[D_1], [D_2] \in \mathcal{D}(\vec{B})$, and suppose that $[D_1]$ is composable with $[D_2]$ as in Def.\ref{def:tensor_product_and_composition_of_oriented_tangle_diagrams_in_a_directed_biangle}; we shall show that $\wh{G}([D_1] \circ [D_2]) = \wh{G}([D_1]) \circ \wh{G}([D_2])$. If we write $[D] = [D_1] \circ [D_2]$ (as in Def.\ref{def:tensor_product_and_composition_of_oriented_tangle_diagrams_in_a_directed_biangle}), then by construction we see that ${\rm wr}(D) = {\rm wr}(D_1) + {\rm wr}(D_2)$. Note that there is a bijection between $\partial_{\rm in} D_1$ and $\partial_{\rm out} D_2$ that preserves both the horizontal and vertical orderings. Thus we have $\mathcal{C}(b_{\rm in};\partial_{\rm in}D_1) = \mathcal{C}(b_{\rm out};\partial_{\rm out}D_2)^{-1}$, in view of the boundary-orientations on $b_{\rm in}$ and $b_{\rm out}$. Now, note that $\partial_{\rm in} D$ is naturally identified with $\partial_{\rm in} D_2$, as $\partial_{\rm out} D$ is with $\partial_{\rm out}D_1$. Combining all these, together with the fact $F([D_1]\circ [D_2]) = F([D_1]) \circ F([D_2])$, one obtains the sought-for $\wh{G}([D_1]\circ[D_2]) = \wh{G}([D_1]) \circ \wh{G}([D_2])$.

\vs

It remains to check values on elementary classes. For the case of Prop.\ref{prop:F_D}(RT1), it is easy to see that both $\wh{G}([D])$ and $G([D])$ are identity maps; see Prop.\ref{prop:G_D}(GB1). For Prop.\ref{prop:F_D}(RT2), we have ${\rm wr}(D)=0$, $\mathcal{C}(b_{\rm in};\partial_{\rm in}D)={\rm id}$ and $\mathcal{C}(b_{\rm out};\partial_{\rm out}D)(\xi_\pm \otimes \xi_\mp) = 1$. So, we see that $\wh{G}([D])(1) = \omega^{-1}(\omega\xi_+\otimes\xi_- - \omega^5\xi_-\otimes \xi_+)$, which equals $G([D])$, as desired. For Prop.\ref{prop:F_D}(RT3), we have ${\rm wr}(D)=0$, $\mathcal{C}(b_{\rm out};\partial_{\rm out}D)={\rm id}$ and $\mathcal{C}(b_{\rm in};\partial_{\rm in}D)(\xi_\pm \otimes \xi_\mp) = -1$, so $\wh{G}([D])$ sends $\xi_+ \otimes \xi_-$ to $\omega(-\omega^{-5})$ and $\xi_- \otimes \xi_+$ to $\omega\omega^{-1}$ while sending other basis vectors to zero, just as $G([D])$ does, as desired. 

\vs

For  Prop.\ref{prop:F_D}(RT4) as in Fig.\ref{fig:elementary_bo_oriented_tangle_diagrams_in_a_directed_biangle}(IV), we have ${\rm wr}(D)=0$, and for each state $s_0$ of $D$ we have $\mathcal{C}(b_{\rm in};\partial_{\rm in}D)(\xi^{s_0}_{\rm in}) = s_0(x_1)s_0(x_2)$ and $\mathcal{C}(b_{\rm out};\partial_{\rm out}D)(\xi^{s_0}_{\rm out}) = s_0(y_1) s_0(y_2)$, so from Prop.\ref{prop:F_D}(RT4) we get
\begin{align*}
\wh{G}([D]) : \xi_i \otimes \xi_j \mapsto \left\{
\begin{array}{ll}
\omega^{-1} \omega^{-1} (\omega^2 \, \xi_i \otimes \xi_j), & \mbox{if $i= j$,} \\
\omega \omega (\omega^{-2} ( \xi_+ \otimes \xi_- + (\omega^4-\omega^{-4})\xi_-\otimes \xi_+)), & \mbox{if $i=+$ and $j=-$,} \\	
\omega \omega (\omega^{-2} \xi_- \otimes \xi_+), & \mbox{if $i=-$ and $j=+$,}
\end{array}
\right.
\end{align*}
which equals $G([D])$ as in Prop.\ref{prop:G_D}(GB4). For the case as in Fig.\ref{fig:elementary_bo_oriented_tangle_diagrams_in_a_directed_biangle}(XI) we have ${\rm wr}(D)=0$, and for each state $s_0$ of $D$ we have we have $\mathcal{C}(b_{\rm in};\partial_{\rm in}D)(\xi^{s_0}_{\rm in}) = -s_0(x_1)s_0(x_2)$ and $\mathcal{C}(b_{\rm out};\partial_{\rm out}D)(\xi^{s_0}_{\rm out}) = -s_0(y_1) s_0(y_2)$, so $\wh{G}([D])$ coincides with the inverse of the map $\wh{G}([D])$ for Fig.\ref{fig:elementary_bo_oriented_tangle_diagrams_in_a_directed_biangle}(IV), hence also matches $G([D])$ as in Prop.\ref{prop:G_D}(GB4). For Prop.\ref{prop:F_D}(RT5), we have ${\rm wr}(D)=1$. For each state $s_0$ of $D$, note $\mathcal{C}(b_{\rm in};\partial_{\rm in} D)(\xi^{s_0}_{\rm in}) = s_0(x_1) s_0(x_2)$ and $\mathcal{C}(b_{\rm out};\partial_{\rm out} D)(\xi^{s_0}_{\rm out}) = - s_0(y_1) s_0(y_2)$. So, from Prop.\ref{prop:F_D}(RT5) we get
\begin{align*}
\wh{G}([D]) : \xi_i \otimes \xi_j \mapsto \left\{
\begin{array}{ll}
\omega^{-2} \omega \omega^{-1} (\omega^{-2} \xi_i \otimes \xi_j), & \mbox{if $i= j$,} \\
\omega^{-2} \omega^{-1} \omega (\omega^2 \, \xi_- \otimes \xi_+), & \mbox{if $i=+$ and $j=-$,} \\
\omega^{-2} \omega^{-1}\omega (\omega^2( \xi_+ \otimes \xi_- + (\omega^{-4}-\omega^4) \xi_- \otimes \xi_+)), & \mbox{if $i=-$ and $j=+$,}
\end{array}
\right.
\end{align*}
which equals $G([D])$ as in Prop.\ref{prop:G_D}(GB5) as desired. For Prop.\ref{prop:F_D}(RT6), we have ${\rm wr}(D)=-1$. For each state $s_0$ of $D$, note $\mathcal{C}(b_{\rm in};\partial_{\rm in} D)(\xi^{s_0}_{\rm in}) = s_0(x_1) s_0(x_2)$ and $\mathcal{C}(b_{\rm out};\partial_{\rm out} D)(\xi^{s_0}_{\rm out}) = - s_0(y_1) s_0(y_2)$.  So one can observe that the map $\wh{G}([D])$ coincides with the inverse of the map $\wh{G}([D])$ for the case of Prop.\ref{prop:F_D}(RT5), and hence  equals $G([D])$ as in Prop.\ref{prop:G_D}(GB6), as desired. \qed

\vs

The above proof of Lem.\ref{lem:relationship_between_the_operator_invariants_G_and_F} in particular provides a proof of Prop.\ref{prop:G_D}, the well-definedness of the invariant $G([D])$, as promised. Now, putting together Def.\ref{def:Gabella_biangle_quantum_holonomy_as_matrix_elements}, Lem.\ref{lem:relationship_between_the_operator_invariants_G_and_F} and Prop.\ref{prop:Bonahon-Wong_biangle_quantum_trace_as_matrix_elements}, we obtain a precise relationship between the biangle factors of the two sides:
\begin{lemma}[relationship between Gabella's biangle quantum holonomy and Bonahon and Wong's biangle quantum trace]
\label{lem:relationship_between_biangle_factors}
Let $\vec{B}$, $B=(\Sigma',\mathcal{P}')$, $\frak{S}' = \Sigma'\setminus\mathcal{P}'$ and $\omega$ be as in Prop.\ref{prop:Bonahon-Wong_biangle_quantum_trace_as_matrix_elements}. If $[D,s_0] \in \mathcal{D}_{\rm s}(\vec{B})$ is the equivalence class of stated boundary-ordered tangle diagram for the VH-isotopy class $[k,s_0]$ of stated VH-tangles in $\frak{S}'\times(-1,1)$, one has
\begin{align}
\label{eq:relationhip_between_biangle_factors}
{\rm Tr}^\omega_B ([k,s_0]) = \omega^{2{\rm wr}(D)} \, \omega^{\partial \mathcal{C}_B(k,s_0)} \, {\rm TrHol}^\omega_B ([D, s_0])
\end{align}
where $\partial \mathcal{C}_B(k,s_0)$ is the boundary signed order correction amount defined in Def.\ref{def:signed_order_correction_amount}.
\end{lemma}

{\it Proof.} For convenience, write $s=s_0$ for now. Note that
\begin{align*}
& {\rm Tr}^\omega_B([k,s_0]) & \\ 
& = \langle \xi^s_{\rm out}, \, F([D]) \xi^s_{\rm in} \rangle & \mbox{(by {\rm Prop}.\ref{prop:Bonahon-Wong_biangle_quantum_trace_as_matrix_elements})} \\
& = \langle \xi^s_{\rm out}, \, \omega^{2{\rm wr}(D)} \, \mathcal{C}(b_{\rm out};\partial_{\rm out}D) \circ G([D]) \circ \mathcal{C}(b_{\rm in};\partial_{\rm in}D) \xi^s_{\rm in}\rangle & \mbox{(by {\rm Lem}.\ref{lem:relationship_between_the_operator_invariants_G_and_F})} \\
& = \langle \xi^s_{\rm out}, \, \omega^{2{\rm wr}(D)} \, \mathcal{C}(b_{\rm out};\partial_{\rm out}D) \circ G([D]) (\omega^{(*)_1} \xi^s_{\rm in}) \rangle & \mbox{(by {\rm Def}.\ref{def:signed_order_correction_operator})} \\
& & \hspace{-60mm} \mbox{(with $(*)_1 = \mathcal{C}(b_{\rm in}; \partial_{\rm in} D, s_{\rm in})$, where $s_{\rm in} = s|_{\partial_{\rm in}D}$)} \\
& = \omega^{2{\rm wr}(D) + (*)_1} \left\langle \xi^s_{\rm out}, \, \mathcal{C}(b_{\rm out}; \partial_{\rm out} D)({\textstyle \sum}_{\vec{\varepsilon}} \, \langle \xi^{\vec{\varepsilon}}, \, G([D]) \xi^s_{\rm in} \rangle \, \xi^{\vec{\varepsilon}} \right\rangle & \\
& = \omega^{2{\rm wr}(D) + (*)_1} \left\langle \xi^s_{\rm out}, \,({\textstyle \sum}_{\vec{\varepsilon}} \, \langle \xi^{\vec{\varepsilon}}, \, G([D]) \xi^s_{\rm in} \rangle \, \omega^{\mathcal{C}(b_{\rm out}; \partial_{\rm out} D, \vec{\varepsilon})} \xi^{\vec{\varepsilon}} \right\rangle & \mbox{(by {\rm Def}.\ref{def:signed_order_correction_operator})} \\
& = \omega^{2{\rm wr}(D) + \mathcal{C}(b_{\rm in};\partial_{\rm in}D, s_{\rm in})} \langle \xi^s_{\rm out}, \, G([D]) \xi^s_{\rm in} \rangle \, \omega^{\mathcal{C}(b_{\rm out}; \partial_{\rm out} D, s_{\rm out})} & \mbox{(by {\rm Def}.\ref{def:basis_of_tensor_space_and_bilinear_form})} \\
& = \omega^{2{\rm wr}(D)} \, \omega^{\mathcal{C}(b_{\rm in}; \partial_{\rm in}D, s_{\rm in}) + \mathcal{C}(b_{\rm out}; \partial_{\rm out}D, s_{\rm out})} \, {\rm TrHol}^\omega_B ([D,s]), & \mbox{(by {\rm Def}.\ref{def:Gabella_biangle_quantum_holonomy_as_matrix_elements})}
\end{align*}
where the seventh line is with respect to the horizontal ordering on $\partial_{\rm out}D$, and $s_{\rm out} = s|_{\partial_{\rm out}D}$ in the eighth line. \qed

\vs

One consequence of the above lemma is a proof of Lem.\ref{eq:G_D_useful_properties}. Before proving Lem.\ref{eq:G_D_useful_properties}(1), the direction independence of Gabella biangle quantum holonomy, the ${\rm TrHol}^\omega_B$ in the right hand side of eq.\eqref{eq:relationhip_between_biangle_factors} should really be written as ${\rm TrHol}^\omega_{\vec{B}}$. However, since all other factors appearing in eq.\eqref{eq:relationhip_between_biangle_factors} are independent of choice of direction on $B$, it follows that so is ${\rm TrHol}^\omega_{\vec{B}}([D,s_0])$, as asserted in Lem.\ref{eq:G_D_useful_properties}(1). For Lem.\ref{eq:G_D_useful_properties}(2), the charge conservation property of Gabella biangle quantum holonomy ${\rm TrHol}^\omega_B$, first recall Lem.\ref{lem:charge_conservation_of_F} which is the charge conservation property of the Bonahon-Wong biangle quantum trace ${\rm Tr}^\omega_B$. Note that the remaining factors in eq.\eqref{eq:relationhip_between_biangle_factors} are powers of $\omega$, hence nonzero. The sought-for charge conservation property of ${\rm TrHol}^\omega_B$ follows.

\subsection{Finishing the proof of main theorem}

We have only to assemble everything established so far to prove the sought-for term-by-term equality eq.\eqref{eq:term_equality}. We first exclude the easy case. For each edge $e$ of $\wh{\Delta}$, let $b_e(J)$ to be the net sum of signs assigned by $J$ to the junctures lying over $e$. If there is a biangle $B_i$ (whose sides are $e_i$ and $e_i'$) such that $b_{e_i}(J) \neq b_{e_i'}(J)$, then by the charge conservation properties (Lem.\ref{lem:simple_properties_of_G_D}(2) and Lem.\ref{lem:charge_conservation_of_F}) it follows that the biangle factors ${\rm Tr}^\omega_{B_i}([L_i, J|_{\partial L_i}])$ and ${\rm TrHol}^\omega_{B_i}([L_i,J|_{\partial L_i}])$ are zero, hence both sides of eq.\eqref{eq:term_equality} are zero. From now on, we may assume
\begin{align}
\label{eq:balanced_condition}
b_{e_i}(J) = b_{e_i'}(J) \quad\mbox{for all $i=1,\ldots,n$.}
\end{align}
The triangle part of the left hand side ${\rm BW}^\omega_{\wh{\Delta}}(K';J) = {\rm BW}^\omega_{\wh{\Delta}}(K;J)$ of eq.\eqref{eq:term_equality} is 
\begin{align*}
\hspace{-5mm} \bigotimes_{j=1}^m {\rm Tr}^\omega_{\wh{t}_j}([K_j, J|_{\partial K_j}]) & = \omega^{\sum_{j=1}^m {\rm dev}_{\wh{t}_j}(K_j,J|_{\partial K_j}) } \, {\textstyle \prod}_{j=1}^m \left[ {\rm Tr}^\omega_{\wh{t}_j}( [k_{j,1}, J|_{\partial k_{j,1}}] ) \cdots {\rm Tr}^\omega_{\wh{t}_j}( [k_{j,l_j}, J|_{\partial k_{j,l_j}}] ) \right]_{\rm Weyl} \quad \begin{array}{ll} \mbox{($\because$ eq.\eqref{eq:triangle-factor} } \\ \mbox{and Def.\ref{def:deviation_of_BW_traingle_factor})} \end{array} \\
& = \omega^{\sum_{j=1}^m {\rm dev}_{\wh{t}_j}(K_j,J|_{\partial K_j}) } \, \left[ {\textstyle \prod}_{j=1}^m  {\rm Tr}^\omega_{\wh{t}_j}( [k_{j,1}, J|_{\partial k_{j,1}}] ) \cdots {\rm Tr}^\omega_{\wh{t}_j}( [k_{j,l_j}, J|_{\partial k_{j,l_j}}] ) \right]_{\rm Weyl}\\
& = \omega^{\sum_{j=1}^m {\rm dev}_{\wh{t}_j}(K_j,J|_{\partial K_j}) } \left[{\textstyle \prod}_{i=1}^n \wh{Z}_{e_i}^{b_{e_i}(J)} \right]_{\rm Weyl} \quad \begin{array}{cc} \mbox{($\because$ Prop.\ref{prop:BW_full}(2)(a), Def.\ref{def:Gabella_quantum_holonomy}(G1),} \\ \mbox{eq.\eqref{eq:balanced_condition}, Def.\ref{def:square-root_Chekhov-Fock_algebra})} \end{array} \\
& = \omega^{\sum_{j=1}^m {\rm dev}_{\wh{t}_j}(K_j,J|_{\partial K_j}) } \, \wh{Z}_{\til{K}^J} \qquad (\because \mbox{Def.\ref{def:Gabella_quantum_holonomy}(G1)}) \\
& = \omega^{-4\sum_{j=1}^m {\rm wr}_{\wh{t}_j}(\til{K}^J_j)}  \omega^{2\sum_{j=1}^m {\rm wr}(K_j)} \omega^{\sum_{j=1}^m \partial \mathcal{C}_{\wh{t}_j}(K_j,J|_{\partial K_j})} \, \wh{Z}_{\til{K}^J} \qquad (\because \mbox{Lem.\ref{lem:equality_of_the_triangle_factors}})
\end{align*}
Note $\omega^{-4\sum_{j=1}^m {\rm wr}_{\wh{t}_j}(\til{K}^J_j)} = \omega^{-4 {\rm wr}_{\wh{\Delta}}(\til{K}^J)} = q^{-{\rm wr}_{\wh{\Delta}}(\til{K}^J)}$, in view of Def.\ref{def:Gabella_quantum_holonomy}(G2). From Lem.\ref{lem:relationship_between_biangle_factors} we see that the biangle part of ${\rm BW}^\omega_{\wh{\Delta}}(K;J)$ is as follows, if we let $D_i$ be the tangle diagram in $B_i$ for the tangle $L_i$ in the thickening of $B_i$:
\begin{align*}
{\textstyle \prod}_{i=1}^n {\rm Tr}^\omega_{B_i} ([L_i, J|_{\partial L_i}]) & = \omega^{2\sum_{i=1}^n {\rm wr}(D_i)} \omega^{\sum_{i=1}^n \partial \mathcal{C}_{B_i}(L_i, J|_{\partial L_i})} {\textstyle \prod}_{i=1}^n {\rm TrHol}^\omega_{B_i} ([L_i, J|_{\partial L_i}]).
\end{align*}
Now observe $\sum_{j=1}^m {\rm wr}(K_j) + \sum_{i=1}^n {\rm wr}(D_i) = {\rm wr}(K)$, and also:
\begin{lemma}
$
\sum_{j=1}^m \partial \mathcal{C}_{\wh{t}_j} (K_j, J|_{\partial K_j}) + \sum_{i=1}^n \partial \mathcal{C}_B(L_i, J|_{\partial L_i}) = \partial \mathcal{C}_{(\Sigma,\mathcal{P})}(K,s)
$.
\end{lemma}
{\it Proof.} Each edge of $\wh{\Delta}$ that is not a boundary arc of $(\Sigma,\mathcal{P})$ appear as a side of exactly one triangle $\wh{t}_j$ and one biangle $B_i$, with opposite orientations, hence the corresponding terms cancel each other. What remain are the boundary signed order correction amounts for the boundary arcs of $(\Sigma,\mathcal{P})$, which form exactly $\partial \mathcal{C}_{(\Sigma,\mathcal{P})}(K,s)$. \qed

\vs

Combining everything:
\begin{align*}
{\rm BW}^\omega_{\wh{\Delta}}(K;J) & = \left({\textstyle \prod}_{i=1}^n {\rm Tr}^\omega_{B_i}([L_i, J|_{\partial L_i}])\right) \left( {\textstyle \bigotimes}_{j=1}^m {\rm Tr}^\omega_{\wh{t}_j}([K_j, J|_{\partial K_j}]) \right) \qquad (\because \mbox{eq.\eqref{eq:BW_term}}) \\
& = \omega^{2{\rm wr}(K)} \omega^{\partial \mathcal{C}_{(\Sigma,\mathcal{P})}(K,s)} \left( {\textstyle \prod}_{i=1}^n {\rm TrHol}^\omega_{B_i} ([L_i, J|_{\partial L_i}]) \right) q^{-{\rm wr}_{\wh{\Delta}}(\til{K}^J)} \, \wh{Z}_{\til{K}^J} \\
& = \omega^{2{\rm wr}(K)} \, \omega^{\partial \mathcal{C}_{(\Sigma,\mathcal{P})}(K,s)} \, \ol{\ul{\Omega}}(\til{K}^J;\omega) \, \wh{Z}_{\til{K}^{\hspace{-0,3mm}J}}. \quad (\because \mbox{eq.\eqref{eq:Gabella_coefficient}})
\end{align*}
This proves the term-by-term equality eq.\eqref{eq:term_equality}, hence the main theorem of the paper.

\vspace{1mm}

\end{document}